\documentclass[a4paper]{article}
\usepackage[top=25truemm,bottom=25truemm,left=20truemm,right=20truemm]{geometry}
\usepackage{amsmath,amsthm,amsfonts,amssymb,stmaryrd,indentfirst,tabularray,makecell,hhline,comment,mathtools}
\usepackage[mathscr]{euscript}
\usepackage[shortlabels]{enumitem}
\usepackage[labelfont=rm]{caption}
\usepackage{titlesec}
\usepackage{tikz-cd}
\usepackage{url}
\usepackage[backref]{hyperref}
\hypersetup{
    colorlinks = true,
    linkcolor = [rgb]{0.3,0.55,0.25},
    citecolor = [rgb]{0.164,0.39,0.77},
    urlcolor = [rgb]{0.7,0.2,0.4}
}

\DeclareMathOperator{\id}{id}
\DeclareMathOperator{\Ker}{Ker}
\DeclareMathOperator{\Imag}{Im}
\DeclareMathOperator{\Der}{Der}

\DeclareMathOperator{\tDer}{tDer}
\DeclareMathOperator{\Hom}{Hom}

\DeclareMathOperator{\Aut}{Aut}
\DeclareMathOperator{\tAut}{tAut}

\DeclareMathOperator{\ad}{ad}
\DeclareMathOperator{\Exp}{exp}
\DeclareMathOperator{\rot}{rot}

\newcommand{\gr}[0]{_\mathrm{gr}}
\newcommand{\sdiv}[0]{\mathsf{div}}
\newcommand{\K}[1]{\mathbb{K}[[#1]]}

\theoremstyle{plain}
\newtheorem{theorem}{Theorem}
\newtheorem*{theorem*}{Theorem}
\newtheorem*{maintheorem*}{Main Theorem}
\newtheorem{proposition}[theorem]{Proposition}
\newtheorem{deflem}[theorem]{Definition-Lemma}
\newtheorem{lemma}[theorem]{Lemma}

\theoremstyle{definition}
\newtheorem{definition}[theorem]{Definition}
\newtheorem{remark}[theorem]{Remark}

\newtheorem{problem}[theorem]{Problem}

\numberwithin{theorem}{section}
\everymath{\displaystyle}
\allowdisplaybreaks

\usepackage{apptools}
\titleformat{\section}{\large\scshape}{\IfAppendix{\appendixname}{} \thesection.}{3pt}{}
\titleformat{\subsection}{\scshape}{\thesubsection.}{3pt}{}

\usepackage{titling}  
\pretitle{
	\begin{center}\Large\bfseries
}
\posttitle{
	\end{center}\ \\[-20pt]
}
\preauthor{
	\begin{center}\large
}
\postauthor{
	\end{center}\ \\[-60pt]
}

\title{Drinfeld Associators and Kashiwara--Vergne Associators\\in Higher Genera}
\author{Toyo TANIGUCHI \thanks{Graduate School of Mathematical Sciences, The University of Tokyo. 3-8-1, Komaba, Meguro-ku, Tokyo, 153-8914, Japan. E-mail: \texttt{toyo(at)ms.u-tokyo.ac.jp}}}
\date{}

\begin{document}
\maketitle

\begin{abstract}
\noindent For $g\geq 0$, a genus $g$ Kashiwara--Vergne associator, introduced by Alekseev--Kawazumi--Kuno--Naef as a solution to the generalised KV equations in relation to the formality problem of the Goldman--Turaev Lie bialgebra on an oriented surface with a framing, is directly constructed from a genus $g$ analogue of a Drinfeld associator formulated by Gonzalez, which we call a Gonzalez--Drinfeld associator. The proof is based on Massuyeau's work in genus $0$. The framing is determined from the choice of a Gonzalez--Drinfeld associator, and in the case of genus $1$, we show that only particular framings are realised by our construction.\\
\end{abstract}

\noindent{\textit{2020 Mathematics Subject Classification: 16T05, 16W70, 18M75, 20F36, 20F40, 57K20, 57M05.}

\section{Introduction}
A \textit{Drinfeld associator} was first defined in his original paper \cite{drinfeld} as an associator for the category of representations of a quasi-Hopf algebra. Over a fixed field $\mathbb{K}$ of characteristic zero, a Drinfeld associator is expressed as the exponential of a Lie series in two variables satisfying a certain system of equations, and one example (for $\mathbb{K} = \mathbb{C}$) is obtained from a solution to the Knizhnik--Zamolodchikov equation. Apart from the application to quasi-tensored categories, a Drinfeld associator appears in low-dimensional topology; most notably, it is used in the construction of the Kontsevich knot invariant (see \cite{kontint}, \cite{cartier} or \cite{kastur}, for example).

On the other hand, there is a very closely related object in Lie theory: a solution to the Kashiwara--Vergne (KV) equations in \cite{kv}, which we will call a \textit{Kashiwara--Vergne associator}. The relation to Drinfeld associators is given by Alekseev--Torossian in \cite{at}, where they constructed a solution to the (generalised) KV equations from a Drinfeld associator. The explicit form of their map was later described in the paper \cite{aet} by Alekseev--Enriquez--Torossian.

KV associators have several low-dimensional interpretations, such as \cite{weldedfoams} via welded foams, but we only deal with the relation to the Goldman--Turaev Lie bialgebra. In the paper \cite{akkn1} by Alekseev--Kawazumi--Kuno--Neaf, it is shown that a solution to the \textit{formality problem} of the Goldman--Turaev Lie bialgebra on the pair of pants $\Sigma_{0,3}$ is almost equivalent to a KV associator; namely, these is a map from the set of solutions to the KV equations to the set of formality morphisms for $\Sigma_{0,3}$, which is surjective up to inner automorphisms of the completed group algebra (see the end of Section \ref{sec:kv} for details). At that time, Massuyeau had already constructed a formality morphism directly from a Drinfeld associator in \cite{mas}, realising the Alekseev--Torossian map in terms of the Goldman--Turaev Lie bialgebra. The commutative diagram below summarises the relations of these works.\\[-10pt]
\[\begin{tikzcd}[cramped]
	\{\textrm{Drinfeld associators}\} \arrow[r, "\substack{\scriptsize\cite{at}\\\scriptsize\cite{aet}}"]\arrow[rd, "\scriptsize\cite{mas}"'] &\{\textrm{Solutions to the generalised KV equations}\}=: \mathrm{SolKV}^\mathsf{fr}_{0,3}\arrow[d,"\scriptsize\cite{akkn1}"]\\
	& \left\{\begin{gathered}\textrm{Formality morphisms of the Goldman--Turaev}\\\textrm{Lie bialgebra on the pair of pants}\end{gathered}\right\}=: \mathrm{Form}^\mathsf{fr}_{0,3}
\end{tikzcd}\]
Furthermore, the genus $g$ version of the KV equations was introduced in their sequel \cite{akkn} and it was shown that there is an analogous map $\mathrm{SolKV}_{g,n+1} \to \mathrm{Form}_{g,n+1}$ with $\mathrm{SolKV}_{g,n+1}$ non-empty. Based on the result of \cite{akkn1}, we call a solution to the genus $g$ KV equation a \textit{genus} $g$ \textit{KV associator}. Another proof of the existence of the formality morphism for any genus $g$ was also obtained by Hain \cite{hain} via the theory of mixed Hodge structure.

Back to Drinfeld associators, there are higher genus analogues too: for genus $1$, one version was introduced by Enriquez \cite{enriquez} and is called \textit{elliptic associators}. For an arbitrary genus, several generalisations are proposed: \cite{gonzalez} by Gonzalez, \cite{felder} by Felder, and \cite{ciw} by Campos--Idrissi--Willwacher. Their relation is explained in  \cite{gonzalez} for the genus $1$ case, and one example of an elliptic associator is obtained from the universal KZB equation by the work of Calaque--Enriquez--Etingof \cite{cee}. For higher genera, however, the relation between them (especially the existence and whether they agree) is still an open question. It is expected to be given by a solution to the higher genus version of the KZB equation; see Conjecture 3.22 of \cite{gonzalez}.\\

In this paper, we adopt the definition by Gonzalez and we call them \textit{genus} $g$ \textit{Gonzalez--Drinfeld associators}. First, we construct a map analogous to the diagonal arrow in the above diagram for any $g\geq 0$. Let $\mathbf{Ass}'_g$ be the set of genus $g$ Gonzalez--Drinfeld associators such that the coupling constant is $1$ and the associated graded map is the identity map, and $\mathrm{Form}_{g,n+1}$ the set of solutions to the formality problem for the Goldman--Turaev Lie bialgebra on the surface $\Sigma_{g,n+1}$ of genus $g$ with $n+1$ boundary components with any framing.\begin{theorem*}[Theorem \ref{thm:assoc}]
For $g,n\geq 0$ and $f \in\mathbb{K}^n$, we have a map $I_{g,n+1}^f \colon \mathbf{Ass}'_g \to \mathrm{Form}_{g,n+1}$.
\end{theorem*}
\noindent The vector $f$ corresponds to the choice of rotation numbers along the boundary components. Our proof is based on Massuyeau's method \cite{mas} using his three-dimensional formulae for loop operations.

Next, we construct a horizontal arrow by combining the result in \cite{akkn} with the map $I_{g,n+1}^f$. Let $\mathrm{SolKV}_{g,n+1}$ be the set of solutions to the KV problem of type $(g, n+1)$ in the sense of \cite{akkn}, recalled in Section \ref{sec:kv}.

\begin{theorem*}[Theorem \ref{thm:solkv}]
For $g,n\geq 0$ and $f \in\mathbb{K}^n$, we have a map $\tilde{I}_{g,n+1}^f\colon \mathbf{Ass}'_g \to \mathrm{SolKV}_{g,n+1}$, which is a lift of $I_{g,n+1}^f$.
\end{theorem*}

These maps are summarised into the following commutative diagram.
\[\begin{tikzcd}[cramped]
	\mathbf{Ass}'_g \arrow[r, "\tilde I_{g,n+1}^f"]\arrow[rd, "I_{g,n+1}^f"'] & \mathrm{SolKV}_{g,n+1}\arrow[d,"\scriptsize\cite{akkn}"]\\
	& \mathrm{Form}_{g,n+1}
\end{tikzcd}\]

The set of genus $g$ Gonzalez--Drinfeld associators is a torsor (if it is non-empty) over the \textit{genus} $g$ \textit{Grothendieck--Teichm\"uller group}. The subgroup $\widehat{\mathbf{GT}}{}'_g$ is defined such that it acts on $\mathbf{Ass}'_g$, and indeed we can show the following:
\begin{theorem*}[Theorem \ref{thm:equivariant}]
There is an action of $\widehat{\mathbf{GT}}{}'_g$ on $\mathrm{SolKV}_{g,n+1}$ such that the map $\tilde{I}_{g,n+1}^f$ is $\widehat{\mathbf{GT}}{}'_g$-equivariant.
\end{theorem*}

On another note, Neaf \cite{neaf} independently gave the essentially same construction. He also adopts the same formulation for both Gonzalez--Drinfeld- and Kashiwara--Vergne associators, but his proof is based on the cohomological description of the Goldman--Turaev Lie bialgebra, avoiding explicit calculations.\\

In the construction of the map $I_{g,n+1}^f$, a framing $\mathsf{fr}_{\vec Z, f}$ of the surface $\Sigma = \Sigma_{g,n+1}$ is defined from a genus $g$ Gonzalez--Drinfeld associator $\vec Z$ and $f\in\mathbb{K}^n$. In the special case of genus $1$, this framing does not realise every possible framing (as indicated in \cite{akkn}), and we have the following result:
\begin{theorem*}[Theorem \ref{thm:framing}]
A framing on $\Sigma_{1,n+1}$ associated with a genus $1$ Gonzalez--Drinfeld associator restricts to a constant vector field outside $D$, where $D$ is an embedded disk in $\Sigma_{1,0}$ which contains all the punctures of  $\Sigma_{1,n+1}$.
\end{theorem*}

In genus zero, other recent constructions regarding Drinfeld associators and KV associators can be seen in \cite{dancso1} and \cite{dancso2}.\\

\noindent\textbf{Organisation of the paper.} Sections \ref{sec:drinfeld} and \ref{sec:kv} are brief reviews of genus $g$ Gonzalez--Drinfeld- and Kashiwara--Vergne associators with related materials on loop operations. The map $I_{g,n+1}^f$ is constructed in Section \ref{sec:main}, and the well-definedness is shown in Section \ref{sec:threedim}. The map $\tilde I_{g,n+1}^f$ is constructed in Section \ref{sec:kveq}. Section \ref{sec:framing} is devoted to the calculation regarding the framing associated with a genus $g$ Gonzalez--Drinfeld associator. Section \ref{sec:proof} is occupied with the proof of Lemma \ref{lem:exactness}.\\

\noindent\textbf{Acknowledgements.} The author thanks Gw\'ena\"el Massuyeau for answering the author's question on higher genus associators and his own construction involving Drinfeld associators, Nariya Kawazumi for the comment on the inevitability of $\mathbb{K}$-framings, and Florian Neaf for the fruitful discussion on his similar result. This work was supported by JSPS KAKENHI Grant Number 25KJ0734.\\

\noindent\textbf{Conventions.} $\mathbb{K}$ is a field of characteristic zero. Conventions regarding string diagrams are aligned with those of \cite{gonzalez}: they are read from top to bottom, and the corresponding composition of elements is read from left to right. All the figures are equipped with the blackboard framing.\\

\section{Genus $g$ Gonzalez--Drinfeld Associators}\label{sec:drinfeld}

Bar-Natan \cite{barnatan} reformulated Drinfeld associators and the Grothendieck--Teichm\"uller (GT) groups in terms of operads of (infinitesimal) braids. For higher genera, we adopt the version introduced by Gonzalez \cite{gonzalez}, which is an extension of Bar-Natan's definition, and we would like to call it a \textit{genus} $g$ \textit{Gonzalez--Drinfeld associator}. For $g\geq 0$ and a fixed field $\mathbb{K}$ of characteristic zero, it is an isomorphism of operads
\[
	\widehat{\mathbb{K}\mathbf{PaB}}{}^f \to \mathbf{PaCD}^f
\]
together with a compatible isomorphism of the \textit{operad modules}
\[
	\widehat{\mathbb{K}\mathbf{PaB}}{}^f_g \to \mathbf{PaCD}^f_g\,.
\]
From now on, we briefly recall these objects.\\

Following the formulation by Bellingeri and Gervais \cite{bg}, a \textit{framed braid} with $m\geq 1$ strands on a smooth surface $\Sigma$ (whether it is compact, oriented or not) is an element of $\mathrm{B}^f_{m}(\Sigma):= \pi_1(F_m(\Sigma)/S_m)$. The space $F_m(\Sigma)$ is defined as the pull-back bundle fitting in the diagram
\[\begin{tikzcd}[cramped]
	F_m(\Sigma) \arrow[r]\arrow[d] &(U\Sigma)^{\times m}\arrow[d]\\
	\mathrm{Conf}_m(\Sigma) \arrow[r] & \Sigma^{\times m}
\end{tikzcd}\]
where $U\Sigma = (T\Sigma \setminus 0_\Sigma)/\mathbb{R}_{>0}$ is the unit tangent bundle and
\[
	\mathrm{Conf}_m(\Sigma) = \{(x_1, \dotsc, x_m) \in \Sigma^{\times m}: x_i\neq x_j \mbox{ for } i\neq j\}
\]
is the configuration space of $m$ points on $\Sigma$. The symmetric group $S_m$ acts on $F_m(\Sigma)$ by permutation of the points. Similarly, a \textit{framed pure braid} is an element of $\mathrm{PB}^f_{m}(\Sigma) := \pi_1(F_m(\Sigma))$. These fundamental groups are called \textit{framed (pure) braid groups}, and when the surface $\Sigma$ is the connected oriented closed surface $\Sigma_g$ of genus $g\geq 0$, we denote these groups by $\mathrm{B}^f_{g, m}$ and $\mathrm{PB}^f_{g, m}$, respectively. The usual (non-framed) braid groups are defined by $\mathrm{B}_m(\Sigma) = \pi_1(\mathrm{Conf}_m(\Sigma)/S_m)$ and $\mathrm{PB}_m(\Sigma) = \pi_1(\mathrm{Conf}_m(\Sigma))$.

Next, we briefly recall the operad $\mathbf{PaB}^f$ of framed pure braids. The reader is assumed to be familiar with the language of operads. The material below is taken from \cite{gonzalez}. For a rigorous treatment of the non-framed version, see the very detailed exposition \cite{fresse} by Fresse. First of all, the operad $\mathbf{Pa}$ in the category of sets evaluated at a finite set $I$ comprises parenthesised permutations of $I$. Then, $\mathbf{PaB}^f$ is an operad in the category of groupoids having $\mathbf{Pa}$ as an operad of objects and the morphisms between parenthesised permutations $p$ and $q$ is the set of all framed braids on the unit disk in $\mathbb{C}$ with a fixed set of $\# I$ base points, identified with the set $I$, such that the underlying permutation of the braid connects $p$ and $q$. The operadic composition maps are given by substituting a braid into one strand with a certain rotation specified by the framing.

Notable morphisms in $\mathbf{PaB}^f$ are $F^{1,2}\colon (1\,2) \to (1\,2)$ and $R^{1,2}\colon (1\,2) \to (2\,1)$ depicted below.\\[3pt]
\centerline{\raisebox{42pt}{$F^{1,2} = \;\,$}\includegraphics{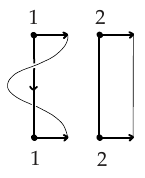}\raisebox{42pt}{,\qquad $R^{1,2} = $}\includegraphics{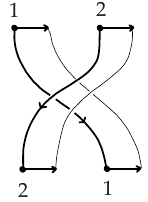}\raisebox{42pt}{.}}
$F^{1,2}$ is the ``negative monogon'' around the point labelled by $1$. An example of operadic composition is shown below.\\
\centerline{\raisebox{84pt}{$(F^{1,2}\circ_2 \varnothing) \circ_1 \id_{(12)} = \;\,$}\includegraphics{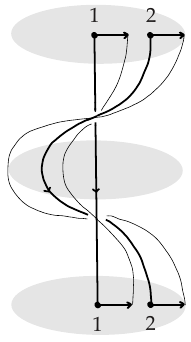}\raisebox{84pt}{$= R^{1,2}R^{2,1}F^{1,2} R^{1,2}F^{2,1}(R^{1,2})^{-1}$}}

The $\mathbb{K}$-linearisation $\mathbb{K}\mathbf{PaB}{}^f$ is an operad in the category of \textit{Hopf groupoids} (the terminology is due to Fresse \cite{fresse}), a multi-object analogue of a Hopf algebra. The completion with respect to the augmentation ideal is denoted by $\widehat{\mathbb{K}\mathbf{PaB}}{}^f$, which lives in the category of completed Hopf groupoids. 

We have another operad $\mathbf{PaCD}^f$ in the same category, which turns out to be the (completed) associated graded quotient of $\widehat{\mathbb{K}\mathbf{PaB}}{}^f$. It has the same operad of objects, and the set of morphisms is explicitly described in terms of the framed Drinfeld--Kohno Lie algebra $\mathfrak{t}^f_I$ for a finite set $I$: the Lie algebra  $\mathfrak{t}^f_I$ over $\mathbb{K}$ is generated by $t_{ij}$ for $i,j\in I$ together with relations
\begin{align*}
\begin{aligned}
	&t_{ij} = t_{ji},&\\
	&[t_{ij}, t_{kl}] = 0& &\mbox{if } \{i,j\}\cap\{k,l\}=\varnothing,\mbox{ and}\\
	&[t_{ij},t_{ik} + t_{jk}] = 0& &\mbox{if } \{i,j\}\cap\{k\}=\varnothing.
\end{aligned}
\end{align*}
Using this Lie algebra, a morphism from any parenthesised permutations $p$ to $q$ is canonically identified with an element of the completed universal enveloping algebra $U(\mathfrak{t}^f_I)$. The operadic composition map is the additive counterpart of that in $\widehat{\mathbb{K}\mathbf{PaB}}{}^f$: it is given by, for finite sets $I$, $J$ and $k\in I$,
\begin{align*}
	\circ_k\colon \mathfrak{t}_I^f \oplus \mathfrak{t}_J^f &\to \mathfrak{t}_{I\sqcup J - \{k\}}^f\\
	(0, t_{\alpha\beta}) &\mapsto t_{\alpha\beta}\\
	(t_{ij}, 0) &\mapsto \left\{\begin{aligned}
		&t_{ij}& &\mbox{if } k \neq i,j,&\\
		&\sum_{l\in J} t_{lj}& &\mbox{if } k=i\neq j,&\\
		&\sum_{l\in J} t_{il}& &\mbox{if } k=j\neq i,&\\
		&\sum_{l, m\in J} t_{lm}& &\mbox{if } k=i=j,&
	\end{aligned}\right.
\end{align*}
and it is naturally extended to their universal enveloping algebras. The mnemonic is that if $k$ appears in the subscript of $t_{ij}$, substitute it with the set $J$ and expand additively over the set $J$ for each subscript. We remark that the formula below Remark 2.7 in \cite{gonzalez} is incorrect: it treats the case $k=i=j$ in the wrong way that the composition would not be associative.

The non-framed versions $\mathfrak{t}_I$ and $\mathfrak{t}_{g, I}$ are obtained as the quotients by modding out the central element $t_{ii}$ for every $i\in I$.\\

\begin{figure}[bt]
\centerline{\includegraphics[scale=0.9]{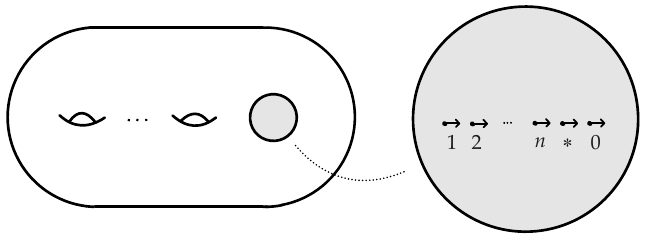}}
\caption{An embedded disk and $n+2$ points with tangent vectors on the surface $\Sigma_g$.}
\label{fig:points}
\end{figure}

On the other hand, the operad module $\mathbf{PaB}^f_g$ over $\mathbf{PaB}^f$ is comprised of braids on the closed surface $\Sigma_g$ of genus $g$. More specifically, fix an embedded disk on $\Sigma_g$ inside which the base points $\mathbf{Pa}$ lie. Then, a ``morphism'' in $\mathbf{PaB}^f_g$ between parenthesised permutations $p$ and $q$ is an element of $\mathrm{B}^f_{g,m}$ where $m$ is the length of $p$ (and thus $q$). The action of $\mathbf{PaB}^f$ is given by the operadic composition of framed braids. The $\mathbb{K}$-linearisation is also filtered by the augmentation ideal, so we can define $\widehat{\mathbb{K}\mathbf{PaB}}{}^f_g$ in a similar manner.

Analogously, we have the operad $\mathbf{PaCD}^f_g$ in which a morphism from $p$ to $q$ in $\mathbf{PaCD}^f_g$ is canoically identified with an element of $U(\mathfrak{t}^f_{g,I})$, where $\mathfrak{t}^f_{g,I}$ is the genus $g$ analogue of $\mathfrak{t}^f_I$. We point out that the definition of $\mathfrak{t}^f_{g, I}$ in Section 3.4.1 of \cite{gonzalez} is incorrect as the original definition does not make the family $\{\mathfrak{t}^f_{g,n}\}_{n\geq 0}$ into an operad module in the category of Lie algebras. The following is the correct one:
\begin{definition}
For $g\geq 0$ and a finite set $I$, the Lie algebra $\mathfrak{t}^f_{g,I}$ is generated by elements
\[
	t_{ij}\; (i,j \in I)\quad\mbox{and}\quad x_i^a, y_i^a\; (i \in I, 1\leq a\leq g)
\]
with the relations given, for $i,j,k,l\in I$ and $1\leq a,b\leq g$, by
\begin{align*}
\begin{aligned}
	&t_{ij} = t_{ji},&\\
	&[t_{ij}, t_{kl}] = 0& &\mbox{if } \{i,j\}\cap\{k,l\}=\varnothing,\\
	&[t_{ij},t_{ik} + t_{jk}] = 0& &\mbox{if } \{i,j\}\cap\{k\}=\varnothing,\\
	&[x_i^a, y_j^b] = \delta_{ab} t_{ij}& &\mbox{if } i\neq j,\\
	&[x_i^a, x_j^b] =[y_i^a, y_j^b] = 0& &\mbox{if } i\neq j,\\
	&[x_k^a, t_{ij}] = [y_k^a, t_{ij}] = 0& &\mbox{if } \{i,j\}\cap\{k\}=\varnothing,\\
	&[x_i^a + x_j^a, t_{ij}] = [y_i^a + y_j^a, t_{ij}] = 0
\end{aligned}
\end{align*}
and, for $i\in I$, 
\[
	\sum_{1\leq a\leq g} [x_i^a, y_i^a] + \sum_{j\in I\setminus\{i\}} t_{ij} = (g-1) t_{ii}\,.
\]
\end{definition}

Only the last relation differs from the original one in \cite{gonzalez}. The module structure over $\mathbf{PaCD}^f$ is given by the natural morphism $\mathfrak{t}_I^f \to \mathfrak{t}_{g,I}^f$ of Lie algebras together with a similar map to $\circ_k$ above.\\

Finally, the set of genus $g$ Gonzalez--Drinfeld associators is defined as the isomorphism set
\[
	\mathbf{Ass}_g = \mathrm{Isom}^+_{\mathrm{Oprd}(\widehat{\mathbf{HGrpd}})}\big((\widehat{\mathbb{K}\mathbf{PaB}}{}^f, \widehat{\mathbb{K}\mathbf{PaB}}{}^f_g), (\mathbf{PaCD}^f, \mathbf{PaCD}^f_g)\big)
\]
of operad modules over the category of complete Hopf groupoids, where the superscript ${}^+$ indicates that we only consider isomorphisms that preserves objects. Since an isomorphism $\widehat{\mathbb{K}\mathbf{PaB}}{}^f \to \mathbf{PaCD}^f$ is equivalent to a Drinfeld associator by the result of Bar-Natan \cite{barnatan} together with Proposition 2.12 of \cite{gonzalez}, a genus $g$ Gonzalez--Drinfeld associator contains a Drinfeld associator as a part of its data. The \textit{coupling constant} $\mu\in \mathbb{K}$ of a Drinfeld associator $Z$ is defined by the formula $Z(R^{1,2}) = \exp(\mu t_{12}/2)$.

Similarly, the Grothendieck--Teichm\"uller groups $\widehat{\mathbf{GT}}{}_g$ and $\mathbf{GRT}_g$ are defined as the automorphism groups of the pairs $(\widehat{\mathbb{K}\mathbf{PaB}}{}^f, \widehat{\mathbb{K}\mathbf{PaB}}{}^f_g)$ and $(\mathbf{PaCD}^f, \mathbf{PaCD}^f_g)$, respectively. By definition, $\mathbf{Ass}_g$ is a bi-torsor over $\widehat{\mathbf{GT}}{}_g$ and $\mathbf{GRT}_g$. The latter has two notable subgroups: the symplectic group $\mathrm{Sp}(2g;\mathbb{K})$ acts on $\mathbf{PaCD}^f_g$, for $\sigma \in \mathrm{Sp}(2g;\mathbb{K})$, by
\begin{align*}
	\sigma(x^a_i) = \sigma(x^a)_i, \quad \sigma(y^a_i) = \sigma(y^a)_i \;\mbox{ and }\; \sigma(t_{ij}) = t_{ij}
\end{align*}
where $(\cdot)_i\colon \mathrm{Span}_\mathbb{K}\{x^a, y^a\} \to \mathfrak{t}_{g,I}^f$ is the linear map specified by $x^a \mapsto x^a_i$ and $y^a \mapsto y^a_i$. The rescaling automorphism $\lambda\in \mathbb{K}^\times$ is defined by
\begin{align*}
	x_i^a \mapsto \lambda x_i^a,\quad y_i^a \mapsto \lambda y_i^a\;\mbox{ and }\; t_{ij} \mapsto \lambda^2 t_{ij}\,.\\
\end{align*}

\section{Loop Operations and Genus $g$ Kashiwara--Vergne Associators} \label{sec:kv}

A genus $g$ Kashiwara--Vergne associator is defined in \cite{akkn} as an automorphism of the (completed) Goldman--Turaev Lie bialgebra. We briefly recall these objects in this section.

Let $g,n\geq 0$ and $\Sigma = \Sigma_{g,n+1}$ a connected compact oriented surface of genus $g$ and $n+1$ boundary components, which is obtained by removing a neighbourhood of the points labelled $1,2,\dotsc,n$ and $0$ from $\Sigma_g$ (see Figure \ref{fig:points}). The base point of $\Sigma$ is the point $*$ in Figure \ref{fig:points}, and consider the fundamental group $\pi = \pi_1(\Sigma, *)$ and the group algebra $\mathbb{K}\pi$.

Next, we recall some loop operations in \cite{mas}. As in Figure \ref{fig:points_uv}, we take two distinct points $u$ and $v$ near $*$, which is obtained by cutting out a small disk containing $*$ and inserting a disk with two points $u$ and $v$. Then, the groups $\pi_1(\Sigma, *)$, $\pi_1(\Sigma, u)$, and $\pi_1(\Sigma, v)$ are identified via paths in this disk neighbourhood. We also take two arcs $\gamma^u$ and $\gamma^v$, shown in the same figure.

\begin{figure}[bt]
\centerline{\includegraphics[scale=0.9]{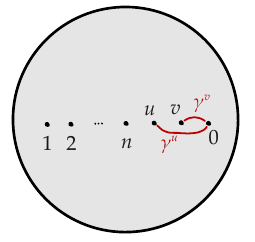}}
\caption{$n+3$ points inside the embedded disk with the arcs $\gamma^u$ and$\gamma^v$.}
\label{fig:points_uv}
\end{figure}

\begin{definition}
We define a linear map $\eta\colon \mathbb{K}\pi\otimes \mathbb{K}\pi \to \mathbb{K}\pi$ by the formula
\[
	\eta(\alpha,\beta) = \sum_{p\in\alpha\cap\beta} \mathrm{sign}(p;\alpha,\beta) \,\alpha_{up}\beta_{pv},
\]
where $\alpha,\beta\in\pi$ are represented by generically immersed curves on $\Sigma$, $\alpha$ is based at $u$, and $\beta$ is based at $v$, and the curves $\alpha$ and $\beta$ are taken so that they are disjoint from the interior of the arcs $\gamma^u$ and $\gamma^v$. We denote by $\mathrm{sign}(p;\alpha,\beta)$ the local intersection number with respect to the fixed orientation of $\Sigma$, and by $\alpha_{up}$ the path from $u$ to $p$ along $\alpha$, and similarly for $\beta_{pv}$. Regarding the path $\alpha_{up}\beta_{pv}$ from $u$ to $v$ as a loop based at $*$, we obtain an element of $\pi$.
\end{definition}

The map $\eta$ is a \textit{Fox pairing}: it satisfies
\begin{align*}
	\eta(\alpha\beta,\gamma) &= \eta(\alpha,\gamma) \varepsilon(\beta) + \alpha\,\eta(\beta,\gamma)\mbox{ and}\\
	\eta(\alpha,\beta\gamma) &= \eta(\alpha,\beta)\,\gamma + \varepsilon(\beta)\eta(\alpha,\gamma)
\end{align*}
where $\varepsilon$ is the augmentation map of the Hopf algebra $\mathbb{K}\pi$. 

We put $\vec\pi = \pi_1(U\Sigma, \vec *)$, where the base point $\vec *$ is taken as in Figure \ref{fig:points}.

\begin{definition}
We define a linear map $\vec\mu\colon \mathbb{K}\vec\pi \to \mathbb{K}\pi$ by the formula
\[
	\vec\mu(\vec\alpha) = \alpha + \hspace{-5pt}\sum_{p\in\mathrm{Self}(\alpha)} \mathrm{sign}(p;\alpha) \,\alpha_{*p}\alpha_{p*},
\]
where $\vec\alpha\in\vec\pi$ is represented by regularly immersed curve on $\Sigma$ based at $-\vec *$ and such that the projected curve $\alpha$ is disjoint from the interior of the short path connecting $*$ and $0$. Via the canonical identification $\vec\pi \cong \pi_1(UT\Sigma, -\vec*)$, we obtain the map.
\end{definition}

The map $\vec\mu$ is a \textit{quasi-derivation ruled by the Fox pairing} $\eta\colon$\!\! it satisfies the formula
\begin{align*}
	\vec\mu(\vec\alpha\vec\beta) = \vec\mu(\vec\alpha)\beta + \alpha\vec\mu(\vec\beta) + \eta(\alpha,\beta)
\end{align*}
for any $\vec\alpha, \vec\beta\in\vec\pi$.

These loop operations have useful expressions in terms of the Fox derivative. Following \cite{mas}, put $\pi^u = \pi_1(\Sigma\setminus\{u\},v)$ and $\pi^v = \pi_1(\Sigma\setminus\{v\},u)$. We have the canonical maps
\begin{align*}
	\iota^u&\colon \pi \cong \pi_1(\Sigma, v) \cong \pi_1(\Sigma\setminus\gamma^u,v) \hookrightarrow \pi^u\mbox{ and }\\
	\iota^v&\colon \pi \cong \pi_1(\Sigma, u) \cong \pi_1(\Sigma\setminus\gamma^v,u) \hookrightarrow \pi^v
\end{align*}
where the arcs $\gamma^u$ and $\gamma^v$ are as in Figure \ref{fig:points_uv} and the first isomorphisms are as remarked above. We also have the natural maps
\begin{align*}
	p^u&\colon \pi^u \to \pi_1(\Sigma,v)\cong \pi\;\mbox{ and }\;p^v\colon \pi^v \to \pi_1(\Sigma, u) \cong \pi
\end{align*}
by filling in the removed points. It is readily seen that $\iota^u$ and $\iota^v$ are injections, and $p^u$ and $p^v$ are surjections.

The group $\pi^u$ is regarded as the normal subgroup of $\mathrm{PB}_{g,12\cdots nuv0}$ (by which we denote the pure braid group with $n+3$ strands labelled by $1,2,\dotsc, n,u,v$ and $0$; we use similar notations throughout this paper) where the strands $1,2,\dotsc,n,u$ and $0$ are stationary. Similarly, $\pi^v$ is regarded as the normal subgroup of $\mathrm{PB}_{g,12\cdots nuv0}$ where the strands $1,2,\dotsc,n,v$ and $0$ are stationary. The elements $z^u\in\pi^u$ and $z^v\in\pi^v$ are both defined as the unique inverse image of $R^{u,v}R^{v,u}\in \mathrm{PB}_{g,12\cdots nuv0}$. We have the canonical decompositions
\begin{align*}
	\pi^u = \iota^u(\pi) * \langle z^u \rangle\;\mbox{ and }\;\pi^v = \iota^v(\pi) * \langle z^v\rangle.
\end{align*}
The Fox derivative $\frac{\partial }{\partial z^v}\colon \mathbb{K}\pi^v \to \mathbb{K}\pi^v$ is defined as the group $1$-cocycle satisfying
\[
	\frac{\partial }{\partial z^v}(z^v) = 1\quad\textrm{and}\quad \frac{\partial }{\partial z^v}(\iota^v(\alpha)) = 0 \textrm{ for any }\alpha\in\pi,
\]
which is well-defined by the above decomposition.

\begin{definition}
Define $c\colon\vec\pi \to \mathrm{PB}_{g,12\cdots nuv0}$ by that $c(\vec\alpha)$ is the pure braid where the strand $1,2,\dotsc,n$ and $0$ are stationary, $u$ traces $\alpha$ and $v$ traces the tip of the unit tangent vector of $\vec\alpha$. 
\end{definition}
\noindent Note that the convention for our tangent vector $\vec*$ is the opposite to the Massuyeau's one; see Section 3.2 of \cite{mas}. Accordingly, our definition of the map $c$ above is also adjusted so that the following remains true:

\begin{theorem}[Theorem 4.3 and 4.6 in \cite{mas}]
We have the following expressions:
\begin{enumerate}[(1)]
	\item $\eta(\alpha,\beta) = \left(p^v\circ\frac{\partial }{\partial z^v}\right) (\iota^u(\beta^{-1})\iota^v(\alpha)\iota^u(\beta))$;
	\item $\vec\mu(\vec\alpha) = \left(p^v\circ\frac{\partial }{\partial z^v}\right) (\iota^u(\alpha^{-1})c(\vec\alpha)\iota^v(\alpha)^{-1})$.
\end{enumerate}
\end{theorem}
\noindent The Fox derivatives $\frac{\partial}{\partial z^v}$\vspace{3pt} above make sense since the products $\iota^u(\beta^{-1})\iota^v(\alpha)\iota^u(\beta)$ and $\iota^u(\alpha^{-1})c(\vec\alpha)\iota^v(\alpha)^{-1}$ are elements of $\pi^v$: each defines the trivial braid by removing the $v$-th strand.\\

Now we define the (completed) Goldman--Turaev Lie bialgebra using the above maps based on \cite{mas} and \cite{akkn}. We have the \textit{weight filtration} on $\mathbb{K}\pi$, whose definition will be recalled in Definition \ref{def:tder}. We remark that the completion with respect to the weight filtration is isomorphic, in the category of topological Hopf algebras, to the completion by the augmentation ideal $\mathbb{I}\pi := \Ker(\varepsilon \colon \mathbb{K}\pi \to \mathbb{K})$. Therefore, as long as we do not take the associated graded quotient, we can identify these completions and we abusively denote them by $\widehat{\mathbb{K}\pi}$.

The map $\eta$ defines the map
\begin{align*}
	\kappa\colon \mathbb{K}\pi \otimes \mathbb{K}\pi &\to \mathbb{K}\pi \otimes \mathbb{K}\pi\\
	\alpha\otimes \beta &\mapsto \beta S(\eta(\alpha, \beta)) \alpha \otimes \eta(\alpha,\beta)
\end{align*}
for $\alpha,\beta\in\pi$ where $S$ is the antipode. This map is continuous with respect to the weight filtration with the filtration degree $(-2)$, so this induces the map on the completions:
\begin{align*}
	\kappa\colon \widehat{\mathbb{K}\pi} \,\hat\otimes\, \widehat{\mathbb{K}\pi} &\to \widehat{\mathbb{K}\pi} \,\hat\otimes \,\widehat{\mathbb{K}\pi}.
\end{align*}
This further induces the \textit{Goldman bracket} on the \textit{trace space} $|\widehat{\mathbb{K}\pi}| := \widehat{\mathbb{K}\pi}/[\widehat{\mathbb{K}\pi},\widehat{\mathbb{K}\pi}]$ by
\begin{align*}
	[\cdot, \cdot]_\mathrm{G}\colon |\widehat{\mathbb{K}\pi}| \otimes |\widehat{\mathbb{K}\pi}| &\to |\widehat{\mathbb{K}\pi}|\\
	|\alpha|\otimes |\beta| &\mapsto |\mathrm{mult}(\kappa(\alpha, \beta))|
\end{align*}
where $\mathrm{mult}\colon \mathbb{K}\pi\otimes \mathbb{K}\pi \to \mathbb{K}\pi$ is the multiplication map. On the other hand, the map $\vec\mu$ defines the map
\begin{align*}
	\delta_{\vec\mu}\colon \mathbb{K}\vec\pi &\to \mathbb{K}\pi \otimes \mathbb{K}\pi\\
	\vec\alpha &\mapsto \tau(\alpha S(\vec\mu(\vec\alpha)') \otimes \vec\mu(\vec\alpha)'')
\end{align*}
where we used Sweedler's notation $\Delta(x) = x'\otimes x''$ and the anti-symmetrisation map $\tau(\alpha\otimes \beta) = \alpha \otimes \beta - \beta \otimes \alpha$. The map $\delta_{\vec\mu}$ is also continuous with the weight-filtration degree $(-2)$, so we obtain the map
\begin{align*}
	\delta_{\vec\mu}\colon \widehat{\mathbb{K}\vec\pi} &\to \widehat{\mathbb{K}\pi} \,\hat\otimes\, \widehat{\mathbb{K}\pi}
\end{align*}
on the completions. This further induces the map 
\begin{align*}
	\vec\delta\colon |\widehat{\mathbb{K}\vec\pi}| &\to |\widehat{\mathbb{K}\pi}| \hat\otimes |\widehat{\mathbb{K}\pi}|
\end{align*}

\begin{definition}\label{def:framingiota}
A $\mathbb{K}$\textit{-framing} on $\Sigma$ is a group homomorphism $\mathsf{fr}\colon \vec\pi \to \mathbb{K}$ such that $\mathsf{fr}(F^*) = -1$ where $F^*$ denotes a negative monogon with respect to the orientation of $\Sigma$. We denote by $\mathrm{Fr}(\Sigma; \mathbb{K})$ the set of all $\mathbb{K}$-framings on $\Sigma$. Given a $\mathbb{K}$-framing, we define the map $\iota^\mathsf{fr}$ by
\begin{align*}
	\iota^\mathsf{fr}\colon \widehat{\mathbb{K}\pi} &\to \widehat{\mathbb{K}\vec\pi}\\
	\alpha &\mapsto \vec\alpha
\end{align*}
where $\vec\alpha$ is the rotation-free lift of $\alpha$: $\mathsf{fr}(\vec\alpha) = 0$.
\end{definition}

Then, for a $\mathbb{K}$-framing, we define the \textit{Turaev cobracket} associated with $\mathsf{fr}$ as the composition
\begin{align*}
	\delta^\mathsf{fr}\colon |\widehat{\mathbb{K}\pi}| \xrightarrow{\iota^\mathsf{fr}} |\widehat{\mathbb{K}\vec\pi}| \xrightarrow{\vec\delta} |\widehat{\mathbb{K}\pi}| \hat\otimes |\widehat{\mathbb{K}\pi}|.
\end{align*}
The triple $( |\widehat{\mathbb{K}\pi}|, [\cdot, \cdot]_\mathrm{G}, \delta^\mathsf{fr})$ constitutes a Lie bialgebra, the \textit{Goldman--Turaev Lie bialgebra}.\\

The \textit{formality problem} for the Goldman--Turaev Lie algebra asks if it is isomorphic to its associated graded quotient as a Lie bialgebra with respect to the weight filtration, and, if so, to determine the set of such isomorphisms. More precisely, a solution to the formality problem is a continuous Hopf algebra isomorphism
\[
	\theta \colon \widehat{\mathbb{K}\pi} \to \mathrm{gr}\,\widehat{\mathbb{K}\pi}
\]
such that $\mathrm{gr}(\theta) = \id$ (which we call an \textit{expansion} of $\mathbb{K}\pi$) and induces an isomorphism of Lie bialgebras $|\widehat{\mathbb{K}\pi}| \to |\mathrm{gr}\,\widehat{\mathbb{K}\pi}|$. Let $\mathrm{Form}^\mathsf{fr}_{g,n+1}$ be the set of all solutions for the formality problem associated with a $\mathbb{K}$-framing $\mathsf{fr}$, and we put
\[
	\mathrm{Form}_{g,n+1} := \bigsqcup_{\mathsf{fr}\in\mathrm{Fr}(\Sigma; \mathbb{K})} \hspace{-10pt}\mathrm{Form}_{g,n+1}^\mathsf{fr}\,.
\]

The formality problem itself is completely solved in \cite{akkn} (whose older version was published in 2018) by considering a set of equations, the \textit{Kashiwara--Vergne equations} of type $(g,n+1)$ for $\theta$, and in the special case of $\Sigma = \Sigma_{0,3}$ with a certain framing, this set of equations is surprisingly equivalent to the Kashiwara--Vergne problem in Lie theory. Furthermore, the solution set is acted on by the analogue of the KV and KRV groups, forming a bi-torsor. We note that another proof of the existence of the formality morphism for any genus $g$ was also obtained by Hain \cite{hain} via the theory of mixed Hodge structure.

On the other hand, the surprising relation between the KV problem (of type $(0,3)$) and (genus $0$) Drinfeld associators in the last section is described in the paper by Alekseev and Torossian \cite{at}: there is an \textit{inclusion}
\[
	\{\mbox{Drinfeld associators}\} \to \{\mbox{solutions to the KV equations}\},
\]
which is a morphism of bi-torsors. The purpose of this paper is to construct an analogous map for any $g\geq 0$.\\

From now on, we recall the KV equation of type $(g,n+1)$ introduced in \cite{akkn} together with ther weight filtration on $\mathbb{K}\pi$.

\begin{definition}\label{def:tder}
Assume $n\geq 0$.
\begin{itemize}
	\item Let $\mathcal{C} = (\alpha_i,\beta_i,\gamma_j)_{1\leq i\leq g,1\leq j\leq n}$ be a free-generating system of $\pi$ so that $\alpha_i$ and $\beta_i$ form a genus pair, $\gamma_j$ is a boundary loop representing the $j$-th boundary $\partial_j\Sigma$ and
	\[
		\alpha_1\beta_1\alpha_1^{-1}\beta_1^{-1}\cdots \alpha_g\beta_g\alpha_g^{-1}\beta_g^{-1}\gamma_1\cdots\gamma_n
	\]
	represents a boundary loop based at $*\in\partial_0\Sigma$ (see Figure 2 of \cite{toyo2}). Denote by $(x_i,y_i,z_j)_{1\leq i\leq g,1\leq j\leq n}$ the corresponding basis of $H = H_1(\Sigma;\mathbb{K})$.
	\item $\hat T(H) \cong \mathbb{K}\langle\!\langle x_i,y_i,z_j \rangle\!\rangle_{1\leq i\leq g, 1\leq j\leq n}$ is the completed free associative algebra over $H$ with respect to the \textit{weight grading} on $H$ defined by $\mathrm{wt}(x_i) = \mathrm{wt}(y_i) = 1$ and $\mathrm{wt}(z_j) = 2$. Then, $\mathrm{gr}\,\widehat{\mathbb{K}\pi}$ is canonically isomorphic to $\hat T(H)$ as a completed Hopf algebra.
	\item Consider the morphism of $\mathbb{K}$-algebras
	\[
		\theta_\mathrm{exp}\colon\mathbb{K}\pi \to \hat T(H)\colon \alpha_i \mapsto e^{x_i}, \beta_i \mapsto e^{y_i} \textrm{ and } \gamma_j \mapsto e^{z_j}.
	\]
	Define the \textit{weight filtration} on $\mathbb{K}\pi$ by the pull-back of the weight filtration by $\theta_\mathrm{exp}$. This induces an isomorphism of Hopf algebras on the completion $\widehat{\mathbb{K}\pi}$.
	\end{itemize}

	Next, we have some spaces:
	\begin{itemize}
	\item $\hat L(H) = L(\!(x_i,y_i,z_j)\!)_{1\leq i\leq g, 1\leq j\leq n}$ is the completed free Lie algebra over $H$, so that $\hat T(H)$ is identified with the (completed) universal enveloping algebra $U\!\hat L(H)$,
	\item $\Der^+(\hat L(H)) = \{u\colon \hat L(H)\to \hat L(H): \textrm{a continuous Lie algebra derivation, degree}\geq 1\}$,
	\item $\Aut^+(\hat L(H)) = \Exp(\Der^+(\hat L(H)))$,
	\item the space of \textit{tangential derivations}:
	\begin{align*}
		\tDer^+(\hat L(H)) = \{\tilde u=(u;u_1,\dotsc,u_n): u\in \Der^+(\hat L(H)), \,u_i\in \hat L(H), \,u(z_j) = [z_j,u_j]\},
	\end{align*}
	\item the space of \textit{tangential automorphisms}:
	\begin{align*}
		\tAut^+(\hat L(H)) &= \Exp(\tDer^+(\hat L(H)))\\
		& = \{\tilde G=(G;g_1,\dotsc,g_n):G\in \Aut^+(\hat L(H)),\, g_j\in\Exp(\hat L(H)), \,G(z_j) = g_j^{-1}z_jg_j\},
	\end{align*}
\end{itemize}
and some elements:
\begin{itemize}
	\item $\xi = \log\left(\prod_i(e^{x_i}e^{y_i}e^{-x_i}e^{-y_i})\prod_j e^{z_j}\right)\in \hat L(H)$,
	\item $\omega = \sum_i [x_i,y_i] + \sum_j z_j\in \hat L(H)$,
	\item $r(s) = \log\left(\frac{e^s-1}{s}\right)\in s\K s$, $\mathbf{r} = \sum_i |r(x_i) + r(y_i)|$,
	\item $a_i = \rot^\mathsf{fr}(\alpha_i),\,b_i = \rot^\mathsf{fr}(\beta_i),\,c_j = \rot^\mathsf{fr}(\gamma_j) \in \mathbb{K}$, and $\mathbf{p}^\mathsf{fr} = \sum_i |a_iy_i - b_ix_i|\in |\hat T(H)|$.
	\end{itemize}
\end{definition}

Now recall the standard divergence for a free Lie algebra and the integration of a $1$-cocycle.
\begin{definition}\ 
\begin{itemize}
	\item $\hat T(H)$ is regarded as an $\hat L(H)$-module by the left multiplication. For $w = x_i,y_i,z_j$, $d_w\colon \hat L(H) \to \hat T(H)$ is a continuous Lie algebra $1$-cocycle specified by $d_w(w') = \delta_{ww'}$ for $w' = x_i,y_i,z_j$ using Kronecker's delta.
	\item We define the \textit{single divergence} $\sdiv_{x,y,z} \colon \Der(\hat L(H)) \to |\hat T(H)|$ by 
	\[
		\sdiv_{x,y,z}(u) = \sum_{w = x_i,y_i,z_j} |d_wu(w)|.
	\]
	This is extended to $\tDer(\hat L(H))$ by the composition
	\[
		\sdiv_{x,y,z}\colon \tDer^+(\hat L(H)) \to \Der(\hat L(H)) \xrightarrow{\sdiv_{x,y,z}} |\hat T(H)|.
	\]
	The single divergence is itself a Lie algebra $1$-cocycle.
	\item For a pro-nilpotent Lie algebra $\mathfrak g$, a continuous $\mathfrak g$-module $V$ and a $1$-cocycle $\psi\colon \mathfrak g \to V$, its \textit{integration} is a group $1$-cocycle $\Psi\colon \exp(\mathfrak g) \to V$ given by, for $u\in \mathfrak g$,
	\[
		\Psi(e^u) = \frac{e^u-1}{u}\cdot \psi(u).
	\]
	The correspondence $\psi \mapsto \Psi$ is $\mathbb{K}$-linear. For the details, see Appendix A of \cite{akkn}.
	\item Since $\tDer^+(\hat L(H))$ is pro-nilpotent, we denote the integration of $\sdiv_{x,y,z}\colon \tDer^+(\hat L(H)) \to |\hat T(H)|$ by $$\mathsf{j}_{x,y,z}\colon \tAut^+(L) \to |\hat T(H)|.$$
\end{itemize}
\end{definition}

\noindent We also have many $1$-cocycles:
\begin{itemize}
	\item $\mathsf{b}^\mathsf{fr}\colon\tDer^+(\hat L(H))\to |\hat T(H)|:\tilde u\;\mapsto\;\sum_j c_j|u_j|$ and $\mathsf{c}^\mathsf{fr}\colon\tAut^+(\hat L(H))\to |\hat T(H)|$ its integration, \\[-10pt]
	\item $\sdiv^\mathsf{fr}\colon \tDer^+(\hat L(H))\to |\hat T(H)|: \tilde u\;\mapsto\; \sdiv_{x,y,z}(u)- \mathsf{b}^\mathsf{fr}(\tilde u) + u(\mathbf{r} - \mathbf{p}^\mathsf{fr})$ and $\mathsf{j}^\mathsf{fr}\colon \tAut^+(\hat L(H))\to |\hat T(H)|$ its integration, and
	\item $\sdiv^\mathsf{fr}_\mathrm{gr}\colon \tDer^+(\hat L(H))\to |\hat T(H)|: \tilde u\;\mapsto\; \sdiv_{x,y,z}(u)- \mathsf{b}^\mathsf{fr}(\tilde u)$ and $\mathsf{j}^\mathsf{fr}_\mathrm{gr}\colon \tAut^+(\hat L(H))\to |\hat T(H)|$ its integration.\\
\end{itemize}

Finally, we recall the definition of the KV groups and associators.
\begin{definition} For $g,n\geq 0$ and a $\mathbb{K}$-framing $\mathsf{fr}$, the \textit{Kashiwara--Vergne group} $\mathrm{KV}^\mathsf{fr}_{g,n+1}$, the graded version $\mathrm{KRV}^\mathsf{fr}_{g,n+1}$, and the set of the \textit{Kashiwara--Vergne associators} $\mathrm{SolKV}^\mathsf{fr}_{g,n+1}$ is defined by the followings:
\begin{itemize}
	\item $\mathrm{KV}^\mathsf{fr}_{g,n+1} = \Big\{\tilde G \in\tAut^+(\hat L(H)): G(\xi) = \xi,\,\mathsf{j}^\mathsf{fr}(\tilde G) \in \Big|\sum_j z_j\K{z_j} + \xi^2\K\xi\Big|\Big\}$,
	\item $\mathrm{KRV}^\mathsf{fr}_{g,n+1} = \Big\{\tilde G \in\tAut^+(\hat L(H)): G(\omega) = \omega,\,\mathsf{j}^\mathsf{fr}_\mathrm{gr}(\tilde G) \in \Big|\sum_j z_j\K{z_j} + \omega^2\K\omega\Big|\Big\}$, and
	\item $\mathrm{SolKV}^\mathsf{fr}_{g,n+1} = \Big\{\tilde G\in\tAut^+(\hat L(H)): G(\omega)=\xi,\,\mathsf{j}^\mathsf{fr}_\mathrm{gr}(\tilde G) -\mathbf{r} + \mathbf{p}^\mathsf{fr} \in \Big|\sum_j z_j\K{z_j} + \xi^2\K\xi\Big|\Big\}$.
\end{itemize}
The set $\mathrm{SolKV}^\mathsf{fr}_{g,n+1} $ is a bi-torsor over the groups $\mathrm{KV}^\mathsf{fr}_{g,n+1}$ and $\mathrm{KRV}^\mathsf{fr}_{g,n+1}$, which is apparent from their defining equations.
\end{definition}

One of the main results in \cite{akkn} is the following.

\begin{theorem}[\cite{akkn}, Theorem 6.27]\label{thm:akkn}
For $\Sigma = \Sigma_{g,n+1}$ with $n\geq 0$, an isomorphism of filtered Hopf algebras $\theta\colon\widehat{\mathbb{K}\pi} \to \hat T(H)$ with $\mathrm{gr}(\theta) = \id$ gives a solution to the formality problem if and only if $\theta\circ\theta_\mathrm{exp}^{-1}$ lifts to an element in $\mathrm{SolKV}_{g,n+1}$ up to conjugation by an element of $\mathrm{exp}(\hat L(H))$.
\end{theorem}
We remark that all the calculations in \cite{akkn} are done for $\mathbb{Z}$-framings, but they remain valid for $\mathbb{K}$-framings. Now consider the diagram
\[\begin{tikzcd}[cramped]
	\mathrm{SolKV}^\mathsf{fr}_{g,n+1} \times \mathrm{exp}(\hat L(H)) \arrow[r]& \tAut^+(\hat L(H))\arrow[d]\\
	\mathrm{Form}^\mathsf{fr}_{g,n+1}  \arrow[r, hookrightarrow]& \mathrm{Isom}^+_\mathrm{Hopf}(\mathbb{K}\pi, \hat T(H))
\end{tikzcd}\]
where $\mathrm{Isom}^+_\mathrm{Hopf}(\mathbb{K}\pi, \hat T(H))$ is the set of all expansions of $\mathbb{K}\pi$, the map of the first row is given by $(\tilde G, g) \mapsto \tilde G\circ\mathrm{Ad}_{g}$ and the vertical map is given by $\tilde G \mapsto G^{-1}\circ \theta_\mathrm{exp}$. The theorem above states that this induces a surjective map
\[
	\mathrm{SolKV}^\mathsf{fr}_{g,n+1} \times \mathrm{exp}(\hat L(H)) \twoheadrightarrow \mathrm{Form}^\mathsf{fr}_{g,n+1}.
\]
We denote the set of all KV associators by
\begin{align*}
	\mathrm{SolKV}_{g,n+1} := \bigsqcup_{\mathsf{fr}\in\mathrm{Fr}(\Sigma; \mathbb{K})} \hspace{-10pt}\mathrm{SolKV}_{g,n+1}^\mathsf{fr}\,.
\end{align*}

\section{The Main Construction}\label{sec:main}
Let $g,n\geq 0$. The purpose of this section is the construction of the map
\[
	\{\mbox{genus } g \mbox{ Gonzalez--Drinfeld associators}\} \to \mathrm{Form}_{g,n+1},
\]
which is done by extending the method by Massuyeau \cite{mas}. We only consider when the coupling constant $\mu$ is equal to 1; if not, we may apply the rescaling automorphism $\lambda = \sqrt\mu$ to normalise (if a square root exists).\\

Let $\vec Z\colon \widehat{\mathbb{K}\mathbf{PaB}}{}^f_g \to \mathbf{PaCD}^f_g$ be a genus $g$ Gonzalez--Drinfeld associator and $\varepsilon^i$ the operadic composition of the empty braid/chord into the $i$-th strand. Fix a maximal parenthesisation $\ell_n$ of the sequence $1\;2\;\cdots\; n$ and let $p = (\ell_n)(* \;0)$ be a maximal parenthesisation of $n+2$ letters $1,2,\dotsc,n$, $*$ and $0$, which represents the points on the surface $\Sigma_g$ in Figure \ref{fig:points}. Put $\Sigma = \Sigma_{g,n+1}$, $\pi = \pi_1(\Sigma,*)$ and $\vec\pi = \pi_1(U\Sigma,\vec*)$ as before.

For a group $G$, we denote by $\hat G$ the Malcev completion over $\mathbb{K}$, which is defined as the group-like part of the completion $\widehat{\mathbb{K}G}$ of the group algebra by the augmentation ideal.

\begin{lemma}\label{lem:exactness}
Denoting by $\vec Z_p$ the evaluation of $\vec Z$ at $p$, we have the following commutative diagram of (completed) groups with rows exact:
\[\begin{tikzcd}[cramped]
	1 \arrow[r]& \widehat{\vec\pi} \arrow[r]\arrow[d,"\vec Z_p"]& \widehat{\mathrm{PB}}^f_{g,12\cdots n*0}\arrow[r,"\varepsilon^*"]\arrow[d,"\vec Z_p"]& \widehat{\mathrm{PB}}^f_{g,12\cdots n0} \arrow[r]\arrow[d,"\vec Z_{\varepsilon^*(p)}"] &1\phantom{\,,}\\
	1 \arrow[r]& \exp(\hat L(H)\oplus\mathbb{K}t_{**})\arrow[r]& \exp(\mathfrak{t}^f_{g,12\cdots n*0}) \arrow[r,"\varepsilon^*"]& \exp(\mathfrak{t}^f_{g,12\cdots n0}) \arrow[r]&1\,,
\end{tikzcd}\]
where $\vec\alpha\in\vec\pi$ is identified with the braid such that the strands $1,2,\dotsc, n$ and $0$ are stationary and the $*$-th strand traces the framed loop specified with $\vec\alpha$, $H = \mathrm{Span}_\mathbb{K}\{t_{j*}, x^i_*, y^i_*\}_{1\leq j\leq n, 1\leq i\leq g}$ is (isomorphic to) the first homology group $H_1(\Sigma_{g,n+1};\mathbb{K})$, and $\hat L(H)\oplus\mathbb{K}t_{**}$ is the direct sum of the complete free Lie algebra over $H$ and one-dimensional Lie algebra $\mathbb{K}t_{**}$ spanned by $t_{**}$.
\end{lemma}

\noindent The lengthy proof is deferred to Section \ref{sec:proof}. Therefore, the map $\vec Z_p$ satisfies the axiom of an expansion on $\vec\pi$ except for $\mathrm{gr}(\vec Z_p) = \id$, which will turn out to give a solution to the KV problem associated with the surface $\Sigma$. Note that the map $\vec Z_p$ induces $Z_p\colon \widehat{\mathrm{PB}}_{g,12\cdots n*0} \xrightarrow{\cong} \exp(\mathfrak{t}_{g,12\cdots n*0})$ on the quotient spaces of non-framed pure braids.\\

From now on, we shall do some preparatory computations involving $\vec Z$. First of all, the element $\vec\omega\in \mathfrak{t}^f_{g,12\cdots n*0}$ is defined by 
\[
	\vec\omega = \sum_{1\leq i\leq g} [x^i_*, y^i_*] + \sum_{1\leq j\leq n} t_{j*}.
\]
The image by the quotient map $\mathfrak{t}^f_{g,12\cdots n*0} \twoheadrightarrow \mathfrak{t}_{g,12\cdots n*0}$ is denoted by $\omega$.

\begin{lemma}\label{lem:special}
Let $\gamma_0\in\pi$ be a simple loop traversing the $0$-th boundary in the opposite direction and $\vec\gamma_0$ be the rotation-free lift with respect to the framing $\mathsf{fr}$ with $\mathsf{fr}(\gamma_i) = -1$ for all $1\leq i\leq n$. Then, we have $\vec Z_p(\vec\gamma_0) = e^{\vec\omega}$ in $\exp(\mathfrak{t}^f_{g,12\cdots n*0})$.
\end{lemma}
\noindent Proof. We denote by $F^*\in\Hom_{\mathbf{PaB}^f}(p,p)$ the braid such that the $*$-th strand draws the negative monogon and other strands are stationary. Then, we have $\vec\gamma_0 = (F^*)^s (R^{*,0})^{-1}(R^{0,*})^{-1}$ with $s = 2g - 2$ by the Poincar\'e--Hopf theorem. Next, we have $\vec Z_p(R^{*,0}R^{0,*}) = e^{t_{*0}}$ since $*$ and $0$ are closely adjecent in $p$. Therefore, we have
\begin{align*}
	\vec Z_p(\vec\gamma_0) &= \vec Z_p((F^*)^s(R^{*,0})^{-1}(R^{0,*})^{-1})\\
	&= e^{s t_{**}/2}e^{-t_{*0}}\\
	&= e^{[s/2+ (1-g)]t_{**} + \vec\omega}.
\end{align*}
The last equality comes from the defining relation of $\mathfrak{t}^f_{g,12\cdots n*0}$ and the fact that $t_{**}$ is central. Since $s/2+ (1-g) = 0$, we have $\vec Z_p(\vec \gamma_0) = e^{\vec \omega}$ as claimed. \qed

\begin{remark}
By the above lemma, we have $\mathrm{gr}\,\vec Z_p(\vec \omega) = \vec\omega$. Furthermore, from the assumption $\mu = 1$, we have $\mathrm{gr}\,\vec Z_p(t_{i*}) = t_{i*}$ for $i \in \{1,2,\dotsc, n, *\}$ and hence $\mathrm{gr}\,\vec Z_p$ preserves the sum $\sum_{1\leq i\leq g} [x^i_*, y^i_*]$. Therefore, by applying the $\mathrm{Sp}(2g;\mathbb{K})$-action if necessary, we can choose $\vec Z$ so that $\mathrm{gr}\,\vec Z_p = \id$. We require this condition from now on, although this is used only to guarantee $\mathrm{gr}\,\vec Z_p$ to be an expansion, aligning with the convention in \cite{akkn}. The calculation below holds as is without this requirement.
\end{remark}

To apply Massuyeau's description, we consider another object $p' = p\circ_* (u\;v) = (\ell_n)((u\;v)\;0)$ in $\mathbf{PaB}_g^f$. Under the identification $(12\cdots nu0) = (12\cdots n\!*\!0) = (12\cdots nv0)$, put
\begin{align*}
	P^u &= U(\Ker(\varepsilon^u\colon \mathfrak{t}_{g,12\cdots nuv0} \to \mathfrak{t}_{g,12\cdots n*0})),\\
	P^v &= U(\Ker(\varepsilon^v\colon \mathfrak{t}_{g,12\cdots nuv0} \to \mathfrak{t}_{g,12\cdots n*0}))\mbox{ and }\\
	P^{u,v} &= U(\Ker\varepsilon^u\cap\Ker\varepsilon^v).
\end{align*}
Since $p^u$ is the restriction of $\varepsilon^u$ on $\mathrm{PB}_{g,12\cdots nuv0}$, the map $Z_{p'}$ restricts to $Z_{p'}\colon \widehat{\mathbb{K}\pi^u} \xrightarrow{\cong} P^v$ and similarly for $\pi^v$.

\begin{lemma}[see Lemma 8.1 in \cite{mas}]
Define the map $D^v$ so that the following diagram is commutative:
\[\begin{tikzcd}[cramped]
	\widehat{\mathbb{K}\pi^v} \arrow[r,"\frac{\partial}{\partial z^v}"]\arrow[d,"Z_{p'}"] & \widehat{\mathbb{K}\pi^v} \arrow[r,"p^v"] & \widehat{\mathbb{K}\pi} \arrow[d,"Z_p"]\phantom{\,.}\\
	P^u \arrow[rr,"D^v"] && U(\mathfrak{t}_{g,12\cdots n*0})\,.
\end{tikzcd}\]
Then, there is a constant $\phi \in U(\mathfrak{t}_{g,12\cdots n*0})$ depending on $\vec Z$ such that
\begin{align}
	&D^v(xy) = D^v(x)\varepsilon(y) + \varepsilon^v(x)D^v(y) \quad \mbox{for } x,y\in P^u, \label{eq:vfoxder}\\
	&D^v(x^\times) = S(D^v(x)) \quad \mbox{for } x\in P^{u,v}, \label{eq:vantipode}\\
	&D^v(t_{uv}) = 1,\quad D^v(t_{ju}) = -t_{j*}\phi,\quad D^v(x^i_u) = -x^i_*\phi,\quad D^v(y^i_u) = -y^i_*\phi, \label{eq:values_Dv}\\
	& \phi - S(\phi) = \frac12 + s(\omega)\,. \label{eq:difference}
\end{align}
Here, $\varepsilon$ is the augmentation map, $S$ is the antipode, $s(\omega) = \frac{e^\omega}{1 - e^\omega} + \frac{1}{\omega}$ and the superscript ${}^\times$ indicates the letters $u$ and $v$ are swapped within the element.
\end{lemma}
\noindent Proof. \eqref{eq:vfoxder}: For $x,y\in P^v$, put $\alpha = Z_{p'}^{-1}(x)$ and $\beta = Z_{p'}^{-1}(y)$. Then, since $Z_p$ and $p^v$ are (complete) Hopf algebra homomorphisms, we have
\begin{align*}
	D^v(xy) &= D^v(Z_{p'}(\alpha\beta))\\
	&= \left(Z_p\circ p^v\circ\frac{\partial}{\partial z^v} \right)(\alpha\beta)\\
	&= (Z_p\circ p^v)\left(\frac{\partial \alpha}{\partial z^v} \varepsilon(\beta) + \alpha\frac{\partial \beta}{\partial z^v}\right)\\
	&= D^v(x) \varepsilon(y) + (Z_p\circ p^v)(\alpha) D^v(y).
\end{align*}
As we noted, the map $p^v$ is the restriction of the operadic map $\varepsilon^v$, so we have
\begin{align*}
	\varepsilon^v\circ Z_{p'} = Z_{\varepsilon^v(p')}\circ \varepsilon^v =  Z_{p}\circ p^v
\end{align*}
and hence
\begin{align*}
	(Z_p\circ p^v)(\alpha) D^v(y) = (\varepsilon^v\circ Z_{p'})(\alpha) D^v(y) = \varepsilon^v(x) D^v(y).
\end{align*}

\noindent\eqref{eq:vantipode}: For a parenthesisation $q$ with the underlying sequence of letters $|q|$ and $\beta\in\mathrm{PB}^f_{g,|q|}$, denote by $\beta_q\in \Hom_{\mathbf{PaB}^f_g}(q,q)$ the braid $\beta$ with the parenthesisation given by $q$. Now take $x\in P^{u,v}$ and put $\alpha = Z_{p'}^{-1}(x)$. By Lemma 4.1 of \cite{mas}, we have
\[
	 \left(Z_p\circ p^v\circ\frac{\partial}{\partial z^v} \right)(\sigma\alpha\sigma^{-1}) =  \left(Z_p\circ S\circ p^v\circ\frac{\partial}{\partial z^v} \right)(\alpha) = S(D^v(x)).
\]
On the other hand, by \eqref{eq:vfoxder},
\begin{align*}
	 \left(Z_p\circ p^v\circ\frac{\partial}{\partial z^v} \right)(\sigma\alpha\sigma^{-1}) &= (D^v\circ Z)((R^{v,u})^{-1}\alpha_{(\ell_n)((v\;u)\;0)}(R^{v,u}))\\
	 &= D^v(e^{-t_{uv}/2}) + D^v(x^\times) + D^v(e^{t_{uv}/2})\\
	 &= D^v(x^\times).
\end{align*}

\noindent\eqref{eq:values_Dv}: Putting $T  = Z_{p'}(z^v) = Z_{p'}(R^{u,v}R^{v,u}) = e^{t_{uv}}$, we have $\varepsilon(T-1)= \varepsilon^v(T-1) = 0$. By \eqref{eq:vfoxder}, we have
\begin{align*}
	D^v(t_{uv}) &= D^v(\log T) = D^v\left(\sum_{m\geq 1} \frac{(-1)^{m-1}}{m} (T-1)^m \right)\\
	&= D^v(T-1)\\
	&= \left(Z_p\circ p^v\circ\frac{\partial}{\partial z^v} \right)(z^v - 1)\\
	&= (Z_p\circ p^v)(1) = 1.
\end{align*}

Next, let $p'' = (\ell_n)(u\;(v\;0))$, $\Phi = \id_{((\ell_n)\;\star)}\circ_\star\,\Phi^{u,v,0}\in \Hom_{\mathbf{PaB}^f_g}(p',p'')$ and put $\varphi = Z(\Phi^{-1})$ and $\phi = D^v(\varphi)$. We have $\varepsilon(\varphi) = \varepsilon^v(\varphi) = 1$ since $\varepsilon^v(p') = \varepsilon^v(p'') = p$. For any operadic element $X$ with inputs marked $a,b,\dotsc,c$, we put $X_{a,b,\dotsc,c}^{a',b',\dotsc,c'}$ the operadic element obtained by composing $\id_{x'}$ to the $x$-th slot for each letter $x = a, b ,\dotsc, c$. Then, for $\alpha\in\pi$, we have $\iota^v(\alpha)_{p''} = (\alpha_p)_{*,0}^{u,(v\;0)}$ and hence
\begin{align*}
	0 &= \left(Z_p\circ p^v\circ\frac{\partial}{\partial z^v} \right)(\iota^v(\alpha))\\
	&= (D^v\circ Z)(\Phi\iota^v(\alpha)_{p''}\Phi^{-1})\\
	&= D^v\Big(\varphi^{-1} Z((\alpha_p)_{*,0}^{u,(v\;0)})\varphi\Big)\\
	&= D^v\Big(\varphi^{-1} (Z(\alpha_p)_{*,0}^{u,(v\;0)})\varphi\Big)\\
	&= -\phi + D^v(Z(\alpha_p)_{*,0}^{u,(v\;0)}) + \varepsilon^v(Z(\alpha_p)_{*,0}^{u,(v\;0)}) \phi\\
	&= D^v(Z(\alpha_p)_{*,0}^{u,(v\;0)}) + (Z(\alpha_p) - 1)\phi.
\end{align*}
Therefore, for $A\in \exp(\hat L(H))$, we have
\[
	D^v(A_{*,0}^{u,(v\;0)}) = (\varepsilon(A) - A)\phi \quad\mbox{in } U(\mathfrak{t}_{g,12\cdots n*0}),
\]
which, in turn, implies the same equality for any $A\in U(\hat L(H))$ since $\exp(\hat L(H))$ is linearly dense in $U(\hat L(H))$. Substituting $A = t_{j*}, x^i_*, y^i_*$ yields
\begin{align*}
	D^v(t_{ju}) = -t_{j*}\phi,\quad D^v(x^i_u) = -x^i_*\phi\;\mbox{ and }\; D^v(y^i_u) = -y^i_*\phi.
\end{align*}

\noindent\eqref{eq:difference}: Put $\tau = \iota^v(\gamma_0)z^v \in \pi^v$. Then, we have
\[
	\tau_{p'} = (R^{v,u})^{-1}\Phi^\times(R^{0,u})^{-1}(R^{u,0})^{-1}(\Phi^\times)^{-1}R^{v,u}.
\]
On the one hand, we have
\begin{align*}
	\left(Z_p\circ p^v\circ\frac{\partial}{\partial z^v} \right)(\tau) = (Z_p\circ p^v)(\iota^v(\gamma_0)) = Z_p(\gamma_0) = e^\omega
\end{align*}
by Lemma \ref{lem:special}. On the other hand, we have $\varepsilon(\varphi^\times) = \varepsilon^v(\varphi^\times) = 1$ since $\varepsilon^u(p') = \varepsilon^u(p'') = p$ and hence
\begin{align*}
	(D^v\circ\, Z)(\tau_{p'}) &= D^v\left(e^{-t_{uv}/2} (\varphi^\times)^{-1}e^{-t_{u0}}\varphi^\times e^{t_{uv}/2}\right)\\
	&= -\frac12 - S(\phi) + D^v(e^{-t_{u0}}) + e^{-t_{*0}}S(\phi) + \frac{e^{-t_{*0}}}2\\
	&= \left(\frac12 + S(\phi)\right) (e^\omega - 1) + D^v(e^{-t_{u0}}).
\end{align*}
Here we used $D^v(\varphi^\times) = S(\phi)$, which is a consequence of $\varphi\in P^{u,v}$ and \eqref{eq:vantipode}. Now we compute $D^v(e^{-t_{u0}})$. Since $\varepsilon(t_{u0}) = 0$, we have
\begin{align*}
	D^v(e^{-t_{u0}}) &= D^v\left(\sum_{m\geq 0} \frac{(-t_{u0})^m}{m!}\right)\\
	&= \sum_{m\geq 1}\frac{(-t_{*0})^{m-1}}{m!} D^v(-t_{u0})\\
	&= \frac{e^\omega - 1}{\omega} D^v\left( \sum_{1\leq i\leq g} [x^i_u, y^i_u] + \sum_{1\leq j\leq n} t_{ju} + t_{uv} \right)\\
	&= \frac{e^\omega - 1}{\omega} (-\omega\phi + 1).
\end{align*}
The last equality comes from \eqref{eq:values_Dv}. Therefore, we have
\[
	e^\omega = \left(\frac12 + S(\phi)\right) (e^\omega - 1) + \frac{e^\omega - 1}{\omega} (-\omega\phi + 1),
\]
which is equivalent to
\[
	 \phi - S(\phi) = \frac{e^\omega}{1 - e^\omega} + \frac12 + \frac{1}{\omega}\,.
\]
This concludes the proof.\qed\\

\begin{lemma}[see Lemma 8.2 in \cite{mas}]
Define the map $D^v$ so that the following diagram is commutative:
\[\begin{tikzcd}[cramped]
	\widehat{\mathbb{K}\pi^u} \arrow[r,"\frac{\partial}{\partial z^u}"]\arrow[d,"Z_{p'}"] & \widehat{\mathbb{K}\pi^u} \arrow[r,"p^u"] & \widehat{\mathbb{K}\pi} \arrow[d,"Z_p"]\phantom{\,.}\\
	P^v \arrow[rr,"D^u"] && U(\mathfrak{t}_{g,12\cdots n*0})\,.
\end{tikzcd}\]
Then, there is a constant $\phi \in U(\mathfrak{t}_{g,12\cdots n*0})$ depending on $\vec Z$ such that
\begin{align}
	&D^u(xy) = D^u(x)\varepsilon(y) + \varepsilon^u(x)D^u(y) \quad \mbox{for } x,y\in P^u, \label{eq:ufoxder}\\
	&D^u(x) = S(D^v(x)) \quad \mbox{for } x\in P^{u,v}, \label{eq:uantipode}\\
	&D^u(t_{uv}) = 1,\quad D^v(t_{jv}) = -t_{j*}\bar\phi,\quad D^v(x^i_v) = -x^i_*\bar\phi\;\mbox{ and }\; D^v(y^i_v) = -y^i_*\bar\phi. \label{eq:values_Du}
\end{align}
Here we put $\bar\phi = \phi + \frac12$.
\end{lemma}
\noindent Proof. \eqref{eq:ufoxder} is similarly done as the previous lemma.

\noindent\eqref{eq:uantipode}: Take $x\in P^{u,v}$ and put $\alpha = Z_{p'}^{-1}(x)$. By Lemma 4.1 of \cite{mas}, we have
\begin{align*}
	D^u(x) &= \left(Z_p\circ p^u\circ\frac{\partial}{\partial z^u} \right)(\alpha)\\
	&= \left(S\circ Z_p\circ p^v\circ\frac{\partial}{\partial z^v} \right)(\alpha)\\
	&= S(D^v(x)).
\end{align*}

\noindent\eqref{eq:values_Du}: First, since $t_{uv}\in P^{u,v}$, we have
\begin{align*}
	D^u(t_{uv}) = S(D^v(t_{uv})) = S(1) = 1.
\end{align*}
Next, we have $\iota^u(\alpha) = (R^{v,u})^{-1}\Phi^\times(\alpha_p)_{*,0}^{v,(u\;0)}(\Phi^\times)^{-1}R^{v,u}$ for $\alpha\in\pi$ and hence
\begin{align*}
	Z_{p'}(\iota^u(\alpha)) &= Z((R^{v,u})^{-1}\Phi^\times(\alpha_p)_{*,0}^{v,(u\;0)}(\Phi^\times)^{-1}R^{v,u})\\
	&= e^{-t_{uv}/2} (\varphi^\times)^{-1} Z(\alpha)_{*,0}^{v,(u\;0)} \varphi^\times e^{t_{uv}/2}.
\end{align*}
Therefore,
\begin{align*}
	0 &= \left(Z_p\circ p^u\circ\frac{\partial}{\partial z^u} \right)(\iota^u(\alpha))\\
	&= D^u(Z_{p'}(\iota^u(\alpha)))\\
	&= D^u(e^{-t_{uv}/2} (\varphi^\times)^{-1} Z(\alpha)_{*,0}^{v,(u\;0)} \varphi^\times e^{t_{uv}/2})\\
	&= D^u(Z(\alpha)_{*,0}^{v,(u\;0)}) + (Z(\alpha)_{*,0}^{v,(u\;0)} - 1)\bar\phi
\end{align*}
By the density argument, we have, for $A\in U(\hat L(H))$,
\begin{align*}
	D^u(A_{*,0}^{v,(u\;0)}) = (\varepsilon(A) - A_{*,0}^{v,(u\;0)})\bar\phi.
\end{align*}
Substituting $A = t_{j*}, x^i_*, y^i_*$ yields
\begin{align*}
	D^v(t_{jv}) = -t_{j*}\bar\phi,\quad D^v(x^i_v) = -x^i_*\bar\phi\;\mbox{ and }\; D^v(y^i_v) = -y^i_*\bar\phi.
\end{align*}
This completes the proof.\qed\\

\section{Proof of the Formality}\label{sec:threedim}
With these calculations done, we show that $Z_p$ gives the expansion we want.

\begin{definition}
The Fox pairing $(-\!\odot\!-)$ on $U(\hat L(H))$ is defined on generators by
\begin{align*}
	x_*^i \odot y_*^j = \delta_{ij},\quad y_*^i \odot x_*^j = -\delta_{ij},\quad t_{i*} \odot t_{j*} = -\delta_{ij}t_{j*},
\end{align*}
and other pairings are zero.
We also define the map $\xi\colon U(\hat L(H)\oplus \mathbb{K}t_{**}) \to U(\hat L(H))$ by
\begin{align*}
	\xi(z_1\cdots z_r (t_{**})^s) = \left\{\begin{aligned}
		&\sum_{1\leq i< r} z_1\cdots z_{i-1} (z_i\odot z_{i+1}) z_{i+2}\cdots z_r & & \mbox{if } s=0\mbox{ and } r\geq 2,&\\
		&0& & \mbox{if } s=0\mbox{ and } r=0,1,&\\
		&2z_1\cdots z_r& & \mbox{if } s=1,&\\
		&0& & \mbox{if } s>1.&
	\end{aligned}\right.
\end{align*}
\end{definition}
\noindent Then, $\xi$ is a quasi-derivation (see Section \ref{sec:kv}) ruled by $(-\!\odot\!-)$.

We need some notation. The inner Fox derivation $\rho_e$ associated with $e\in U(\hat L(H))$ is defined as
\[
	\rho_e(x, y) = (x - \varepsilon(x))e(y - \varepsilon(y)),
\]
while the quasi-derivation $q_{e_1,e_2}\colon U(\hat L(H)\oplus \mathbb{K}t_{**}) \to U(\hat L(H))$ associated with $e_1, e_2\in U(\hat L(H))$, which is ruled by $\rho_{e_1 + e_2}$, is defined as
\[
	q_{e_1, e_2}(\vec x) = (\varepsilon(x) - x)e_1 + e_2(\varepsilon(x) - x)
\]
where $x\in U(\hat L(H))$ is the image of $\vec x \in U(\hat L(H)\oplus \mathbb{K}t_{**})$ by the natural projection.\\

\begin{definition}\label{def:EN} 
Define the maps $E$ and $N$ so that the following diagrams are commutative:
\[\begin{tikzcd}[cramped]
	\widehat{\mathbb{K}\pi}\otimes \widehat{\mathbb{K}\pi} \arrow[r,"\eta"]\arrow[d,"Z_p\otimes Z_p"] & \widehat{\mathbb{K}\pi}\arrow[d,"Z_p"] & \widehat{\mathbb{K}\vec\pi} \arrow[r,"\vec\mu"]\arrow[d,"\vec Z_{p}"] & \widehat{\mathbb{K}\pi} \phantom{\,.} \arrow[d,"Z_p"]\\
	 U(\hat L(H))\otimes U(\hat L(H)) \arrow[r,"E"] & U(\hat L(H))& U(\hat L(H)\oplus \mathbb{K}t_{**}) \arrow[r,"N"] & U(\hat L(H))\,.
\end{tikzcd}\]
\end{definition}
\noindent Then, $E$ is a Fox pairing and $N$ is a quasi-derivation ruled by $E$, since $Z$ is an isomorphism of Hopf algebras.

\begin{lemma}\label{lem:E}
We have $E = (-\!\odot\!-) + \rho_{s(\omega)}$.
\end{lemma}
\noindent Proof. For $\alpha,\beta\in\pi$, we have $\iota^u(\beta) = (R^{v,u})^{-1}\Phi^\times(\beta_p)_{*,0}^{v,(u\;0)}(\Phi^\times)^{-1}R^{v,u}$ and hence
\begin{align*}
	E(Z_p(\alpha), Z_p(\beta)) &= Z_p(\eta(\alpha,\beta))\\
	&= \left(Z_p\circ p^v\circ \frac{\partial }{\partial z^v}\right) (\iota^u(\beta^{-1})\iota^v(\alpha)\iota^u(\beta))\\
	&= (D^v\circ Z_{p'})(\iota^u(\beta^{-1})\iota^v(\alpha)\iota^u(\beta))\\
	&= D^v\Big(e^{-t_{uv}/2} (\varphi^\times)^{-1} (Z_p(\beta^{-1})_{*,0}^{v,(u\;0)} \varphi^\times e^{t_{uv}/2}\\
	&\hspace{40pt} \varphi^{-1} (Z_p(\alpha)_{*,0}^{u,(v\;0)}\varphi e^{-t_{uv}/2} (\varphi^\times)^{-1} (Z_p(\beta)_{*,0}^{v,(u\;0)} \varphi^\times e^{t_{uv}/2}\Big).
\end{align*}
By the density argument, we have, for $x,y\in U(\hat L(H))$,
\begin{align*}
	E(x,y) &= D^v\Big(e^{-t_{uv}/2} (\varphi^\times)^{-1} (S(y')_{*,0}^{v,(u\;0)} \varphi^\times e^{t_{uv}/2}\\
	&\hspace{40pt} \varphi^{-1} x_{*,0}^{u,(v\;0)}\varphi e^{-t_{uv}/2} (\varphi^\times)^{-1} (y'')_{*,0}^{v,(u\;0)} \varphi^\times e^{t_{uv}/2}\Big).
\end{align*}
Here, we put $\Delta(y) = y'\otimes y''$. In particular, for $y\in H$, we have
\begin{align*}
	E(x,y) &= D^v([\varphi^{-1} x _{*,0}^{u,(v\;0)}\varphi, e^{-t_{uv}/2} (\varphi^\times)^{-1} y_{*,0}^{v,(u\;0)} \varphi^\times e^{t_{uv}/2}]).
\end{align*}
Now put
\begin{align*}
	U &= \varphi^{-1} x _{*,0}^{u,(v\;0)}\varphi - x _{*,0}^{u,(v\;0)} \mbox{ and}\\
	V &= e^{-t_{uv}/2} (\varphi^\times)^{-1} y_{*,0}^{v,(u\;0)} \varphi^\times e^{t_{uv}/2} - y_{*,0}^{v,(u\;0)}.
\end{align*}
Then, we have $U,V\in P^{u,v}$ and
\begin{align*}
	E(x,y) &= D^v([U + x _{*,0}^{u,(v\;0)}, V + y_{*,0}^{v,(u\;0)}]).
\end{align*}
We compute this by parts: for $x,y\in H$,
\begin{align*}
	D^v([U,V]) &= 0,\\
	D^v([U,y_{*,0}^{v,(u\;0)}]) &= S(D^u([U,y_{*,0}^{v,(u\;0)}])) = S(-y D^u(U)) = S(D^u(U)) y\\
	&= D^v(U) y = D^v(\varphi^{-1} x _{*,0}^{u,(v\;0)}\varphi - x _{*,0}^{u,(v\;0)}) y = x \phi y,\mbox{ and}\\
	D^v([x_{*,0}^{u,(v\;0)},V]) &= x D^v(V) = x S(D^u(V))\\
	&= x S(D^u(e^{-t_{uv}/2} (\varphi^\times)^{-1} y_{*,0}^{v,(u\;0)} \varphi^\times e^{t_{uv}/2} - y_{*,0}^{v,(u\;0)}))\\
	&= x S\left(y\Big(D^u(\varphi^\times) + \frac12\Big) \right) = -x S\left(D^u(\varphi^\times) + \frac12 \right) y\\
	&= -x \left( D^v(\varphi^\times) + \frac12 \right) y = -x \left( S(\phi) + \frac12 \right) y.
\end{align*}
Then, we have
\begin{align*}
	D^v([U,V]) + D^v([U,y_{*,0}^{v,(u\;0)}]) + D^v([x_{*,0}^{u,(v\;0)},V]) &= x \left( \phi - S(\phi) - \frac12 \right)y = x s(\omega)y.
\end{align*}
On the other hand, the values of $D^v([x_{*,0}^{u,(v\;0)},y_{*,0}^{v,(u\;0)}])$ is computed as follows. First of all, by the relations in $\mathfrak{t}_{g,12\cdots nuv0}$, we have\\
\centerline{\begin{tabular}{lll}
	$[x_u^i, x_v^k] = 0$, & $[x_u^i, y_v^k] = \delta_{ik}t_{uv}$, & $[x_u^i, t_{lv}] = 0$,\\
	$[y_u^i, x_v^k] = -\delta_{ik}t_{uv}$, & $[y_u^i, y_v^k] = 0$, & $[y_u^i, t_{lv}] = 0$,\\
	$[t_{ju}, x_v^k] = 0$, & $[t_{ju}, y_v^k] = 0$, & $[t_{ju}, t_{lv}] = -\delta_{jl}[t_{ju},t_{uv}]$.\\
\end{tabular}}
Then, we have\\
\centerline{\begin{tabular}{lll}
	$D^v([x_u^i, x_v^k]) = 0$, & $D^v([x_u^i, y_v^k]) = \delta_{ik}$, & $D^v([x_u^i, t_{lv}]) = 0$,\\
	$D^v([y_u^i, x_v^k]) = -\delta_{ik}$, & $D^v([y_u^i, y_v^k]) = 0$, & $D^v([y_u^i, t_{lv}]) = 0$,\\
	$D^v([t_{ju}, x_v^k]) = 0$, & $D^v([t_{ju}, y_v^k]) = 0$, & $D^v([t_{ju}, t_{lv}]) = -\delta_{jl}t_{j*}$,\\
\end{tabular}}
which implies $E = (-\!\odot\!-) + \rho_{s(\omega)}$.\qed\\

\begin{lemma}\label{lem:N}
We have $N = \xi + q_{\phi, -S(\phi) - \frac12}$.
\end{lemma}
\noindent Proof. Firstly, we have
\begin{align*}
	N(t_{**}) = N(e^{t_{**}}) = N(Z_p((F^*)^2)) = Z_p(\vec\mu((F^*)^2)) = 2.
\end{align*}
Next, we compute the value for other generators. For $\vec\alpha\in\vec\pi$, we have
\begin{align*}
	(N\circ \vec Z_p)(\vec\alpha) &= (Z_p\circ\vec\mu)(\vec\alpha)\\
	&= \left(Z_p\circ p^v\circ \frac{\partial }{\partial z^v}\right) (\iota^u(\alpha^{-1})c(\vec\alpha))\\
	&= (D^v\circ Z_{p'})(\iota^u(\alpha^{-1})c(\vec\alpha))\\
	&= D^v(e^{-t_{uv}/2}(\varphi^\times)^{-1} Z_p(\alpha^{-1})_{*,0}^{v,(u\;0)} \varphi^\times e^{t_{uv}/2} \cdot \vec Z_p(\vec\alpha)_*^{(uv)}).
\end{align*}
By the density argument, we have, for $\vec A\in U(\hat L(H)\oplus \mathbb{K}t_{**})$,
\begin{align*}
	N(\vec A) = D^v(e^{-t_{uv}/2}(\varphi^\times)^{-1} S(A')_{*,0}^{v,(u\;0)} \varphi^\times e^{t_{uv}/2} \cdot (\vec A'')_*^{(uv)}).
\end{align*}
Substituting $A \in H$, we have
\begin{align*}
	N(\vec A) &= D^v(-e^{-t_{uv}/2}(\varphi^\times)^{-1} A_{*,0}^{v,(u\;0)} \varphi^\times e^{t_{uv}/2} + (A_*^u + A_*^v))\\
	&= S(D^u(-e^{-t_{uv}/2}(\varphi^\times)^{-1} A_{*,0}^{v,(u\;0)} \varphi^\times e^{t_{uv}/2} + A_*^v)) + D^v(A_*^u)\\
	&= S(D^u(0 + A_*^v)) + D^v(A_*^u)\\
	&= S(-A\bar\phi) - A\phi\\
	&= \left(S(\phi) + \frac12\right)A - A\phi.
\end{align*}
This implies $N = \xi + q_{\phi, -S(\phi) - \frac12}$ since $N$ and $\xi + q_{\phi, -S(\phi) - \frac12}$ are both quasi-derivations ruled by the Fox pairing $E = (-\!\odot\!-) + \rho_{s(\omega)}$.\qed\\

Let $\varpi\colon U(L(H) \oplus \mathbb{K}t_{**}) \to U(L(H))$ be the natural projection and $f\colon U(L(H)) \to U(L(H) \oplus \mathbb{K}t_{**})$ a morphism of graded Hopf algebras being a section of $\varpi$. We define $f_j\in\mathbb{K}$ by $f(t_{j*}) = t_{j*} + f_j \frac{t_{**}}{2} =: z_{j*}^f$. Denote by $\vec\gamma_i\in\vec\pi$ a simple loop representing the $i$-th boundary $\partial_i\Sigma\, (1\leq i\leq n)$ with the induced orientation and the velocity vectors. 

\begin{lemma}\label{lem:deltagr}
Let $f$ and $f_i$ be as above. For any framing satisfying $\mathsf{fr}(\vec\gamma_i) + 1= f_i\, (1\leq i\leq n)$, we have
\[
	\delta_\xi \circ f(X) = \delta\gr^\mathsf{fr} (|X|)
\]
for $X\in U(\hat L(H))$.
\end{lemma}
\noindent Proof. Let $X = w_1\cdots w_r$ with $w_i\in H$. By Proposition 3.22 of \cite{akkn}, we have
\begin{align*}
	 \delta\gr^\mathsf{fr} (|w_1\cdots w_r|) &= \sum_{1\leq i\leq r} q^f(w_i) |1\wedge w_{i+1}\cdots w_{i-1}| + \sum_{1\leq i<j\leq r} |\tau\big(\tilde\Delta(w_i \odot w_j) (w_{i+1}\cdots w_{j-1} \otimes w_{j+1} \cdots w_{i-1})\big) |
\end{align*}
where $\tau$ is the anti-symmetrisation map, $\tilde\Delta = (S\otimes\id)\circ \Delta$ and $q^f\colon H \to \mathbb{K}$ is defined by
\begin{align*}
	q^f(x_*^a) = q^f(y_*^a ) = 0, \quad q^f(t_{*j}) = \mathsf{fr}(\vec\gamma_j) + 1 = f_j.
\end{align*}
On the other hand, we have
\begin{align*}
	f(w_1\cdots w_r) = w_1\cdots w_r + \sum_{1\leq i\leq r} \frac{q^f(w_i)}{2} w_1\cdots w_{i-1}w_{i+1}\cdots w_r t_{**} + Y t_{**}^2
\end{align*}
for some $Y \in U(L(H) \oplus \mathbb{K}t_{**})$. For the first term, analogous to the calculation in Section 8.5 of \cite{mas}, we have,  denoting $[a,b] = \{x\in \mathbb{Z}: a\leq x\leq b\}$,
\begin{align*}
	&\delta_\xi(w_1\cdots w_r)\\
	&= |\tau(\mu\otimes \id)(\id\otimes S\otimes \id)(\id \otimes \Delta) (\id \otimes \xi) \Delta (w_1\cdots w_r)|\\
	&= \sum_{I\sqcup J = [1,r]} |\tau(\mu\otimes \id)(\id\otimes S\otimes \id)(\id \otimes \Delta) (\id \otimes \xi) (w_I\otimes w_J)|\\
	&= \sum_{I\sqcup J = [1,r]} \sum_{J = \{J_1 < i < j < J_2\}} |\tau(\mu\otimes \id)(\id\otimes S\otimes \id)(\id \otimes \Delta) (w_I\otimes w_{J_1} (w_i \odot w_j) w_{J_2})|\\
	&= \sum_{I\sqcup J = [1,r]} \sum_{J = \{J_1 < i < j < J_2\}} |\tau(\mu\otimes \id)(\id\otimes S\otimes \id)(w_I\otimes w_{J_1}' (w_i \odot w_j)' w_{J_2}' \otimes w_{J_1}'' (w_i \odot w_j)'' w_{J_2}'')|\\
	&= \sum_{1\leq i<j\leq r} \sum_{I_1\sqcup J_1 = [1,i-1]} \sum_{I_2\sqcup J_2 = [j+1,r]} |\tau\big(w_{I_1}w_{[i+1, j-1]}w_{I_2} S(w_{J_1}' (w_i \odot w_j)' w_{J_2}' )\otimes w_{J_1}'' (w_i \odot w_j)'' w_{J_2}''\big)|\\
	&= \sum_{1\leq i<j\leq r} \sum_{I_1\sqcup J_1 = [1,i-1]} \sum_{I_2\sqcup J_2 = [j+1,r]} |\tau\big(w_{[i+1, j-1]}w_{I_2} S(w_{J_2}') S\big((w_i \odot w_j)'\big) S(w_{J_1}')w_{I_1} \otimes w_{J_1}'' (w_i \odot w_j)'' w_{J_2}''\big)|\\
	&= \sum_{1\leq i<j\leq r} |\tau\big(w_{[i+1, j-1]} S\big((w_i \odot w_j)'\big) \otimes w_{[1,i-1]} (w_i \odot w_j)'' w_{[j+1,r]}\big)|\\
	&= \sum_{1\leq i<j\leq r} |\tau\big(\tilde\Delta(w_i \odot w_j) (w_{i+1}\cdots w_{j-1} \otimes w_{j+1} \cdots w_{i-1})\big) |.
\end{align*}
As for the second term, we have
\begin{align*}
	&\delta_\xi\Big(\frac{q^f(w_i)}{2} w_1\cdots w_{i-1}w_{i+1}\cdots w_r t_{**} \Big)\\
	&= \frac{q^f(w_i)}{2} |\tau(\mu\otimes \id)(\id\otimes S\otimes \id)(\id \otimes \Delta) (\id \otimes \xi) \Delta( w_1\cdots w_{i-1}w_{i+1}\cdots w_r t_{**})|.
\end{align*}
If the first tensor component of $\Delta( w_1\cdots w_{i-1}w_{i+1}\cdots w_r t_{**}) \in U(L(H)\oplus \mathbb{K}t_{**})^{\otimes 2}$ has $t_{**}$ in the factor, the whole expression in the right-hand side vanishes due to the implicit projection $\varpi$, we only have to consider the contribution
\begin{align*}
	&\frac{q^f(w_i)}{2}  |\tau(\mu\otimes \id)(\id\otimes S\otimes \id)(\id \otimes \Delta) (\id \otimes \xi)\big( \Delta(w_1\cdots w_{i-1}w_{i+1}\cdots w_r)(1\otimes t_{**})\big)|\\
	&= q^f(w_i) |\tau(\mu\otimes \id)(\id\otimes S\otimes \id)(\id \otimes \Delta) \Delta(w_1\cdots w_{i-1}w_{i+1}\cdots w_r)|\\
	&= q^f(w_i) |\tau(1\otimes w_1\cdots w_{i-1}w_{i+1}\cdots w_r)|\\
	&= q^f(w_i) |1\wedge w_1\cdots w_{i-1}w_{i+1}\cdots w_r|.
\end{align*}
Tha last term $Y t_{**}^2$ does not contribute to the sum due to the projection $\varpi$, so we obtained $\delta_\xi \circ f(X) = \delta\gr^\mathsf{fr} (|X|)$ as claimed.\qed\\
 
We put $\mathbf{Ass}'_g = \{\vec Z \in \mathbf{Ass}_g: \mathrm{gr}(\vec Z_p) = \id, \mu = 1\}$. This set does not depend on the choice of $p$ or $n$ since $\vec Z$ is a morphism between (symmetric) operads. Now we can summarise the result into the following:

\begin{theorem}\label{thm:assoc}
Let $g,n\geq 0$ and $f$ as above. Then, we have a map $I^f_{g,n+1}\colon \mathbf{Ass}'_g \to \mathrm{Form}_{g,n+1}$.
\end{theorem}
\noindent Proof. The framing $\mathsf{fr}_{\vec Z, f}$ associated with a genus $g$ associator $\vec Z$ and the morphism $f$ is defined as the unique one making the diagram
\[\begin{tikzcd}[cramped]
	\widehat{\mathbb{K}\pi}\arrow[d,"Z_p"]\arrow[r,"\iota^{\mathsf{fr}_{\vec Z, f}}"] & \widehat{\mathbb{K}\vec\pi} \arrow[d,"\vec Z_{p}"]\\
	U(\hat L(H))\arrow[r, "f"]& U(\hat L(H)\oplus \mathbb{K}t_{**})
\end{tikzcd}\]
commute. This framing is described as follows: express $\vec Z_p(\vec \alpha)$ in terms of the new set of generators $x_*^a$, $y_*^a$, $z_{j*}^f$ and $t_{**}$. For $\vec\alpha\in\vec\pi$, the rotation number is computed as
\[
	\mathsf{fr}_{\vec Z, f}(\vec\alpha) = -2 \cdot (\mbox{the coefficient of }t_{**}\mbox{ in }\log\vec Z_p(\vec \alpha)).
\]

We shall check that $Z_p$ is the solution for the formality problem of the Goldman--Turaev Lie bialgebra with the $\mathbb{K}$-framing $\mathsf{fr}_{\vec Z, f}$. First of all, since we have assumed $\mathrm{gr}(\vec Z_p) = \id$, $Z_p$ is also an expansion. The equality $E = (-\!\odot\!-) + \rho_{s(\omega)}$ implies $Z_p$ preserves the Goldman bracket. Nextly, the commutative square involving $N$ in Definition \ref{def:EN} in the Lemma above induces another commutative square
\[\begin{tikzcd}[cramped]
	 \widehat{\mathbb{K}\vec\pi} \arrow[r,"\delta_{\vec\mu}"]\arrow[d,"\vec Z_{p}"] & \vert\widehat{\mathbb{K}\pi}\vert^{\otimes 2} \phantom{\,,} \arrow[d,"\vert Z_p\vert^{\otimes 2}"]\\
	 U(\hat L(H)\oplus \mathbb{K}t_{**}) \arrow[r,"\delta_N"] & \vert U(\hat L(H))\vert^{\otimes 2}\,,
\end{tikzcd}\]
where $\delta_{\vec\mu}$ is defined in Section \ref{sec:kv}, and $\delta_N$ is obtained similarly. By Equation (3.5) in \cite{mas}, we have 
\[
	\delta_{\vec\mu} \circ \iota^\mathsf{fr} (\alpha) = \delta^\mathsf{fr}(|\alpha|) + |1\wedge \alpha|\,.
\]

We now show that the adjoint action by $\mathfrak{t}_{12\cdots n*0}$ preserves the direct sum decomposition $\Imag f \oplus \mathbb{K}t_{**}$ of $\Ker\varepsilon^*$. For $i,j\in \{1,2,\dotsc, n, 0\}$, we have
\begin{align*}
	[t_{ij}, x_*^a] &= [t_{ij}, y_*^a] = 0,\mbox{ and}\\
	[t_{ij}, z^f_{k*}] &= [t_{ij}, t_{k*}]\\
	&= \delta_{ik} [t_{k*}, t_{j*}] + \delta_{jk} [t_{k*}, t_{i*}]\\
	&= \delta_{ik} [z^f_{k*}, z^f_{j*}] + \delta_{jk} [z^f_{k*}, z^f_{i*}].
\end{align*}
The action by $x_*^a$, $y_*^a$ or $t_{j*}$ obviously sends $\Imag f$ to itself, so $\Imag f$ is closed under the action. The action on $\mathbb{K}t_{**}$ is trivial, so the claim follows. 

The loop $\vec\gamma_i$ is a conjugate of $R^{i,*}R^{*,i}F^*$ by some product of $R^{j,k}$'s and $\Phi^{j,k,l}$'s, where $j,k,l \in \{1,2,\dotsc, n, *, 0\}$. This implies that $\vec Z_p(\vec \gamma_i)$ is a conjugate to $e^{t_{*j} + t_{**}/2}$ by an element of $\exp(\mathfrak{t}_{12\cdots n*0})$. Since the direct sum decomposition is preserved by the $\exp(\mathfrak{t}_{12\cdots n*0})$-conjugation, the coefficient of $t_{**}$ in $\log\vec Z_p(\vec \gamma_i)$ with respect to the generators $(x_*^a, y_*^a, z^f_{*j})$ only depends on the conjugacy class. Therefore, we have
\begin{align*}
	\mathsf{fr}_{\vec Z, f}(\vec\gamma_i) &= \mathsf{fr}_{\vec Z, f}(R^{i,*}R^{*,i} F^*)\\
	&= -2 \cdot (\mbox{the coefficient of }t_{**}\mbox{ in }\log (e^{z^f_{*i} - (f_i - 1)t_{**}/2}))\\
	&= f_i - 1.
\end{align*}
and hence
\[
	\delta_\xi \circ f(X) = \delta\gr^{\mathsf{fr}_{\vec Z, f}} (|X|)
\]
for $X\in U(\hat L(H))$ by the lemma above. In addition, we have
\[
	\delta_{q_{\phi, -S(\phi) - \frac12}}\circ f(X) = |XS(e')\otimes e'' + Xe''\otimes S(e') - S(e')\otimes e''X - e''\otimes S(e')X|
\]
by Lemma 2.5 of \cite{mas}. Here we put $e = s(\omega) = \frac{e^\omega}{1-e^\omega} + \frac1\omega$ and $\Delta(e) = e'\otimes e''$ is the coproduct of $e$. Since we can write as $s(\omega) = -\frac12 + (\mbox{an odd function in }\omega)$ and $\omega$ is primitive, we have
\[
	S(e')\otimes e'' + e''\otimes S(e') = -1\otimes 1
\]
and therefore
\[
	\delta_{q_{\phi, -S(\phi) - \frac12}}\circ f(X) = |-X\otimes 1 + 1\otimes X| = |1\wedge X|.
\]
Hence, we have,
\begin{align*}
	|Z_p|^{\otimes 2}\circ \delta^{\mathsf{fr}_{\vec Z,f}}(|\alpha|) &= |Z_p|^{\otimes 2}(\delta_{\vec\mu} \circ \iota^{\mathsf{fr}_{\vec Z,f}}(\alpha) - |1\wedge \alpha|)\\
	&= \delta_N\circ\vec Z_p\circ\iota^{\mathsf{fr}_{\vec Z,f}} (\alpha) - |1\wedge Z_p(\alpha)|\\
	&= \delta_N\circ f\circ Z_p (\alpha) - |1\wedge Z_p(\alpha)|\\
	&= \delta\gr^{\mathsf{fr}_{\vec Z,f}} (|Z_p (\alpha)|) + |1\wedge Z_p (\alpha)| - |1\wedge Z_p(\alpha)|\\
	&=  \delta\gr^{\mathsf{fr}_{\vec Z,f}} (|Z_p (\alpha)|).
\end{align*}
This shows that $|Z_p|$ also preserves the Turaev cobracket associated with $\mathsf{fr}_{\vec Z, f}$. Thus, we obtained the map $I_{g,n+1}^f(\vec Z) = (\mathsf{fr}_{\vec Z,f}, Z_p)$. \qed

\begin{remark}
The construction of the map $I_{g,n+1}^f$ also depends on the choice of $p\in\mathrm{Ob}(\mathbf{PaB}^f_g)$. This choice corresponds to the choice of a tree in Theorem 8.19 of \cite{akkn}.\\
\end{remark}

\section{The KV Equations}\label{sec:kveq}

In this section, we construct a solution to the KV equations in the sense of \cite{akkn} as a lift of the formality morphism obtained above. Namely, we will define the map
\[
	\tilde{I}_{g,n+1}^f\colon \mathbf{Ass}'_g \to \mathrm{SolKV}_{g,n+1}
\]
as a lift of $I_{g,n+1}^f$ along the natural map $\mathrm{SolKV}_{g,n+1} \to \mathrm{Form_{g,n+1}}$, connecting Gonzalez--Drinfeld associators to KV associators. Recall that the set $\mathrm{SolKV}^\mathsf{fr}_{g,n+1}$ is defined as
\[
	\mathrm{SolKV}^\mathsf{fr}_{g,n+1} = \Big\{\tilde G\in\tAut^+(\hat L(H)): G(\omega)=\xi,\,\mathsf{j}^\mathsf{fr}_\mathrm{gr}(\tilde G) -\mathbf{r} + \mathbf{p}^\mathsf{fr} \in \Big|\sum_j z_j\K{z_j} + \xi^2\K\xi\Big|\Big\}.
\]
We refer to the first equation $G(\omega)=\xi$ as (KVI), and the second as (KVII). A tangential automorphism is better understood as an isomorphism of Hopf groupoids in view of \cite{toyo4}, so we start with an embedding of the fundamental groupoid of the surface to $\mathbf{PaB}^f_g$.\\

The tangential base points in Figure \ref{fig:points} define the set $V = \{*_i\}_{0\leq i\leq n} \subset \partial\Sigma$ of the base points of $\Sigma$, where the point $*_i$ is on the $i$-th boundary component $\partial_i\Sigma$. Now set $\mathscr{G} = \pi_1(\Sigma, V)$, the fundamental groupoid of $\Sigma$ with the base points $V$, which is a free groupoid, and $\vec{\mathscr{G}} = \pi_1(U\Sigma, V)$ with the tangent vectors at $V$ fixed once and for all. We identify $\vec\pi$ with the endomorphism group $\vec{\mathscr{G}}(*_0, *_0)$.

We have a Hopf groupoid $\mathbb{K}\mathscr{G}$ and its completion $\widehat{\mathbb{K}\mathscr{G}}$ with respect to the multiplicative filtration given, for $X\in\mathbb{K}\mathscr{G}(*_i, *_j)$, by
\[
	\mathrm{wt}(X) := \mathrm{wt}(\alpha X\beta)
\]
for some $\alpha\in\mathscr{G}(*_0, *_i)$ and $\beta\in\mathscr{G}(*_j, *_0)$, using the weight filtration on $\mathbb{K}\pi$ in Section \ref{sec:kv}. Since every invertible element has weight $0$ under any multiplicative filtration (parametrised by non-negative integers), this is well-defined. 

Now we construct an embedding of $\vec{\mathscr{G}}$ into $\mathbf{PaB}^f_g(n+2)$. Recall that we have fixed an object $p = (\ell_n)(* \;0)$. We define the object $p_i \in \mathrm{Ob}(\mathbf{PaB}^f_g)$ for $0\leq i\leq n$ by
\[
	p_i = \varepsilon^*(p) \circ_i (*\;i)
\]
where $\varepsilon$ is the deletion of $*$ and $\circ_i$ denotes the operadic composition to the letter $i$. We have $p_0 = p$ by definition. Set $P = \{p_i\}_{0\leq i\leq n}$. Next, we consider the subgroupoid $\mathscr{S}$ of $\mathbf{PaB}^f_g(n+1)$ with the sole object $\varepsilon^*(p)$ and the only morphism being $\id_{\varepsilon^*(p)}$.

\begin{lemma}
The full subgroupoid $(\varepsilon^*)^{-1}(\mathscr{S})|_P$ of the fibre of $\mathscr{S}$ along the groupoid homomorphism
\[
	\varepsilon^* \colon \mathbf{PaB}^f_g(n+2) \to \mathbf{PaB}^f_g(n+1)
\]
is naturally isomorphic to $\vec{\mathscr{G}}$ by identifying $*_i \in \mathrm{Ob}(\vec{\mathscr{G}})$ with $p_i\in \mathrm{Ob}(\mathbf{PaB}^f_g)$.
\end{lemma}
\noindent Proof. The map $\vec{\mathscr{G}} \to \mathbf{PaB}^f_g(n+2)$ of groupoids is induced by the natural embedding $\Sigma_{g,n+1} \to \Sigma_g$. Since we identify $\mathrm{Ob}(\vec{\mathscr{G}})$ with $P$, and both $\vec{\mathscr{G}}$ and $(\varepsilon^*)^{-1}(\mathscr{S})|_P$ are groupoids such that every hom-set is non-empty, we only check the isomorphism on one object, which we take to be $p$. Then, it is equivalent to the exactness of the sequence
\[\begin{tikzcd}
	1 \arrow[r]& \vec\pi \arrow[r, "\iota"]& \mathrm{PB}^f_{g,12\cdots n*0} \arrow[r, "\varepsilon^*"]& \mathrm{PB}^f_{g,12\cdots n0} \arrow[r]& 1,
\end{tikzcd}\]
whose proof will be given in Lemma \ref{lem:exactgrp}.\qed\\

Similarly, we have the groupoid $\mathscr{S}'$ of the group-like part $\mathbb{G}(\mathbf{PaCD}^f_g)$ corresponding to $\mathscr{S}$. The morphism set $(\varepsilon^*)^{-1}(\mathscr{S}')(p_i,p_j)$ is isomorphic to $\mathrm{exp}(\hat L(H)\oplus \mathbb{K}t_{**})$ by Lemma \ref{lem:exactness}, so we have $(\varepsilon^*)^{-1}(\mathscr{S}')|_P = \mathbb{G}(\mathrm{gr}\,\widehat{\mathbb{K}\vec{\mathscr{G}}})$.

\begin{definition}
A \textit{special tangential automorphism} $F$ is an element of $\mathrm{Isom}^+_\partial(\widehat{\mathbb{K}\mathscr{G}}, \mathrm{gr}\,\widehat{\mathbb{K}\mathscr{G}})$, the set of continuous isomorphisms of complete Hopf groupoids with $\mathrm{gr}\,F = \id$ and preserving boundary: $\partial\widehat{\mathbb{K}\mathscr{G}} \xrightarrow{\cong} \mathrm{gr}\,\partial\widehat{\mathbb{K}\mathscr{G}}$. Here, $\partial\widehat{\mathbb{K}\mathscr{G}}$ is the Hopf subgroupoid of $\widehat{\mathbb{K}\mathscr{G}}$ (topologically) generated by boundary loops $\partial_i\Sigma\in \mathscr{G}(*_i, *_i)$.
\end{definition}

The Hopf subgroupoid $\mathrm{gr}\,\partial\widehat{\mathbb{K}\mathscr{G}}$ is generated by $t_{i*} \in U(\mathfrak{t}_{g,n+2}) = \varpi(\mathbf{PaCD}^f_g(p_i,p_i))$ for $0\leq i\leq n$. Therefore, the boundary-preserving condition amounts to the equality $F(\partial_i\Sigma) = e^{t_{i*}}$ in $\mathfrak{t}_{g,n+2}$ under the condition $\mathrm{gr}\,F = \id$ (otherwise the coefficient of $t_{i*}$ would not be $1$). 

Once we fix a set of paths $\delta_j\in\mathscr{G}(*_j, *_0)$ for $1\leq j\leq n$ such that $\delta_j\gamma_j\delta_j^{-1} = \partial_j\Sigma$, we have a bijection
\begin{align*}
	J\colon \mathrm{Isom}^+_\partial(\widehat{\mathbb{K}\mathscr{G}}, \mathrm{gr}\,\widehat{\mathbb{K}\mathscr{G}}) &\to \{ \tilde G\in \tAut^+(\hat L(H)): G(\omega) =\xi \}\\
	F &\mapsto \big(F|_{*_0}\circ \theta_\mathrm{exp}^{-1}; F(\delta_1), \dotsc, F(\delta_n) \big)^{-1}
\end{align*}
For more details, see Section 4 of $\cite{toyo4}$.

\begin{lemma}
A genus $g$ Gonzalez--Drinfeld associator $\vec Z\in\mathbf{Ass}'_g$ induces a special tangential automorphism $Z\colon \widehat{\mathbb{K}\mathscr{G}}\to \mathrm{gr}\,\widehat{\mathbb{K}\mathscr{G}}$.
\end{lemma}
\noindent Proof. Since $\vec Z(\mathscr{S}) = \mathscr{S}'$ and the operad isomorphism $\vec Z$ was assumed to preserve objects, we have
\[
	\vec Z (\vec{\mathscr{G}}) = \vec Z((\varepsilon^*)^{-1}(\mathscr{S})|_P) = (\varepsilon^*)^{-1}(\mathscr{S}')|_P = \mathbb{G}(\mathrm{gr}\,\widehat{\mathbb{K}\vec{\mathscr{G}}}).
\]
This induces the isomorphism $Z\colon \widehat{\mathbb{K}\mathscr{G}}\to \mathrm{gr}\,\widehat{\mathbb{K}\mathscr{G}}$ of the quotient spaces.

We check the property $\mathrm{gr}(Z) = \id$. For $X\in\widehat{\mathbb{K}\mathscr{G}}(*_i, *_j)$, take any $\alpha\in\mathscr{G}(*_0, *_i)$ and $\beta\in\mathscr{G}(*_j, *_0)$ so that
\[
	\alpha X \beta \in \widehat{\mathbb{K}\mathscr{G}}(*_0, *_0) = \widehat{\mathbb{K}\pi}.
\]
Since we take $\vec Z$ from $\mathbf{Ass}'_g$, whose element satisfies $\mathrm{gr}(Z) = \id$ on $\mathbb{K}\pi$, we have, modulo weight $\geq (\mathrm{wt}(X)+1)$-part,
\[
	\alpha X\beta \equiv Z(\alpha X\beta) =  Z(\alpha) Z(X) Z(\beta)
\]
and therefore
\[
	Z(X) \equiv Z(\alpha^{-1})\alpha X\beta Z(\beta^{-1}).
\]
Since multiplication by an invertible element preserves filtration, we conclude $Z(X) \equiv X$. This shows $\mathrm{gr}(Z) = \id$ on the whole $\widehat{\mathbb{K}\mathscr{G}}$.

Lastly, we check the boundary-preserving condition. Since in the object $p_i$, the letters $i$ and $*$ are placed inside the same innermost parenthesis, we can take the lift of the boundary loop $\partial_i\Sigma\in\mathscr{G}(*_i,*_i)$ as $R^{*i}R^{i*}\in\widehat{\mathbb{K}\mathbf{PaB}}{}^f_g(p_i,p_i)$. Therefore, we have
\[
	Z(\partial_i\Sigma) = \varpi(\vec Z(R^{*i}R^{i*})) = e^{t_{i*}},
\]
which is exactly the boundary-preserving condition. This completes the proof.\qed

\begin{theorem}\label{thm:solkv}
The above construction defines the map $\tilde{I}_{g,n+1}^f\colon \mathbf{Ass}'_g \to \mathrm{SolKV}_{g,n+1}$ by setting
\[
	\tilde{I}_{g,n+1}(\vec Z) = (\mathsf{fr}_{\vec Z, f}, J(Z))
\]
which is a lift of $I_{g,n+1}^f\colon \mathbf{Ass}'_g \to \mathrm{Form}_{g,n+1}$.
\end{theorem}

Before the proof, we need the following:
\begin{lemma}\label{lem:cocycj}
Let $\tilde F \in \tAut^+(\hat L(H))$ and suppose that the composition $\theta = F\circ \theta_\mathrm{exp}$ is
\begin{enumerate}[(1)]
	\item a solution to the formality problem with respect to a $\mathbb{K}$-framing $\mathsf{fr}$, and
	\item $\theta$ is special: $\theta(\partial_0\Sigma) = e^{\omega}$.
\end{enumerate}
Then, we have $\tilde F^{-1} \in \mathrm{SolKV}_{g,n+1}^\mathsf{fr}$.
\end{lemma}
\noindent Proof. Since $\theta$ is a formality morphism for $\mathsf{fr}$, we can take a lift $\tilde G\in\tAut^+(\hat L(H))$ of $F$ (meaning $F=G$) and $g\in \mathrm{exp}(\hat L(H))$ such that $(\tilde G \circ \mathrm{Ad}_g)^{-1} \in \mathrm{SolKV}^\mathsf{fr}_{g,n+1}$ by Theorem \ref{thm:akkn} (the reader may consult the diagram below Theorem \ref{thm:akkn}). By assumption (2), we have $F(\xi) = \omega$, which is equivalent to the equation (KVI). Therefore, it remains to verify (KVII); we shall show the stronger statement
\[
	\mathsf{j}^\mathsf{fr}\gr(\tilde F^{-1}) \equiv \mathsf{j}^\mathsf{fr}\gr((\tilde G \circ \mathrm{Ad}_g)^{-1}),
\]
where $\equiv$ denotes an equality modulo $\Big|\sum_j z_j\K{z_j}\Big|$. By (KVI) for $(\tilde G \circ \mathrm{Ad}_g)^{-1}$, we have
\[
	\omega = G \circ \mathrm{Ad}_g(\xi) = G(g \xi g^{-1}) = F(g \xi g^{-1}) = F(g) \omega F(g^{-1}) = G(g) \omega G(g^{-1})
\]
in $\hat L(H)$. Therefore, $\omega$ and $G(g)$ commute and we have $G(g) = e^{-\lambda\omega}$ for some $\lambda\in \mathbb{K}$. We put $\tilde u= -\log \tilde F$ and $\tilde v = -\log \tilde G$, and denote by $u$ and $v$ the underlying derivations on $\hat L(H)$, and by $(u_i)_{1\leq i\leq n}$ and $(v_i)_{1\leq i\leq n}$ their tangential coefficients. Since we have $F = G$, we also have $u = v$ and hence
\[
	[z_i, u_i] = u(z_i) = v(z_i) = [z_i, v_i]
\]
in $\hat L(H)$ so that $u_i - v_i = \lambda_i z_i$ for some $\lambda_i \in \mathbb{K}$. This implies $|u_i| \equiv |v_i|$ on the cyclic quotient $|\hat T(H)|$ and hence
\[
	\mathsf{b}^\mathsf{fr}(\tilde u) = \sum_i c_i|u_i| \equiv \sum_i c_i|v_i| = \mathsf{b}^\mathsf{fr}(\tilde v).
\]
Since $\tilde u$ and $\tilde v$ are tangential derivations, the action on $\Big|\sum_j z_j\K{z_j}\Big|$ is trivial and therefore preserves the equivalence $\equiv$. This implies
\[
	\mathsf{c}^\mathsf{fr}(\tilde F^{-1}) = \frac{e^u - 1}{u}\cdot \mathsf{b}^\mathsf{fr}(\tilde u)  \equiv \frac{e^v - 1}{v}\cdot \mathsf{b}^\mathsf{fr}(\tilde v)  = \mathsf{c}^\mathsf{fr}(\tilde G^{-1}),
\]
using $u = v$.

Since the action of $\mathrm{Ad}_{e^{\lambda\omega}}$ is identity on $|T(H)|$, we have
\begin{align*}
	\mathsf{j}_{x,y,z}(\mathrm{Ad}_{e^{\lambda\omega}}) &= \frac{e^{\mathrm{ad}_{\lambda\omega}} -1}{\mathrm{ad}_{\lambda\omega}}\cdot \mathsf{div}_{x,y,z}(\mathrm{ad}_{\lambda\omega})\\
	&= \mathsf{div}_{x,y,z}(\mathrm{ad}_{\lambda\omega})\\
	&= \lambda\, \Big| \sum_{1\leq a\leq g} d_{x_a}[\omega, x_a] + d_{y_a}[\omega, y_a] + \sum_{1\leq i\leq g} d_{z_i}[\omega, z_i] \Big|\\
	&= \lambda\, \Big| \sum_{1\leq a\leq g} (\omega  - x_a d_{x_a}(\omega)) + (\omega  - y_a d_{y_a}(\omega))+ \sum_{1\leq i\leq g} (\omega  - z_i d_{z_i}(\omega))\Big|\\
	&\equiv -\lambda\, \Big| \sum_{1\leq a\leq g} (x_a d_{x_a}(\omega) + y_a d_{y_a}(\omega)) + \sum_{1\leq i\leq g} z_i d_{z_i}(\omega)\Big|\\
	&\equiv -\lambda\, \Big| \sum_{1\leq a\leq g} (x_a (-y_a) + y_a x_a) + \sum_{1\leq i\leq g} z_i \cdot 1\Big|\\
	&\equiv 0
\end{align*}
using $|\omega| \equiv 0$. In addition, we have
\begin{align*}
	 \mathsf{c}^\mathsf{fr}(\mathrm{Ad}_{e^{\lambda\omega}}) =  \mathsf{b}^\mathsf{fr}(\mathrm{ad}_{\lambda\omega}) = \sum_i c_i |\lambda\omega| \equiv 0.
\end{align*}
Combining the above with the fact that $\mathsf{j}^\mathsf{fr}\gr$ is a group $1$-cocycle, we have
\begin{align*}
	\mathsf{j}^\mathsf{fr}\gr((\tilde G \circ \mathrm{Ad}_g)^{-1}) &= \mathsf{j}^\mathsf{fr}\gr((\mathrm{Ad}_{G(g)}\circ \tilde G)^{-1})\\
	&= \mathsf{j}^\mathsf{fr}\gr(\tilde G^{-1} \circ \mathrm{Ad}_{G(g^{-1})})\\
	&= \mathsf{j}^\mathsf{fr}\gr(\tilde G^{-1}) + G^{-1}\cdot \mathsf{j}^\mathsf{fr}\gr(\mathrm{Ad}_{e^{\lambda\omega}})\\
	&= \big( \mathsf{j}_{x,y,z}(G^{-1}) - \mathsf{c}^\mathsf{fr}(\tilde G^{-1})\big) + G^{-1}\cdot \big( \mathsf{j}_{x,y,z}(\mathrm{Ad}_{e^{\lambda\omega}}) - \mathsf{c}^\mathsf{fr}(\mathrm{Ad}_{e^{\lambda\omega}}) \big)\\
	&\equiv  \big( \mathsf{j}_{x,y,z}(F^{-1}) - \mathsf{c}^\mathsf{fr}(\tilde F^{-1})\big) + G^{-1}\cdot ( 0-0)\\
	&= \mathsf{j}^\mathsf{fr}\gr(\tilde F^{-1}) .
\end{align*}
This shows $\tilde{F}^{-1}$ satisfies (KVII) associated with $\mathsf{fr}$; this completes the proof.\qed\\

\noindent \textbf{Proof of Theorem \ref{thm:solkv}.} First of all, we show $Z \in \mathrm{SolKV}_{g,n+1}^{\mathsf{fr}_{\vec Z, f}}$. By Theorem \ref{thm:assoc}, $\theta := Z_p$ is a formality morphism with respect to $\mathsf{fr}_{\vec Z, f}$. We also have $\theta(\gamma_0) = \theta(\partial_0\Sigma) =e^\omega$ by Lemma \ref{lem:special}. Now we can apply Lemma \ref{lem:cocycj} for $\tilde F = J(Z)^{-1}$ to conclude $J(Z) \in \mathrm{SolKV}_{g,n+1}^{\mathsf{fr}_{\vec Z, f}}$.

Since the natural map $\mathrm{SolKV}_{g,n+1} \to \mathrm{Form}_{g,n+1}$ sends $J(Z)$ to its restriction $Z_p$ to the object $p = p_0 = {*_0}$, $\tilde{I}_{g,n+1}^f$ is a lift of $I_{g,n+1}^f$ as claimed. This completes the proof.\qed\\

Finally, we deal with the Grothendieck--Teichm\"uller groups. Recall that the group $\widehat{\mathbf{GT}}_g$ is defined as the automorphism group of the pair $(\widehat{\mathbb{K}\mathbf{PaB}}{}^f, \widehat{\mathbb{K}\mathbf{PaB}}{}^f_g)$. We define $\widehat{\mathbf{GT}}{}'_g$ as the subgroup consisting of automorphisms such that the coupling constant is $1$ and the associated graded map is the identity map. We will show that this acts on $\mathrm{SolKV}_{g,n+1}$ via the map $\tilde I_{g,n+1}^f$.

\begin{deflem}
For $\vec G\in \widehat{\mathbf{GT}}{}'_g$, we define $h_{\vec G}\in H^1(\Sigma_{g,0};\mathbb{K})$ as the unique element making the diagram
\[\begin{tikzcd}[cramped]
	\widehat{\mathbb{K}\pi}\arrow[d,"G_p"]\arrow[r,"\iota^{\mathsf{fr} + h_{\vec G}}"] & \widehat{\mathbb{K}\vec\pi} \arrow[d,"\vec G_p"]\\
	\widehat{\mathbb{K}\pi} \arrow[r, "\iota^\mathsf{fr}"]& \widehat{\mathbb{K}\vec\pi}
\end{tikzcd}\]
commutes, where $\mathsf{fr}$ is any $\mathbb{K}$-framing and we regard $h_{\vec G}$ as an element of $H^1(\Sigma_{g,n+1};\mathbb{K})$ via the pull-back by the inclusion $\Sigma_{g,n+1} \to \Sigma_{g,0}$.
\end{deflem}
\noindent Proof. Since the coupling constant is $1$, we have $\vec G(R^{j,*}R^{*,j}) = R^{j,*}R^{*,j}$ and $\vec G(F^*) = F^*$; the former shows $h_{\vec G}(\gamma_j) = 0$ and therefore $h_{\vec G}\in H^1(\Sigma_{g,0};\mathbb{K})$. Now take $h\in H^1(\Sigma;\mathbb{K}) \cong \Hom_\mathrm{grp}(\pi,\mathbb{K})$ so that $\mathsf{fr} + h$ is another $\mathbb{K}$-framing. Then, we have
\begin{align*}
	\iota^{\mathsf{fr} + h}(\alpha) = \iota^\mathsf{fr}(\alpha) (F^*)^{h(\alpha)}
\end{align*}
for $\alpha\in\pi$, so that
\begin{align*}
	\iota^{\mathsf{fr} + h}(G_p(\alpha)) &= \iota^\mathsf{fr}(G_p(\alpha)) (F^*)^{h(G_p(\alpha))}\\
	&= \vec G_p\big(\iota^{\mathsf{fr} + h_{\vec G}}(\alpha) (F^*)^{h(G_p(\alpha)}\big).
\end{align*}
If $h(\alpha) = h(G_p(\alpha))$ holds, the right-hand side is equal to $\vec G_p\big(\iota^{\mathsf{fr} + h + h_{\vec G}}(\alpha)\big))$ so we are done. To prove $h(\alpha) = h(G_p(\alpha))$, we use the assumption $\mathrm{gr}(\vec G) = \id$ to conclude that the induced map on the abelianisation $H = \pi^\mathrm{ab}$ by $G_p$ is also the identity. Since $h$ factors through the quotient map $\hat\pi \to H$, we have $h(\alpha) = h(G_p(\alpha))$ as claimed. This shows that $h_{\vec G}$ does not depend on the choice of $\mathsf{fr}$.\qed

\begin{proposition}
Assuming the existence of a genus $g$ Gonzalez--Drinfeld associator, $\widehat{\mathbf{GT}}{}'_g$ acts on $\mathrm{SolKV}_{g,n+1}$ from right by
\[
	(\mathsf{fr}, J(F))\cdot \vec G = (\mathsf{fr} + h_{\vec G}, J(F\circ G))\,.
\]
Here, $G$ is the induced isomorphism on $\widehat{\mathbb{K}\mathscr{G}}$ from $\vec G$ on the quotient.
\end{proposition}
\noindent Proof. To show $(\mathsf{fr} + h_{\vec G}, J(F\circ G)) \in \mathrm{SolKV}_{g,n+1}$, it is enough to prove that: (1) $F\circ G \in \mathrm{Isom}^+_\partial(\widehat{\mathbb{K}\mathscr{G}}, \mathrm{gr}\,\widehat{\mathbb{K}\mathscr{G}})$ and (2) $F\circ G$ induces a formality morphism for $\mathsf{fr} + h_{\vec G}$, by Lemma \ref{lem:cocycj}.

(1) follows from
\[
	G(\partial_j\Sigma) = G(R^{j,*}R^{*,j}) = R^{j,*}R^{*,j} = \partial_j\Sigma
\]
for $0\leq j\leq n$. For (2), we first show that $\vec G$ induces the isomorphism
\begin{align}\label{eq:Gp}
	G_p \colon (|\widehat{\mathbb{K}\pi}|, \mathsf{fr} +  h_{\vec G}) \to (|\widehat{\mathbb{K}\pi}|, \mathsf{fr})
\end{align}
of Lie bialgebras for any $\mathsf{fr}$. Take $\vec Z \in \mathbf{Ass}_g$. Then, we know that $\vec Z\circ \vec G\in \mathbf{Ass}_g$ by definition. We apply Theorem \ref{thm:assoc} to both $\vec Z$ and $\vec Z\circ \vec G$ and the same $f$ to obtain the Lie bialgebra isomorphism
\[
	G_p \colon (|\widehat{\mathbb{K}\pi}|, \mathsf{fr}_{\vec Z\circ\vec G, f}) \to (|\widehat{\mathbb{K}\pi}|, \mathsf{fr}_{\vec Z, f}).
\]
We have $\mathsf{fr}_{\vec Z\circ\vec G, f} =  \mathsf{fr}_{\vec Z, f} + h_{\vec G}$ by the definition of $h_{\vec G}$.

We now show that we can replace $\mathsf{fr}_{\vec Z, f}$ with an arbitrary $\mathbb{K}$-framing. Let $h\in H^1(\Sigma;\mathbb{K}) \cong \Hom_\mathrm{grp}(\pi,\mathbb{K})$. We have 
\[
	\delta^{\mathsf{fr} + h} (\alpha) = \delta^\mathsf{fr}(\alpha) + h(\alpha) |1\wedge \alpha| 
\]
for $\alpha\in\pi$ and any $\mathbb{K}$-framing $\mathsf{fr}$. The map $h$ naturally extends to the Malcev completion, so we obtain $h\colon \hat\pi \to \mathbb{K}$. Substituting $\mathsf{fr} = \mathsf{fr}_{\vec Z\circ\vec G, f}$ and applying $G_p$ to this equality, we have
\begin{align*}
	G_p(\delta^{ \mathsf{fr}_{\vec Z\circ\vec G, f} + h} (\alpha)) &=  G_p(\delta^{ \mathsf{fr}_{\vec Z\circ\vec G, f}}(\alpha) + h(\alpha) |1\wedge \alpha| )\\
	&= \delta^{ \mathsf{fr}_{\vec Z, f}}(G_p(\alpha)) + h(\alpha) |1\wedge G_p(\alpha)| 
\end{align*}
Combining with $h(\alpha) = h(G_p(\alpha))$, the right-hand side is equal to $\delta^{ \mathsf{fr}_{\vec Z, f} + h}(G_p(\alpha))$, so we obtain the isomorphism $\eqref{eq:Gp}$ by setting $h = \mathsf{fr} - \mathsf{fr}_{\vec Z, f} $.

At this point, we have the sequence of isomorphisms
\[
	(|\widehat{\mathbb{K}\pi}|, \mathsf{fr}_{\vec Z\circ\vec G, f}) \xrightarrow{G_p} (|\widehat{\mathbb{K}\pi}|, \mathsf{fr}_{\vec Z, f}) \xrightarrow{F_p} (\mathrm{gr}|\widehat{\mathbb{K}\pi}|, \mathsf{fr}_{\vec Z, f}).
\]
Since $\delta_\mathrm{gr}^\mathsf{fr}$ only depends on the values of $\mathsf{fr}(\gamma_j)$ as can be seen in the proof of Lemma \ref{lem:deltagr}, $(\mathrm{gr}|\widehat{\mathbb{K}\pi}|, \mathsf{fr}_{\vec Z, f})$ is equal to $(\mathrm{gr}|\widehat{\mathbb{K}\pi}|, \mathsf{fr}_{\vec Z\circ\vec G, f})$ using $h_{\vec G}\in H^1(\Sigma_{g,0};\mathbb{K})$. This completes the proof.\qed

\begin{theorem}\label{thm:equivariant}
The map $\tilde I_{g,n+1}^f$ is $\widehat{\mathbf{GT}}{}'_g$-equivariant.
\end{theorem}
\noindent Proof. This is straightforward from the construction.\qed\\

Similarly, we have a subgroup $\mathbf{GRT}{}'_g$ of $\mathbf{GRT}_g$ and one can show that there is a natural $\mathbf{GRT}{}'_g$-action under which the map $\tilde I_{g,n+1}^f$ is equivariant, without the assumption on the existence of the genus $g$ Gonzalez--Drinfeld associator.

\begin{remark}
By the proposition above, we have the group homomorphism
\[
	\widehat{\mathbf{GT}}{}'_g\cap \{h_{\vec G} = 0\} \to \mathrm{KV}^\mathsf{fr}_{g,n+1}
\]
for each framing $\mathsf{fr}$ such that $\mathrm{Form}^\mathsf{fr}_{g,n+1}$ is non-empty. The domain can be seen as an analogue of the Chillingworth subgroup (see the original paper \cite{chillingworth}) of the mapping class group of a surface since such $\vec G$ fixes every $\mathbb{K}$-framing.
\end{remark}

As mentioned in the introduction, the Alekseev--Torossian map
\[
	\{\mbox{Drinfeld associators}\} \to \{\mbox{solutions to the KV equations}\}
\]
is injective. In this virtue, we pose the following:
\begin{problem}\label{problem:inj}
Describe the image of the map $\{\tilde{I}_{g,n+1}^f\}_{n\geq 0}\colon \mathbf{Ass}'_g \to \prod_{n\geq 0}\mathrm{SolKV}_{g,n+1}$.\\
\end{problem}

\section{Computation on the Associated Framing}\label{sec:framing}
In this section, we compute the associated framing $\mathsf{fr}_{\vec Z, f}$ in terms of the coefficients appearing in $\vec Z$ to see which framing appears. We begin with the case of $n=0$ so that $\mathfrak{t}_{g,12\cdots n*0}^f = \mathfrak{t}_{g,2}^f$, where we have only one choice for $f\colon U(L(H)) \to U(L(H) \oplus \mathbb{K}t_{**})$, namely, $f(t_{*0}) = t_{*0} + (1-g)t_{**}$.

Consider the relation $(\mathrm{D}_g)$ in \cite{gonzalez}:
\[
	C_a^{(12)3,\varnothing} = \Phi^{1,2,3} C_a^{1,23}R^{1,23} \Phi^{2,3,1} C_a^{2,31}R^{2,31} \Phi^{3,1,2} C_a^{3,12}R^{3,12}
\]
for $C = A,B$ and $1\leq a\leq g$. By removing the third strand, it is reduced to 
\begin{align}
	C_a^{(12),\varnothing} = C_a^{1,2}R^{1,2} C_a^{2,1}R^{2,1}. \label{eq:Dg_red}
\end{align}
Let $\vec Z$ be a genus $g$ Gonzalez--Drinfeld associator. The relation $A^{\varnothing,1}_a = \id$ says
\[
	\vec Z(A^{1,2}_a) \in \Ker(\exp(\mathfrak{t}^f_{g,2}) \xrightarrow{\varepsilon^1} \exp(\mathfrak{t}^f_{g,1})) \cong \exp( \hat L(x_1^a, y_1^a)_{1\leq a\leq g} \oplus \mathbb{K}t_{11}),
\]
so we can uniquely write
\begin{align*}
	\vec Z(A^{1,2}_a) = \exp(\xi^a_1+ s^a \cdot t_{11})
\end{align*}
with $\xi^a_1 \in \hat L(x_1^a, y_1^a)_{1\leq a\leq g}$ and $s^a\in\mathbb{K}$. Put
\begin{align*}
	\xi_1^a &= \lambda^a_b x_1^b + \mu^a_b y_1^b +(\nu^{xx})^a_{bc} [x_1^b,x_1^c]  + (\nu^{xy})^a_{bc} [x_1^b,y_1^c] + (\nu^{yy})^a_{bc} [y_1^b,y_1^c]\\
	&\qquad + (\pi^{xxx})^a_{bcd} [x_1^b,[x_1^c,x_1^d]] + (\pi^{xxy})^a_{bcd} [x_1^b,[x_1^c,y_1^d]] + (\pi^{yyx})^a_{bcd} [y_1^b,[y_1^c,x_1^d]] + (\pi^{yyy})^a_{bcd} [y_1^b,[y_1^c,y_1^d]]\\
	&\qquad + (\mbox{tot.deg.}\geq 4)
\end{align*}
where the coefficients $\lambda^a_b, \mu^a_b, \dotsc$ are elements of $\mathbb{K}$ and summations are implicit. By applying $\vec Z$ to \eqref{eq:Dg_red} with $C=A$, we have
\[
	\exp(\xi^a_{12} + s^a\cdot t_{12,12}) = \exp(\xi^a_1 + s^a\cdot t_{11})\exp(\frac{t_{12}}2) \exp(\xi^a_2 + s^a\cdot t_{22})\exp(\frac{t_{12}}2) \quad\mbox{in } \exp(\mathfrak{t}^f_{g,2}).
\]
Taking the logarithm, we have
\[
	\xi^a_{12} + (2s^a-\frac12)\cdot t_{12} = \mathrm{bch}(\xi^a_1,\frac{t_{12}}2,\xi^a_2).
\]
Skipping a tedious calculation, the right-hand side yields
\begin{align*}
	&\mathrm{bch}(\xi^a_1,\frac{t_{12}}2,\xi^a_2)\\
	&= \frac12 t_{12} + \lambda^a_b x_1^b + \mu^a_b y_1^b +(\nu^{xx})^a_{bc} [x_1^b,x_1^c]  + (\nu^{xy})^a_{bc} [x_1^b,y_1^c] + (\nu^{yy})^a_{bc} [y_1^b,y_1^c]\\
	&\qquad + (\pi^{xxx})^a_{bcd} [x_1^b,[x_1^c,x_1^d]] + (\pi^{xxy})^a_{bcd} [x_1^b,[x_1^c,y_1^d]] + (\pi^{yyx})^a_{bcd} [y_1^b,[y_1^c,x_1^d]] + (\pi^{yyy})^a_{bcd} [y_1^b,[y_1^c,y_1^d]]\\
	& + \lambda^a_b x_2^b + \mu^a_b y_2^b +(\nu^{xx})^a_{bc} [x_2^b,x_2^c]  + (\nu^{xy})^a_{bc} [x_2^b,y_2^c] + (\nu^{yy})^a_{bc} [y_2^b,y_2^c]\\
	&\qquad + (\pi^{xxx})^a_{bcd} [x_2^b,[x_2^c,x_2^d]] + (\pi^{xxy})^a_{bcd} [x_2^b,[x_2^c,y_2^d]] + (\pi^{yyx})^a_{bcd} [y_2^b,[y_2^c,x_2^d]] + (\pi^{yyy})^a_{bcd} [y_2^b,[y_2^c,y_2^d]]\\
	& -  \Big(\sum_b\lambda^a_b(\nu^{xy})^a_{db} +\mu^a_b((\nu^{xx})^a_{bd} - (\nu^{xx})^a_{db}) - \frac12\lambda^a_d \Big)[x_1^d, t_{12}] \\
	&\qquad +\Big(  \sum_b \lambda^a_b((\nu^{yy})^a_{bd} - (\nu^{yy})^a_{db}) - \mu^a_b(\nu^{xy})^a_{bd} + \frac12\mu^a_d \Big)[y_1^d,t_{12}]\\
	&\qquad + (\mbox{tot.deg.}\geq 4).
\end{align*}
For the left-hand side, we have
\begin{align*}
	&\xi^a_{12} + (2s^a-\frac12)\cdot t_{12}\\
	&= (2s^a-\frac12)\cdot t_{12} + \lambda^a_b (x_1^b + x_2^b) + \mu^a_b (y_1^b + y_2^b)\\
	&\qquad  +(\nu^{xx})^a_{bc} ([x_1^b ,x_1^c] + [x_2^b, x_2^c]) + (\nu^{xy})^a_{bc} ([x_1^b,y_1^c] + 2\delta^{bc}t_{12} + [x_2^b,y_2^c]) + (\nu^{yy})^a_{bc} ([y_1^b,y_1^c] + [y_2^b,y_2^c])\\
	&\qquad + (\pi^{xxx})^a_{bcd} ([x_1^b,[x_1^c ,x_1^d]] + [x_2^b,[x_2^c, x_2^d]])
	+ (\pi^{xxy})^a_{bcd} ([x_1^b,[x_1^c,y_1^d]] + [x_2^b,[x_2^c,y_2^d]])\\
	&\qquad + (\pi^{yyx})^a_{bcd} ([y_1^b,[y_1^c ,x_1^d] ] + [y_2^b, [y_2^c,x_2^d]]) + (\pi^{yyy})^a_{bcd} ([y_1^b,[y_1^c ,y_1^d]] + [y_2^b,[y_2^c, y_2^d]]) \\
	&\qquad + (\mbox{tot.deg.}\geq 4).
\end{align*}
Since these two are equal, we obtain
\begin{align*}
	&\frac12 t_{12}  -  \Big(\sum_b\lambda^a_b(\nu^{xy})^a_{db} +\mu^a_b((\nu^{xx})^a_{bd} - (\nu^{xx})^a_{db}) - \frac12\lambda^a_d \Big)[x_1^d, t_{12}] \\
	&\qquad + \Big(  \sum_b\lambda^a_b((\nu^{yy})^a_{bd} - (\nu^{yy})^a_{db}) - \mu^a_b(\nu^{xy})^a_{bd} + \frac12\mu^a_d \Big)[y_1^d,t_{12}]\\
	&= (2s^a-\frac12)\cdot t_{12} + 2(\nu^{xy})^a_{bc}\delta^{bc} t_{12}
\end{align*}
In particular, we have, for each $1\leq d\leq g$,
\begin{itemize}
	\item bi-degree $(1,1)$: $ 2s^a + \sum_b 2(\nu^{xy})^a_{bb} - 1 = 0$,
	\item bi-degree $(1,2)$: $  \sum_b\lambda^a_b(\nu^{xy})^a_{db} +\mu^a_b((\nu^{xx})^a_{bd} - (\nu^{xx})^a_{db}) - \frac12\lambda^a_d = 0$,
	\item bi-degree $(2,1)$: $\sum_b \lambda^a_b((\nu^{yy})^a_{bd} - (\nu^{yy})^a_{db}) - \mu^a_b(\nu^{xy})^a_{bd} + \frac12\mu^a_d = 0$,
\end{itemize}
Therefore, by the description of $\mathsf{fr}_{\vec Z, f}$ in the proof of Theorem \ref{thm:assoc}, we have the following:
\begin{proposition}
Using the notations above, $\mathsf{fr}_{\vec Z, f}(A_a^{1,2}) = -2s^a = \sum_b 2(\nu^{xy})^a_{bb} - 1$ holds. Therefore, the associated framing on $\Sigma_{g,1}$ is determined only by the coefficients of the quadratic terms $[x_1^b, y_1^b]$.
\end{proposition}

In the special case of $g = 1$, these equations read
\begin{itemize}
	\item bi-degree $(1,1)$: $ 2s^1 + 2(\nu^{xy})^1_{11} - 1 = 0$,
	\item bi-degree $(1,2)$: $ \lambda^1_1 \left( (\nu^{xy})^1_{11} - \frac12\right) = 0$, and 
	\item bi-degree $(2,1)$: $ -\mu^1_1 \left( (\nu^{xy})^1_{11} - \frac12\right) = 0$.
\end{itemize}
Since $\vec Z$ is an isomorphism of complete Hopf algebras and $x_1^1$ and $y_1^1$ have degree $1$, we have $\lambda^1_1\neq 0$ or $\mu^1_1\neq 0$ and therefore $(\nu^{xy})^1_{11} = \frac12$. This implies $s^1 = 0$, and it is similarly done for $B_a^{1,2}$. Hence, we obtained the following:

\begin{proposition}
The only framing on $\Sigma_{1,1}$ coming from a genus $1$ Gonzalez--Drinfeld associator is given by a constant vector field on a flat torus.
\end{proposition}

\begin{remark}
\noindent This framing is called the \textit{adapted framing} in \cite{akkn}. We also note that the above calculation does not use the assumption $\mathrm{gr}\,(Z_p) = \id$ (which is equivalent to $\lambda^a_b = \delta_{ab}$ and $\mu^a_b = 0$ in the above equations).
\end{remark}

Now we discuss the case of arbitrary $n$. Let $\vec\alpha_a\in\vec\pi$ be an $A$-loop on $\Sigma_{g,n+1}$ which is $\mathbf{PaB}^f(n+2)$-conjugate to $A_a^{*, p'}$ in $\mathbf{PaB}^f_g(n+2)$ for some parenthesisation $p'$ of $1\,2\,\cdots \,n \, 0$. This corresponds to the generator in Definition \ref{def:tder}. We have already discussed in the proof of \ref{thm:assoc} that the framing only depends on the $\mathbf{PaB}^f(n+2)$-conjugacy class, so we have $\mathsf{fr}_{\vec Z, f}(\vec\alpha_a) = \mathsf{fr}_{\vec Z, f}(A_a^{*, p'})$. Putting
\[
	\vec Z(A^{*,0}_a) = \exp(\xi^a_* + s^a\cdot t_{**}) \in \exp( \hat L(x_*^a, y_*^a)_{1\leq a\leq g} \oplus \mathbb{K}t_{**})
\]
as before, we have $\vec Z(A^{*,p'}_a) = \vec Z(A^{*,0}_a) \circ_0 \id_{p'} = \exp(\xi^a_* + s^a\cdot t_{**})$ as well, since the expression for $\vec Z(A^{*,0}_a)$, once written only in terms of $x_*^a, y_*^a$ and $t_{**}$, does not involve the index $0$. This shows $\mathsf{fr}_{\vec Z, f}(\vec\alpha_a) = -2s^a$. The same applies for $B$-loops, yielding the following result.

\begin{theorem}\label{thm:framing}
Let $D$ be the embedded disk in Figure \ref{fig:points}. Then, a framing on $\Sigma_{1,n+1}$ coming from a genus $1$ Gonzalez--Drinfeld associator always restricts to a constant vector field outside $D$.\\
\end{theorem}

\section{Proof of Lemma \ref{lem:exactness}}\label{sec:proof}
We begin with the first row of the diagram in Lemma \ref{lem:exactness}.

\begin{lemma}\label{lem:exactgrp}
The sequence
\[\begin{tikzcd}
	1 \arrow[r]& \vec\pi \arrow[r, "\iota"]& \mathrm{PB}^f_{g,12\cdots n*0} \arrow[r, "\varepsilon^*", yshift = 2pt]& \mathrm{PB}^f_{g,12\cdots n0} \arrow[r]\arrow[l, yshift = -2pt, "\circ_{0}\id_{*0}"]& 1
\end{tikzcd}\]
is split, and therefore defines the semi-direct product  $\mathrm{PB}^f_{g,12\cdots n*0} = \mathrm{PB}^f_{g,12\cdots n0}\ltimes \vec\pi$.
\end{lemma}
\noindent Proof. In Section \ref{sec:drinfeld}, we defined $F_m(\Sigma_g)$ as the total space of the pull-back bundle of $(U\Sigma_g)^{\times m} \to \Sigma_g^{\times m}$ along $\mathrm{Conf}_m(\Sigma_g) \hookrightarrow \Sigma_g^{\times m}$. Since we have the locally trivial fibration $F_{m+1}(\Sigma_g) \to F_{m}(\Sigma_g)$ with the fibre $\Sigma_{g,m}\times S^1$ by forgetting the point labelled $m+1$, we have the homotopy exact sequence
\[
	\pi_2(F_{12\cdots n*0}(\Sigma_g)) \xrightarrow{\varepsilon^*} \pi_2(F_{12\cdots n0}(\Sigma_g)) \to \vec\pi \to \mathrm{PB}^f_{g,12\cdots n*0} \to \mathrm{PB}^f_{g,12\cdots n0} \to 1.
\]
Furthermore, the map $\varepsilon^*\colon F_{12\cdots n*0}(\Sigma_g) \to F_{12\cdots n0}(\Sigma_g)$ admits a continuous global section by doubling the point labelled by $*$ in the direction of specified framing (such that it induces $\circ_{0}\id_{*0}$ on the fundamental groups). Therefore, $\varepsilon^*$ is surjective on $\pi_2$ and we obtain the claimed split sequence.\qed

\begin{remark}
For $g\geq 1$, we have an alternative proof: since $\Sigma_{g,m}\times S^1$ is an Eilenberg--Mac Lane space, we have the exact sequence
\[
	1 \to \pi_2(F_{m+1}(\Sigma_g)) \xrightarrow{\varepsilon^*} \pi_2(F_{m}(\Sigma_g)),
\]
and therefore the inclusion $\pi_2(F_{m}(\Sigma_g)) \subset \pi_2(F_{1}(\Sigma_g))$ by induction. The space $F_{1}(\Sigma_g)$ is just $\Sigma_g$, which is also an Eilenberg--Mac Lane space, so we have $\pi_2(F_{m}(\Sigma_g)) = 1$ for all $m\geq 1$.
\end{remark}

The sequence in Lemma \ref{lem:exactness} is obtained by taking the Malcev completion of the above, but it might end up with a non-exact sequence, so we need the following lemma.

\begin{lemma}\label{lem:malcev}
Let $G = H\ltimes K$ be a semi-direct product of groups, and suppose that there exists a subgroup $L\subset K$ satisfying the following conditions:
\begin{enumerate}[(1)]
	\item $[H,K] \subset L$;\label{item:malcev1}
	\item $[H,L] \subset [K, K]$.\label{item:malcev2}
\end{enumerate}
Then, $\hat H$ acts on $\hat K$ and we have $\hat G\cong \hat H\ltimes \hat K$ as topological groups.
\end{lemma}

\begin{remark}
If we can take $L = [K, K]$, the semi-direct product $G = H\ltimes K$ is said to be an \textit{almost direct product}. In that case, the lemma above is reduced to Proposition 8.5.3 in \cite{fresse}. Our proof is based on their method.
\end{remark}

\noindent Proof. We have the natural isomorphism
\begin{align}
	\mathbb{K}G \cong \mathbb{K}H\otimes \mathbb{K}K \label{eq:isomcoalg}
\end{align}
of coalgebras. Recall that the filtration $\mathbb{I}^\bullet G$ on the left-hand side is defined using the augmentation ideal $\mathbb{I}^1G = \mathbb{I}G$ of $G$, and the Malcev completion $G$ over $\mathbb{K}$ is the set of group-like elements of $\widehat{\mathbb{K}G}$. On the other hand, the filtration on the right-hand side is given by
\[
	F^m := \sum_{\substack{p,q\geq 0\\p+q=m}} \mathbb{I}^pH\otimes \mathbb{I}^qK.
\]	
for $m\geq 0$. We check these filtrations lead to the same completion $\widehat{\mathbb{K}G} \cong \widehat{\mathbb{K}H}\hat\otimes \widehat{\mathbb{K}K}$; more precisely, we show $\mathbb{I}^mG \supset F^m\supset \mathbb{I}^{3^m}G$ for $m\geq 0$ under the identification \eqref{eq:isomcoalg}.

Since the $\mathbb{K}$-vector space $\mathbb{I}G$ is spanned by the elemets of the form $g-1$ for $g\in G$, and the same applies for $H$ and $K$, $F^m$ is $\mathbb{K}$-linearly spanned by the elements
\[
	(h_1 - 1)\cdots (h_p-1) (k_1-1)\cdots (k_q-1)
\]
where $h_i\in H$, $k_j\in K$ and $p+q=m$. This is obviously contained in $\mathbb{I}^mG$, so we have $F^m \subset \mathbb{I}^mG$.

Next, we show $\mathbb{I}K \mathbb{I}^{2p}H \subset \mathbb{I}^pG \mathbb{I}K$ for $p\geq 1$ by induction on $p$. Firstly, we have
\begin{align*}
	(k - 1)(h - 1) &=  (h - 1)(k^h - 1) + (k^h - k)\\
	&= (h - 1)(k^h - 1) + (k - 1)(l - 1)+ (l - 1) 
\end{align*}
where $k\in K$, $h\in H$ and we put $k^h = h^{-1}kh$ and $k^h = kl$ with $l\in L$, which is possible by the condition \ref{item:malcev1}. This shows $\mathbb{I}K\mathbb{I}H \subset \mathbb{I}H\mathbb{I}K + \mathbb{I}K\mathbb{I}L + \mathbb{I}L$, which is further contained in $\mathbb{I}H\mathbb{I}K + \mathbb{I}K$. Next, we have
\begin{align*}
	(l - 1)(h-1) &= (h-1)(l^h - 1) + (l^h - l).
\end{align*}
In addition, we have $l^h = l (u_1, v_1)\cdots (u_r, v_r)$ with $u_i, v_i \in K$ by the condition \ref{item:malcev2} where we put $(a,b) = aba^{-1}b^{-1}$. Therefore,
\[
	l^h - l = l\big( (u_1, v_1)\cdots (u_r, v_r) - 1 \big) \in \mathbb{K}L\cdot \mathbb{I}^2K \subset \mathbb{I}^2K
\]
which can be shown as in Proposition 8.5.3 of \cite{fresse}, hence $\mathbb{I}L\mathbb{I}H \subset \mathbb{I}H\mathbb{I}K + \mathbb{I}^2K$. Combining these, we have
\begin{align*}
	\mathbb{I}K\mathbb{I}^2H &\subset \big( \mathbb{I}H\mathbb{I}K + \mathbb{I}K\mathbb{I}L + \mathbb{I}L \big) \mathbb{I}H\\
	&\subset \mathbb{I}H\big( \mathbb{I}K\mathbb{I}H \big) + \mathbb{I}K \big( \mathbb{I}L\mathbb{I}H \big) + \mathbb{I}L\mathbb{I}H\\
	&\subset \mathbb{I}H\big( \mathbb{I}H\mathbb{I}K + \mathbb{I}K \big) + \mathbb{I}K \big( \mathbb{I}H\mathbb{I}K + \mathbb{I}^2K \big) + \mathbb{I}H\mathbb{I}K + \mathbb{I}^2K\\
	&\subset \mathbb{I}G\mathbb{I}K.
\end{align*}
This shows the case of $p=1$. For $p\geq 2$, we have
\begin{align*}
	\mathbb{I}K\mathbb{I}^{2p}H &= \mathbb{I}K\mathbb{I}^2H\mathbb{I}^{2(p-1)}H\\
	&\subset \mathbb{I}G\mathbb{I}K\cdot \mathbb{I}^{2(p-1)}H\\
	&\subset \mathbb{I}^{p-1+1}G\mathbb{I}K.
\end{align*}
Here, we used the induction hypothesis for the last inclusion. This shows $\mathbb{I}K \mathbb{I}^{2p}H \subset \mathbb{I}^pG \mathbb{I}K$ for $p\geq 1$.

With this seen, we show $\mathbb{I}^{3^m}G\subset F^m$ for $m\geq 0$ by induction on $m$. For $m=0$, this is evident since $\mathbb{I}^1G\subset F^0 = \mathbb{K}G$ clearly holds. Next, for $g \in G$, we can write $g = hk$ for some $h\in H$ and $k\in K$ and therefore
\[
	g - 1 = hk - 1 = (h-1)(k-1) + (h-1) + (k-1).
\]
In particular, for $p\geq 1$, an element of $\mathbb{I}^pG$ can be written as a linear combination of the products of $p$ factors of the forms $(h-1)(k-1)$, $(h-1)$ or $(k-1)$. Therefore, we have
\begin{align*}
	\mathbb{I}^pG \subset \mathbb{I}^pH + \sum_{p'\geq p - 1}\mathbb{I}^{p'}G \cdot \mathbb{I}K + \sum_{\substack{p'\geq 0, 1\leq q< p\\p'\geq p-q-1}} \mathbb{I}^{p'}G \cdot \mathbb{I}K\mathbb{I}^qH
\end{align*}
where the first term corresponds to the terms with no factors in $\mathbb{I}K$ appear, the second term to the ones whose rightmost factor is in $\mathbb{I}K$, and the third term is obtained by noticing that there is a factor in $\mathbb{I}K$ somewhere, say, the $(q+1)$-th factor from the right side. For $p\geq 2$, this is contained in
\begin{align*}
	&F^p + \sum_{p'\geq p - 1}\mathbb{I}^{p'}G \cdot \mathbb{I}K + \sum_{p'\geq p-2} \mathbb{I}^{p'}G \cdot \mathbb{I}K\mathbb{I}H + \sum_{\substack{p'\geq 0, 2\leq q< p\\p'\geq p-q-1}} \mathbb{I}^{p'+\lfloor q/2\rfloor}G\mathbb{I}K\\
	&\subset F^p + \mathbb{I}^{p-2}G \cdot \mathbb{I}K + \sum_{\substack{p'\geq 0, 2\leq q< p\\p'\geq p-q-1}} \mathbb{I}^{p'+\lfloor q/2\rfloor}G\mathbb{I}K
\end{align*}
Now let $m\geq 1$ and $p = 3^m$. Since we have $p-2 = 3^m - 2 \geq 3^{m-1}$, the first term is contained in $F^m$ since $m \leq p=3^m$, and the second term is contained in $\mathbb{I}^{3^{m-1}}G \mathbb{I}K \subset F^{m-1}\mathbb{I}K \subset F^m$ by the induction hypothesis. Next, we have
\begin{align*}
	p'+\lfloor q/2\rfloor &\geq p-q-1+\lfloor q/2\rfloor \geq p-(p-1)-1+ (p-1)/2 \\
	&= \frac{3^m - 1}{2}  \geq 3^{m-1}.
\end{align*}
Therefore, we have
\begin{align*}
	\sum_{\substack{p'\geq 0, 2\leq q< p\\p'\geq p-q-1}} \mathbb{I}^{p'+\lfloor q/2\rfloor}G\mathbb{I}K \subset \mathbb{I}^{3^{m-1}}G \mathbb{I}K \subset F^{m-1}\mathbb{I}K \subset F^m.
\end{align*}
This shows $\mathbb{I}^{3^m}G\subset F^m$ for $m\geq 0$ and therefore we have $\widehat{\mathbb{K}G} \cong \widehat{\mathbb{K}H}\hat\otimes \widehat{\mathbb{K}K}$ as topological coalgebras.

Finally, if the action map $\rho\colon \mathbb{K}K\otimes \mathbb{K}H \to \mathbb{K}K\colon k\otimes h \mapsto k^h$ is continuous, we can extend $\rho$ to the completions so that we have $\widehat{\mathbb{K}G} \cong \widehat{\mathbb{K}H}\sharp \widehat{\mathbb{K}K}$ as topological algebras and therefore as topological (complete) Hopf algebras. Here, the symbol $\sharp$ denotes the semi-direct product of Hopf algebras as in \cite{fresse}. Then, taking the group-like element part yields $\hat G\cong \hat H\ltimes \hat K$ as noted in Proposition 8.5.3 of \cite{fresse}. To prove that $\rho$ is continuous with respect to the filtration by the augmentation ideal, notice that the action map is equal to the composition
\[
	\mathbb{K}K\otimes \mathbb{K}H \xrightarrow{\id\otimes\Delta} \mathbb{K}K\otimes \mathbb{K}H \otimes \mathbb{K}H \xrightarrow{\id\otimes S\otimes \id} \mathbb{K}K\otimes \mathbb{K}H \otimes \mathbb{K}H \xrightarrow{\tau\otimes \id}  \mathbb{K}H\otimes \mathbb{K}K \otimes \mathbb{K}H \xrightarrow{\mu} \mathbb{K}G\,,
\]
where $\Delta$ is the coproduct of $\mathbb{K}H$, $S$ is the antipode of $\mathbb{K}H$, $\tau$ is the transposition map and $\mu$ is the multiplication map in $\mathbb{K}G$, all of which are continuous. This concludes the proof.\qed\\

To apply Lemma \ref{lem:malcev} to our case of braid groups, we check the conditions \ref{item:malcev1} and \ref{item:malcev2} for $H = \mathrm{PB}^f_{g,12\cdots n0}$ and $K = \vec\pi$. We use the result of Bellingeri and Gervais (see Theorem 8 in \cite{bg}): the framed pure braid group $\mathrm{PB}^f_{g,m}$ on $\Sigma_g$ for $g\geq 0, m\geq 1$ is generated by
\[
	B_{i,j} \;(1\leq i\leq 2g+m-1, 2g+1\leq j\leq 2g+m, i<j)\;\mbox{ and } \;f_k\; (1\leq k\leq m),
\]
and we can deduce the following (non-exhaustive) set of relations: putting $s = 2g+m-1$, we have
\begin{align}
	\begin{aligned}\label{eq:braidrel}
	(B_{r,s}^{-1}, B_{r+1, j}^{-1}) &\equiv B_{j,s}^{-1} \mbox{ for } r\leq2g \mbox{ and } r \mbox{ odd}, j<s,\\
	(B_{r,s}^{-1}, B_{r-1, j}^{-1}) &\equiv B_{j,s} \mbox{ for } r\leq 2g \mbox{ and } r \mbox{ even}, j<s,\\
	(B_{r,2g+m}^{-1}, B_{r+1, s}^{-1}) &\equiv B_{s,2g+m}^{-1} \mbox{ for } r\leq2g \mbox{ and } r \mbox{ odd},\\
	(B_{r,2g+m}^{-1}, B_{r-1, s}^{-1}) &\equiv B_{s,2g+m} \mbox{ for } r\leq 2g \mbox{ and } r \mbox{ even},\\
	(B_{r,s}^{-1}, B_{i, j}^{-1}) &\equiv 1 \mbox{ otherwise},\\
	f_k &\mbox{ is central, and}\\
	(B_{2g, s}^{-1}, B_{2g-1, s})\cdots (B_{2, s}^{-1}, B_{1, s}) &= B_{2g+1,s}\cdots B_{2g+m-2, s}B_{s, 2g+m} f_{m-1}^{2(g-1)}.
	\end{aligned}
\end{align}
Here we set $m = n+2$ so that $G = \mathrm{PB}^f_{g,12\cdots n*0} = \mathrm{PB}^f_{g,m}$, $K$ is identified with the subgroup of $G$ generated by $B_{i, s} \,(1\leq i< s)$ and $f_{m-1}$, and $\equiv$ indicates an equality modulo $[K,K]$. We remark that, even though their result is for $g\geq 1$, these relations hold for $g=0$; in fact, the $g=0$ case is a quotient of the usual pure braid group.

\begin{lemma}
The subgroup $L := \langle [K,K], B_{r, s} \,(2g+1\leq r < s), f_{m-1}\rangle$ of $K$ is normal.
\end{lemma}
\noindent Proof. $[K,K]$ being closed under conjugation by $K$ is standard. In addition, for $2g + 1\leq  r < s$ and $1\leq i\leq s$, the element $(B_{i,s}, B_{r,s})$ is obviously contained in $[K,K] \subset L$. Finally, since $f_{m-1}$ is central, we deduce that $L$ is normal. \qed

\begin{remark}
The group $L$ corresponds to the degree $\geq 2$ part of the Lie algebra $\mathfrak{u}^f_{g,12\cdots n}$ below.
\end{remark}

\begin{lemma}\label{lem:normal}
Let $G = H\ltimes K$ be a semi-direct product, $N$ a normal subgroup of $K$, $h, h_1, h_2\in H$ and $k, k_1, k_2 \in K$.
\begin{enumerate}[(1)]
	\item If $(h, k_1), (h, k_2) \in N$, then $(h, k_1k_2) \in N$.
	\item If $(h_1, k), (h_2, k) \in N$, then $(h_1h_2, k) \in N$.
\end{enumerate}
\end{lemma}
\noindent Proof. (1) We compute
\begin{align*}
	(h, k_1k_2) = hk_1k_2h^{-1}k_2^{-1}k_1^{-1}
	= hk_1h^{-1}k_1^{-1}\cdot k_1h k_2h^{-1}k_2^{-1}k_1^{-1}
	= (h, k_1) (h, k_2)^{k_1^{-1}} \in N
\end{align*}
since $N$ is a normal subgroup of $K$.

\noindent (2) We compute
\begin{align*}
	(h_1h_2, k) = h_1h_2kh_2^{-1}h_1^{-1}k^{-1} = h_1h_2kh_2^{-1}k^{-1}h_1^{-1} \cdot h_1kh_1^{-1}k^{-1} = (h_2, k)^{h_1^{-1}} (h_1, k) \in N
\end{align*}
since $N$ is a normal subgroup of $K$.\qed\\

\noindent Applying Lemma \ref{lem:normal} to the case $N = L$ and $N = [K,K]$, we only have to check the conditions on generators. 
\begin{itemize}
	\item $[H, K] \subset L \colon$ By the first five relations in \eqref{eq:braidrel}, $[H,K]$ is contained in the subgroup generated by $[K,K]$, $B_{j,s}$ with $2g+1\leq j < s$, and $B_{s, 2g+m}$. By the last relation in \eqref{eq:braidrel}, we can rewrite $B_{s, 2g+m}$ into the product of $B_{j,s}$'s and $f_{m-1}$. This shows $[H, K] \subset L$.
	\item $[H, L] \subset [K, K] \colon$ The inclusion $[H, [K,K]] \subset [K,K]$ is standard. For $r\geq 2g+1$, we have $[H, B_{r,s}]\equiv 1$ and therefore $[H, B_{r,s}]\subset [K,K]$. Finally, since $f_{m-1}$ is central, we deduce $[H, L] \subset [K, K]$.
\end{itemize}
This completes the proof of the first row of Lemma \ref{lem:exactness}.\\

Next, we move on to the second row of the diagram in Lemma \ref{lem:exactness}. We will show in Theorem \ref{thm:split}, for $g,n \geq 0$, that
\[\begin{tikzcd}
	0 \arrow[r]& \mathfrak{u}^f_{g,12\cdots n} \arrow[r, "\iota"]& \mathfrak{t}^f_{g,12\cdots n*0} \arrow[r, "\varepsilon^*", yshift = 2pt]& \mathfrak{t}^f_{g,12\cdots n0} \arrow[r]\arrow[l, yshift = -2pt, "\circ_{0}\id_{*0}"]& 0
\end{tikzcd}\]
is split, where $ \mathfrak{u}^f_{g,12\cdots n} = L(H)\oplus\mathbb{K}t_{**} $, $H = \mathrm{Span}_\mathbb{K}\{t_{i*}, x^a_*, y^a_*\}_{1\leq i\leq n, 1\leq a\leq g}$, and $\iota$ is the natural map (which is not yet shown to be injective). We set
\begin{align*}
	\omega_* = \sum_{1\leq a\leq g} [x^a_*, y^a_*] + \sum_{1\leq j\leq n} t_{j*} - (g-1)t_{**} \quad \in  \mathfrak{u}^f_{g,12\cdots n}
\end{align*}
so that $\iota(\omega_*) = -t_{*0}$ holds in $ \mathfrak{t}^f_{g,12\cdots n*0}$.

\begin{deflem}
We define the action of $\mathfrak{t}^f_{g,12\cdots n0}$ on $\mathfrak{u}^f_{g,12\cdots n}$ by the following table, where $1\leq i,j,k \leq n$ and $1\leq a,b\leq g$, and the value of the $z$-row and the $w$-column computes $\ad_z(w)$.\\

\SetTblrInner{rowsep=0pt}
\centerline{
\begin{tabular}{|c||c|c|c|c|}\hline
	& $t_{k*}$ & $x_*^b$ & $y_*^b$ & $t_{**}$\\\hhline{|=#=|=|=|=|}
	$t_{ij}$ & $[t_{k*}, \delta_{ik}t_{j*} + \delta_{jk}t_{i*}]$ & $0$ & $0$ & $0$ \\\hline
	$t_{i0}$ & $[t_{i*} + \delta_{ik}\omega_*, t_{k*}]$ & $[t_{i*}, x^b_*]$ & $[t_{i*}, y^b_*]$ & $0$ \\\hline
	$t_{00}$ & $2[t_{k*},\omega_*]$ & $2[x_*^b,\omega_*]$ & $2[y_*^b,\omega_*]$ & $0$ \\\hline
	$x_i^a$ & $\delta_{ik}[t_{k*}, x_*^a]$ & $0$ & $\delta_{ab}t_{i*}$ & $0$ \\\hline
	$x_0^a$ & $[x_*^a, t_{k*}]$ & $[x_*^a, x_*^b]$ & $[x_*^a, y_*^b] - \delta_{ab}\omega_*$ & $0$ \\\hline
	$y_i^a$ & $\delta_{ik}[t_{k*}, y_*^a]$ & $-\delta_{ab}t_{i*}$ & $0$ & $0$ \\\hline
	$y_0^a$ & $[y_*^a, t_{k*}]$ & $[y_*^a, x_*^b] + \delta_{ab}\omega_*$ & $[y_*^a, y_*^b]$ & $0$ \\\hline
\end{tabular}}
\end{deflem}
\noindent Proof. We check that all the relations in $\mathfrak{t}^f_{g,12\cdots n0}$ are satisfied. Let $1\leq i,j,k,l,p\leq n$ and $1\leq a,b,c\leq g$. We first calculate
\begin{align*}
	t_{ij}\cdot \omega_* &= \sum_{1\leq b\leq g} t_{ij}\cdot [x^b_*, y^b_*] + \sum_{1\leq k\leq n} t_{ij}\cdot t_{k*} - (g-1) t_{ij}\cdot t_{**}\\
	&= 0 + \sum_{1\leq k\leq n} [t_{k*}, \delta_{ik}t_{j*} + \delta_{jk}t_{i*}] + 0\\
	&=  [t_{i*}, t_{j*}] + [t_{j*}, t_{i*}] = 0,\\
	t_{i0}\cdot \omega_* &= \sum_{1\leq b\leq g} t_{i0}\cdot [x^b_*, y^b_*] + \sum_{1\leq j\leq n} t_{i0}\cdot t_{j*} - t_{i0}\cdot (g-1)t_{**}\\
	&= \sum_{1\leq b\leq g}  [[t_{i*}, x^b_*], y^b_*] + [x^b_*, [t_{i*},y^b_*]] + \sum_{1\leq k\leq n} [t_{i*} + \delta_{ik}\omega_*, t_{k*}] - 0\\
	&= [t_{i*}, \omega_*] +  [\omega_*, t_{i*}] = 0,\\
	x_i^a \cdot \omega_* &= \sum_{1\leq b\leq g} x_i^a\cdot [x^b_*, y^b_*] + \sum_{1\leq j\leq n} x_i^a\cdot t_{j*} - x_i^a\cdot (g-1)t_{**}\\
	&= \sum_{1\leq b\leq g} [x^b_*, \delta_{ab} t_{i*}] + \sum_{1\leq j\leq n} \delta_{ij}[t_{j*}, x_*^a]\\
	&= [x^a_*, t_{i*}] + [t_{i*}, x_*^a] = 0,\\
	x_0^a \cdot \omega_* &= \sum_{1\leq b\leq g} x_0^a\cdot [x^b_*, y^b_*] + \sum_{1\leq j\leq n} x_0^a\cdot t_{j*} - x_0^a\cdot (g-1)t_{**}\\
	&= \sum_{1\leq b\leq g} [[x_*^a, x_*^b], y^b_*] +  [x^b_*, [x_*^a, y_*^b] - \delta_{ab}\omega_*] + \sum_{1\leq j\leq n} [x_*^a, t_{j*}]\\
	&=  - [x^a_*, \omega_*] + \sum_{1\leq b\leq g}  [x_*^a, [x_*^b, y^b_*]] + \sum_{1\leq j\leq n} [x_*^a, t_{j*}] = 0.
\end{align*}
\begin{itemize}
	\item \underline{$[t_{ij}, t_{kl}] = 0$ ($\{i,j\}\cap \{k,l\} = \varnothing$)}: We have
	\begin{align*}
		&t_{ij}\cdot (t_{kl}\cdot t_{p*})\\
		&= t_{ij}\cdot [t_{p*}, \delta_{kp}t_{l*} + \delta_{lp}t_{k*}] \\
		&= [t_{ij}\cdot t_{p*}, \delta_{kp}t_{l*} + \delta_{lp}t_{k*}] + \delta_{kp} [t_{p*}, t_{ij}\cdot t_{l*}] + \delta_{lp} [t_{p*}, t_{ij}\cdot t_{k*}]\\
		&= [[t_{p*}, \delta_{ip}t_{j*} + \delta_{jp}t_{i*}], \delta_{kp}t_{l*} ] + [[t_{p*}, \delta_{ip}t_{j*} + \delta_{jp}t_{i*}], \delta_{lp}t_{k*}]\\
		&\qquad + \delta_{kp} [t_{p*}, [t_{l*}, \delta_{il}t_{j*} + \delta_{jl}t_{i*}]] + \delta_{lp} [t_{p*}, [t_{k*}, \delta_{ik}t_{j*} + \delta_{jk}t_{i*}]]\\
		&= 0.
	\end{align*}
	Here we used $\delta_{ip} \delta_{kp} =  \delta_{ip} \delta_{lp} = \delta_{jp} \delta_{kp} = \delta_{jp}\delta_{lp} = \delta_{il} = \delta_{jl} = \delta_{ik} = \delta_{jk} = 0$ in the last equality, which follow from $\{i,j\}\cap \{k,l\} = \varnothing$. Likewise, we have $t_{kl}\cdot (t_{ij}\cdot t_{p*}) = 0$, which shows $t_{ij}\cdot (t_{kl}\cdot z) - t_{kl}\cdot (t_{ij}\cdot z) = 0$ for any $z\in \mathfrak{u}^f_{g,12\cdots n}$.
	\item \underline{$[t_{ij}, t_{k0}] = 0$ ($\{i,j\}\cap \{k\} = \varnothing$)}: We have
	\begin{align*}
		&t_{ij}\cdot (t_{k0}\cdot t_{p*})\\
		&= t_{ij}\cdot [t_{k*} + \delta_{kp}\omega_*, t_{p*}]\\
		&= [[t_{k*}, \delta_{ik}t_{j*} + \delta_{jk}t_{i*}], t_{p*}] + \delta_{kp} [t_{ij}\cdot \omega_*, t_{p*}] + [t_{k*}, [t_{p*}, \delta_{ip}t_{j*} + \delta_{jp}t_{i*}]] + \delta_{kp}[\omega_*, [t_{p*}, \delta_{ip}t_{j*} + \delta_{jp}t_{i*}]]\\
		&= \delta_{ip} [t_{k*}, [t_{p*},t_{j*}]] + \delta_{jp} [t_{k*}, [t_{p*}, t_{i*}]] + \delta_{kp}\delta_{ip} [\omega_*, [t_{p*}, t_{j*}]] + \delta_{kp}\delta_{jp} [\omega_*, [t_{p*}, t_{i*}]]\\
		&= \delta_{ip} [t_{k*}, [t_{p*},t_{j*}]] + \delta_{jp} [t_{k*}, [t_{p*}, t_{i*}]].
	\end{align*}
	On the other hand, we have
	\begin{align*}
		&t_{k0}\cdot (t_{ij}\cdot t_{p*})\\
		&= t_{k0}\cdot [t_{p*}, \delta_{ip}t_{j*} + \delta_{jp}t_{i*}]\\
		&= [[t_{k*} + \delta_{kp}\omega_*, t_{p*}], \delta_{ip}t_{j*} + \delta_{jp}t_{i*}] + \delta_{ip} [t_{p*},  [t_{k*} + \delta_{kj}\omega_*, t_{j*}]] + \delta_{jp} [t_{p*}, [t_{k*} + \delta_{ki}\omega_*, t_{i*}]]\\
		&= \delta_{ip}( -[t_{j*}, [t_{k*}, t_{p*}]] - [t_{p*},  [t_{j*}, t_{k*}]]) + \delta_{jp}( - [t_{i*}, [t_{k*}, t_{p*}]] - [t_{p*}, [t_{i*}, t_{k*}]])\\
		&= \delta_{ip} [t_{k*}, [t_{p*},t_{j*}]] + \delta_{jp} [t_{k*}, [t_{p*}, t_{i*}]],
	\end{align*}
	so these two are equal and therefore $t_{ij}\cdot (t_{k0}\cdot t_{p*}) - t_{k0}\cdot (t_{ij}\cdot t_{p*}) = 0$. Next, we have
	\begin{align*}
		t_{ij}\cdot (t_{k0}\cdot x_*^a) - t_{k0}\cdot (t_{ij}\cdot x_*^a) &= t_{ij}\cdot [t_{k*}, x_*^a] - 0\\
		&= [[t_{k*}, \delta_{ik}t_{j*} + \delta_{jk}t_{i*}], x_*^a] = 0.
	\end{align*}
	It is similar for $y_*^a$. This shows $t_{ij}\cdot (t_{k0}\cdot z) - t_{k0}\cdot (t_{ij}\cdot z) = 0$ for any $z\in \mathfrak{u}^f_{g,12\cdots n}$.
	\item \underline{$[t_{ij}, t_{00}] = 0$ ($\{i,j\}\cap \{k\} = \varnothing$)}: Since the action of $t_{00}$ is the inner derivation by $-2\omega_*$ and also $t_{ij}\cdot \omega_* = 0$, we have
	\begin{align*}
		t_{ij}\cdot (t_{00}\cdot z) - t_{00} \cdot (t_{ij}\cdot z) &= t_{ij}\cdot 2[z, \omega_*] + 2 [\omega_*, t_{ij}\cdot z]\\
		&= 2[t_{ij}\cdot z, \omega_*] + 2 [\omega_*, t_{ij}\cdot z] = 0
	\end{align*}
	for any $z\in \mathfrak{u}^f_{g,12\cdots n}$.
	\item \underline{$[t_{ij}, t_{ik} + t_{jk}] = 0$ ($\{i,j\}\cap \{k\} = \varnothing$)}: We have
	\begin{align*}
		&t_{ij}\cdot ((t_{ik} + t_{jk}) \cdot t_{p*})\\
		&= t_{ij}\cdot ([t_{p*}, \delta_{ip}t_{k*} + \delta_{kp}t_{i*}] + [t_{p*}, \delta_{jp}t_{k*} + \delta_{kp}t_{j*}])\\
		&\qquad + [t_{ij}\cdot t_{p*}, \delta_{jp}t_{k*} + \delta_{kp}t_{j*}] +  \delta_{jp}[t_{p*}, t_{ij}\cdot t_{k*}]  + \delta_{kp} [t_{p*}, t_{ij}\cdot t_{j*}]\\
		&= [[t_{p*}, \delta_{ip}t_{j*} + \delta_{jp}t_{i*}], \delta_{ip}t_{k*} + \delta_{kp}t_{i*}] + \delta_{ip}[t_{p*}, [t_{k*}, \delta_{ik}t_{j*} + \delta_{jk}t_{i*}]] + \delta_{kp} [t_{p*}, [t_{i*}, \delta_{ii}t_{j*} + \delta_{ji}t_{i*}]]\\
		&\qquad + [[t_{p*}, \delta_{ip}t_{j*} + \delta_{jp}t_{i*}], \delta_{jp}t_{k*} + \delta_{kp}t_{j*}] +  \delta_{jp}[t_{p*}, [t_{k*}, \delta_{ik}t_{j*} + \delta_{jk}t_{i*}]] + \delta_{kp} [t_{p*}, [t_{j*}, \delta_{ij}t_{j*} + \delta_{jj}t_{i*}]]\\
		&= (\delta_{ip} + \delta_{jp}) [[t_{p*}, \delta_{ip}t_{j*} + \delta_{jp}t_{i*}], t_{k*}] + \delta_{kp} ([t_{p*}, [t_{i*}, t_{j*} + \delta_{ji}t_{i*}]] +  [t_{p*}, [t_{j*}, \delta_{ij}t_{j*} + t_{i*}]])\\
		&= (\delta_{ip} + \delta_{jp}) [[t_{p*}, \delta_{ip}t_{j*} + \delta_{jp}t_{i*}], t_{k*}].
	\end{align*}
	On the other hand, we have
	\begin{align*}
		&(t_{ik} + t_{jk}) \cdot (t_{ij}\cdot t_{p*})\\
		&= (t_{ik} + t_{jk}) \cdot [t_{p*}, \delta_{ip}t_{j*} + \delta_{jp}t_{i*}]\\
		&\qquad + [t_{jk}\cdot t_{p*}, \delta_{ip}t_{j*} + \delta_{jp}t_{i*}] +\delta_{ip} [t_{p*},  t_{jk}\cdot t_{j*}] + \delta_{jp}[t_{p*},  t_{jk}\cdot t_{i*}]\\
		&= [[t_{p*}, \delta_{ip}t_{k*} + \delta_{kp}t_{i*}], \delta_{ip}t_{j*} + \delta_{jp}t_{i*}] +\delta_{ip} [t_{p*}, [t_{j*}, \delta_{ij}t_{k*} + \delta_{jk}t_{i*}]] + \delta_{jp}[t_{p*},  [t_{i*}, \delta_{ii}t_{k*} + \delta_{ik}t_{i*}]]\\
		&\qquad + [[t_{p*}, \delta_{jp}t_{k*} + \delta_{kp}t_{j*}], \delta_{ip}t_{j*} + \delta_{jp}t_{i*}] +\delta_{ip} [t_{p*},  [t_{j*}, \delta_{jj}t_{k*} + \delta_{jk}t_{j*}]] + \delta_{jp}[t_{p*},  [t_{i*}, \delta_{ji}t_{k*} + \delta_{ik}t_{j*}]]\\
		&= (\delta_{ip} + \delta_{jp}) [[t_{p*}, t_{k*}], \delta_{ip}t_{j*} + \delta_{jp}t_{i*}] + \delta_{ip} (1 + \delta_{ij}) [t_{p*}, [t_{j*},t_{k*}]] + \delta_{jp} (1 + \delta_{ij})[t_{p*},  [t_{i*}, t_{k*}]].
	\end{align*}
	Since we have
	\begin{align*}
		\delta_{ip} (1 + \delta_{ij}) &= \delta_{ip} + \delta_{ip} \delta_{ij} = \delta_{ip}^2 + \delta_{ip} \delta_{jp} = \delta_{ip}(\delta_{ip} + \delta_{jp})\;\mbox{ and}\\
		\delta_{jp} (1 + \delta_{ij}) &= \delta_{jp} + \delta_{jp} \delta_{ij} = \delta_{jp}^2 + \delta_{jp} \delta_{ip} = \delta_{jp}(\delta_{ip} + \delta_{jp}),
	\end{align*}
	we obtain
	\begin{align*}
		&(t_{ik} + t_{jk}) \cdot (t_{ij}\cdot t_{p*})\\
		&= (\delta_{ip} + \delta_{jp}) [[t_{p*}, t_{k*}], \delta_{ip}t_{j*} + \delta_{jp}t_{i*}] + \delta_{ip} (\delta_{ip} + \delta_{jp}) [t_{p*}, [t_{j*},t_{k*}]] + \delta_{jp} (\delta_{ip} + \delta_{jp})[t_{p*},  [t_{i*}, t_{k*}]]\\
		&= (\delta_{ip} + \delta_{jp}) [[t_{p*}, t_{k*}], \delta_{ip}t_{j*} + \delta_{jp}t_{i*}] + (\delta_{ip} + \delta_{jp}) [t_{p*}, [\delta_{ip}t_{j*} + \delta_{jp} t_{i*},t_{k*}]]\\
		&= (\delta_{ip} + \delta_{jp}) [[t_{p*}, \delta_{ip}t_{j*} + \delta_{jp}t_{i*}], t_{k*}]
	\end{align*}
	by the Jacobi identity. This shows $t_{ij}\cdot ((t_{ik} + t_{jk}) \cdot z) - (t_{ik} + t_{jk}) \cdot (t_{ij}\cdot z) = 0$ for any $z\in \mathfrak{u}^f_{g,12\cdots n}$.
	\item \underline{$[t_{i0}, t_{ik} + t_{k0}] = 0$ ($i\neq k$)}: We have
	\begin{align*}
		&t_{i0}\cdot ((t_{ik} + t_{k0})\cdot t_{p*})\\
		&= t_{i0}\cdot ([t_{p*}, \delta_{ip}t_{k*} + \delta_{kp}t_{i*}] + [t_{k*} + \delta_{kp}\omega_*, t_{p*}])\\
		&\qquad  + [t_{i0}\cdot t_{k*}, t_{p*}] + \delta_{kp} [t_{i0}\cdot \omega_*, t_{p*}] + [t_{k*} + \delta_{kp} \omega_*, t_{i0}\cdot t_{p*}]\\
		&= [[t_{i*} + \delta_{ip}\omega_*, t_{p*}], \delta_{ip}t_{k*} + \delta_{kp}t_{i*}] + \delta_{ip} [t_{p*}, [t_{i*} + \delta_{ik}\omega_*, t_{k*}]] + \delta_{kp} [t_{p*}, [t_{i*} + \delta_{ii}\omega_*, t_{i*}]]\\
		&\qquad  + [[t_{i*} + \delta_{ik}\omega_*, t_{k*}], t_{p*}] + \delta_{kp} [t_{i0}\cdot \omega_*, t_{p*}] + [t_{k*} + \delta_{kp} \omega_*, [t_{i*} + \delta_{ip}\omega_*, t_{p*}]]\\
		&= (\delta_{ip} - 1) [[t_{i*}, t_{p*}], t_{k*}] + (\delta_{ip} - 1) [t_{p*}, [t_{i*}, t_{k*}]] + \delta_{kp} [t_{p*}, [\omega_*, t_{i*}]] + \delta_{kp} [\omega_*, [t_{i*}, t_{p*}]] + \delta_{kp} [[t_{i*}, t_{p*}],t_{i*}]\\
		&=  -(\delta_{ip} - 1) [[t_{p*}, t_{k*}], t_{i*}] - \delta_{kp} [t_{i*}, [t_{p*}, \omega_* - t_{i*}]].
	\end{align*}
	On the other hand, we have
	\begin{align*}
		&(t_{ik} + t_{k0})\cdot (t_{i0}\cdot t_{p*})\\
		&= (t_{ik} + t_{k0})\cdot [t_{i*} + \delta_{ip}\omega_*, t_{p*}]\\
		&\qquad + [t_{k0}\cdot t_{i*}, t_{p*}] + \delta_{ip} [t_{k0}\cdot \omega_*, t_{p*}] + [t_{i*}, t_{k0}\cdot t_{p*}] + \delta_{ip} [\omega_*, t_{k0}\cdot t_{p*}]\\
		&= [[t_{i*}, \delta_{ii}t_{k*} + \delta_{ik}t_{i*}], t_{p*}] + [t_{i*}, [t_{p*}, \delta_{ip}t_{k*} + \delta_{kp}t_{i*}]] + \delta_{ip} [\omega_*, [t_{p*}, \delta_{ip}t_{k*} + \delta_{kp}t_{i*}]]\\
		&\qquad + [[t_{k*} + \delta_{ik}\omega_*, t_{i*}], t_{p*}] + [t_{i*}, [t_{k*} + \delta_{kp}\omega_*, t_{p*}]] + \delta_{ip} [\omega_*, [t_{k*} + \delta_{kp}\omega_*, t_{p*}]]\\
		&= [t_{i*}, [t_{p*}, (\delta_{ip} - 1) t_{k*} + \delta_{kp}(t_{i*} - \omega_*)]],
	\end{align*}
	so these two are equal and therefore $t_{i0}\cdot ((t_{ik} + t_{k0})\cdot t_{p*}) - (t_{ik} + t_{k0})\cdot (t_{i0}\cdot t_{p*}) = 0$. Next, we have
	\begin{align*}
		&t_{i0}\cdot ((t_{ik} + t_{k0})\cdot x_*^a) - (t_{ik} + t_{k0})\cdot (t_{i0}\cdot x_*^a)\\
		&= t_{i0}\cdot [t_{k*}, x_*^a] - (t_{ik} + t_{k0})\cdot [t_{i*}, x_*^a]\\
		&= [[t_{i*} + \delta_{ik}\omega_*, t_{k*}], x_*^a] + [t_{k*}, [t_{i*}, x_*^a]]\\
		&\qquad - ([[t_{i*}, \delta_{ii}t_{k*} + \delta_{ik}t_{k*}], x_*^a] + 0 + [[t_{k*} + \delta_{ik}\omega_*, t_{i*}], x_*^a] + [t_{i*}, [t_{k*}, x_*^a]])\\
		&= [[t_{i*}, t_{k*}], x_*^a] + [t_{k*}, [t_{i*}, x_*^a]] - ([[t_{i*}, t_{k*}], x_*^a] + [[t_{k*}, t_{i*}], x_*^a] + [t_{i*}, [t_{k*}, x_*^a]])\\
		&= 0.
	\end{align*}
	It is similar for $y_*^a$, so we obtain $t_{i0}\cdot ((t_{ik} + t_{k0})\cdot z) - (t_{ik} + t_{k0})\cdot (t_{i0}\cdot z) = 0$ for any $z\in \mathfrak{u}^f_{g,12\cdots n}$.
	\item \underline{$[t_{ij}, t_{i0} + t_{j0}] = 0$}: We have
	\begin{align*}
		&t_{ij}\cdot ((t_{i0} + t_{j0}) \cdot t_{p*})\\
		&= t_{ij}\cdot ([t_{i*} + \delta_{ip}\omega_*, t_{p*}] + [t_{j*} + \delta_{jp}\omega_*, t_{p*}])\\
		&\qquad + [t_{ij}\cdot t_{j*}, t_{p*}] + \delta_{jp} [t_{ij}\cdot\omega_*, t_{p*}] + [t_{j*}, t_{ij}\cdot t_{p*}] + \delta_{jp} [\omega_*, t_{ij}\cdot t_{p*}]\\
		&= [[t_{i*}, \delta_{ii}t_{j*} + \delta_{ij}t_{i*}], t_{p*}] + [t_{i*}, [t_{p*}, \delta_{ip}t_{j*} + \delta_{jp}t_{i*}]] + \delta_{ip} [\omega_*, [t_{p*}, \delta_{ip}t_{j*} + \delta_{jp}t_{i*}]]\\
		&\qquad + [[t_{j*}, \delta_{ij}t_{j*} + \delta_{jj}t_{i*}], t_{p*}] + [t_{j*}, [t_{p*}, \delta_{ip}t_{j*} + \delta_{jp}t_{i*}]] + \delta_{jp} [\omega_*, [t_{p*}, \delta_{ip}t_{j*} + \delta_{jp}t_{i*}]]\\
		&= \delta_{ip}[t_{i*}, [t_{p*}, t_{j*}]]  + \delta_{jp}[t_{i*}, [t_{p*},t_{i*}]] + \delta_{ip} [t_{j*}, [t_{p*}, t_{j*}]] + \delta_{jp} [t_{j*}, [t_{p*}, t_{i*}]]\\
		&\qquad + [\omega_*, [t_{p*}, \delta_{ip}t_{j*} + \delta_{ip} \delta_{jp}t_{i*} + \delta_{ip}\delta_{jp} t_{j*} + \delta_{jp}t_{i*}]].
	\end{align*}
	On the other hand, we have
	\begin{align*}
		&(t_{i0} + t_{j0}) \cdot (t_{ij}\cdot t_{p*})\\
		&= (t_{i0} + t_{j0}) \cdot [t_{p*}, \delta_{ip}t_{j*} + \delta_{jp}t_{i*}]\\
		&\qquad + [t_{j0}\cdot t_{p*}, \delta_{ip}t_{j*} + \delta_{jp}t_{i*}] + \delta_{ip}[t_{p*}, t_{j0}\cdot t_{j*}] + \delta_{jp} [t_{p*}, t_{j0}\cdot t_{i*}]\\
		&= [[t_{i*} + \delta_{ip}\omega_*, t_{p*}], \delta_{ip}t_{j*} + \delta_{jp}t_{i*}] + \delta_{ip}[t_{p*}, [t_{i*} + \delta_{ij}\omega_*, t_{j*}]] + \delta_{jp} [t_{p*}, [t_{i*} + \delta_{ii}\omega_*, t_{i*}]]\\
		&\qquad + [[t_{j*} + \delta_{jp}\omega_*, t_{p*}], \delta_{ip}t_{j*} + \delta_{jp}t_{i*}] + \delta_{ip}[t_{p*}, [t_{j*} + \delta_{jj}\omega_*, t_{j*}]] + \delta_{jp} [t_{p*},[t_{j*} + \delta_{ij}\omega_*, t_{i*}]]\\
		&= \delta_{ip} [[t_{i*}, t_{p*}], t_{j*}] + \delta_{ip} [[t_{j*}, t_{p*}], t_{j*}] + \delta_{jp} [[t_{j*}, t_{p*}], t_{i*}]\\
		&\qquad + \delta_{jp} [[t_{i*}, t_{p*}],t_{i*}] + \delta_{ip}[t_{p*}, [t_{i*}, t_{j*}]] + \delta_{jp} [t_{p*},[t_{j*}, t_{i*}]] \\
		&\qquad + \delta_{ip}(1 + \delta_{jp}) [[\omega_*, t_{p*}], t_{j*}] + \delta_{jp} (1 + \delta_{ij})[t_{p*}, [\omega_*, t_{i*}]]\\
		&\qquad + \delta_{jp}(1 + \delta_{ip}) [[\omega_*, t_{p*}], t_{i*}] + \delta_{ip}(1 + \delta_{ij})[t_{p*}, [\omega_*, t_{j*}]] \\
		&= \delta_{ip} [[t_{j*}, t_{p*}], t_{j*}] + \delta_{jp} [[t_{i*}, t_{p*}],t_{i*}] - \delta_{ip}[t_{i*}, [t_{j*}, t_{p*}]] - \delta_{jp} [t_{j*},[t_{i*}, t_{p*}]] \\
		&\qquad - \delta_{jp} (1 + \delta_{ip}) [\omega_*, [t_{i*}, t_{p*}]] - \delta_{ip}(1 + \delta_{jp}) [\omega_*, [t_{j*}, t_{p*}]] 
	\end{align*}
	so these two are equal and therefore $t_{ij}\cdot ((t_{i0} + t_{j0}) \cdot t_{p*}) - (t_{i0} + t_{j0}) \cdot (t_{ij}\cdot t_{p*}) = 0$. Next, we have
	\begin{align*}
		&t_{ij}\cdot ((t_{i0} + t_{j0}) \cdot x_*^a) - (t_{i0} + t_{j0}) \cdot (t_{ij}\cdot x_*^a)\\
		&= t_{ij}\cdot [t_{i*} + t_{j*}, x_*^a] - 0\\
		&= [t_{ij}\cdot t_{i*} + t_{ij}\cdot t_{j*}, x_*^a] + [t_{i*} + t_{j*}, t_{ij}\cdot x_*^a]\\
		&= 0 + 0.
	\end{align*}
	It is similar for $y_*^a$, so we obtain $t_{ij}\cdot ((t_{i0} + t_{j0}) \cdot z) - (t_{i0} + t_{j0}) \cdot (t_{ij}\cdot z) = 0$ for any $z\in \mathfrak{u}^f_{g,12\cdots n}$.
	\item \underline{$[t_{00}, 2t_{k0}] = 0$}: This is analogous to above since we have $t_{k0}\cdot \omega_* = 0$.
	\item \underline{$[x_i^a, y_j^b] = \delta_{ab}t_{ij}$ ($i\neq j$)}: Since $\delta_{ij} = 0$, we have
	\begin{align*}
		&x_i^a \cdot (y_j^b\cdot t_{p*}) - y_j^b \cdot (x_i^a \cdot t_{p*}) - \delta_{ab} t_{ij}\cdot t_{p*}\\
		&= x_i^a \cdot \delta_{jp}[t_{p*}, y_*^b] - y_j^b \cdot \delta_{ip}[t_{p*}, x_*^a] - \delta_{ab}[t_{p*},\delta_{jp}t_{i*} + \delta_{ip}t_{j*}]\\
		&= \delta_{jp}[\delta_{ip}[t_{p*}, x_*^a], y_*^b] + \delta_{jp}[t_{p*}, \delta_{ab}t_{i*}] - \delta_{ip}[\delta_{jp}[t_{p*}, y_*^b], x_*^a] - \delta_{ip}[t_{p*}, -\delta_{ab}t_{j*}] - \delta_{ab}[t_{p*},\delta_{jp}t_{i*} + \delta_{ip}t_{j*}]\\
		&= 0,\\
		&x_i^a \cdot (y_j^b\cdot x_*^c) - y_j^b \cdot (x_i^a \cdot x_*^c) -  \delta_{ab} t_{ij}\cdot x_*^c\\
		&= x_i^a \cdot (-\delta_{bc} t_{j*}) - 0 - 0\\
		&= -\delta_{bc} \delta_{ij}[t_{j*}, x_*^a] = 0,\mbox{ and}\\
		&x_i^a \cdot (y_j^b\cdot y_*^c) - y_j^b \cdot (x_i^a \cdot y_*^c) - \delta_{ab}t_{ij}\cdot y_*^c\\
		&= 0 - y_j^b \cdot \delta_{ac} t_{i*} - 0\\
		&= -\delta_{ac} \delta_{ij}[t_{i*}, y_*^b] = 0.
	\end{align*}
	This shows $x_i^a \cdot (y_j^b\cdot z) - y_j^b \cdot (x_i^a \cdot z) -  \delta_{ab} t_{ij}\cdot z = 0$ for any $z\in \mathfrak{u}^f_{g,12\cdots n}$.
	\item \underline{$[x_i^a, y_0^b] = \delta_{ab}t_{i0}$}: We have
	\begin{align*}
		&x_i^a\cdot (y_0^b \cdot t_{p*}) - y_0^b\cdot (x_i^a\cdot t_{p*}) - \delta_{ab} t_{i0}\cdot t_{p*}\\
		&= x_i^a\cdot [y_*^b, t_{p*}] - y_0^b\cdot \delta_{ip}[t_{p*}, x_*^a] - \delta_{ab} [t_{i*} + \delta_{ip}\omega_*, t_{p*}]\\
		&= [\delta_{ab}t_{i*}, t_{p*}] + [y_*^b, \delta_{ip}[t_{p*}, x_*^a]] - \delta_{ip} [[y_*^b, t_{p*}], x_*^a] - \delta_{ip}[t_{p*}, [y_*^b, x_*^a] + \delta_{ab}\omega_*] - \delta_{ab} [t_{i*} + \delta_{ip}\omega_*, t_{p*}]\\
		&= 0,\\
		&x_i^a\cdot (y_0^b \cdot x_*^c) - y_0^b\cdot (x_i^a\cdot x_*^c) - \delta_{ab} t_{i0}\cdot x_*^c\\
		&= x_i^a\cdot ([y_*^b, x_*^c] + \delta_{bc}\omega_*) - 0 - \delta_{ab}[t_{i*}, x^c_*]\\
		&= [\delta_{ab} t_{i*}, x_*^c] + 0 + 0 - \delta_{ab}[t_{i*}, x^c_*]\\
		&= 0,\mbox{ and}\\
		&x_i^a\cdot (y_0^b \cdot y_*^c) - y_0^b\cdot (x_i^a\cdot y_*^c) - \delta_{ab} t_{i0}\cdot y_*^c\\
		&= x_i^a\cdot [y_*^b, y_*^c] - y_0^b\cdot \delta_{ac}t_{i*} - \delta_{ab} [t_{i*}, y^c_*]\\
		&= [\delta_{ab} t_{i*}, y_*^c] + [y_*^b, \delta_{ac} t_{i*}]  - \delta_{ac} [y_*^b, t_{i*}] - \delta_{ab} [t_{i*}, y^c_*]\\
		&= 0.
	\end{align*}
	This shows $x_i^a\cdot (y_0^b \cdot z) - y_0^b\cdot (x_i^a\cdot z) - \delta_{ab} t_{i0}\cdot z = 0$ for any $z\in \mathfrak{u}^f_{g,12\cdots n}$.
	\item \underline{$[x_0^a, y_j^b] = \delta_{ab}t_{j0}$}: We have
	\begin{align*}
		&x_0^a\cdot (y_j^b\cdot t_{p*}) - y_j^b\cdot (x_0^a\cdot t_{p*}) - \delta_{ab}t_{j0}\cdot t_{p*}\\
		&= x_0^a\cdot \delta_{jp}[t_{p*}, y_*^b] - y_j^b\cdot [x_*^a, t_{p*}] - \delta_{ab}[t_{j*} + \delta_{jp}\omega_*, t_{p*}]\\
		&= \delta_{jp}[[x_*^a, t_{p*}], y_*^b] + \delta_{jp}[t_{p*},  [x_*^a, y_*^b] - \delta_{ab}\omega_*] - [-\delta_{ab}t_{j*}, t_{p*}] - [x_*^a, \delta_{jp}[t_{p*}, y_*^b]] - \delta_{ab}[t_{j*} + \delta_{jp}\omega_*, t_{p*}]\\
		&= 0,\mbox{ and}\\
		&x_0^a\cdot (y_j^b\cdot x_*^c) - y_j^b\cdot (x_0^a\cdot x_*^c) - \delta_{ab}t_{j0}\cdot x_*^c\\
		&= x_0^a\cdot (-\delta_{bc} t_{j*}) - y_j^b\cdot [x_*^a, x_*^c] - \delta_{ab}[t_{j*}, x^c_*]\\
		&= -\delta_{bc} [x_*^a, t_{j*}] - [-\delta_{ab} t_{j*}, x_*^c] - [x_*^a, -\delta_{bc}t_{j*}] - \delta_{ab}[t_{j*}, x^c_*]\\
		&= 0.
	\end{align*}
	This shows $x_0^a\cdot (y_j^b\cdot z) - y_j^b\cdot (x_0^a\cdot z) - \delta_{ab}t_{j0}\cdot z = 0$ for any $z\in \mathfrak{u}^f_{g,12\cdots n}$.
	\item \underline{$[x_i^a, x_j^b] = 0$ ($i\neq j$)}: We have
	\begin{align*}
		x_i^a\cdot (x_j^b\cdot t_{p*}) - x_j^b\cdot (x_i^a\cdot t_{p*}) &= x_i^a\cdot \delta_{jp}[t_{p*}, x_*^b] - x_j^b\cdot \delta_{ip}[t_{p*}, x_*^a]\\
		&=\delta_{jp}[\delta_{ip}[t_{p*}, x_*^a], x_*^b] + 0 - \delta_{ip}[\delta_{jp}[t_{p*}, x_*^b], x_*^a] - 0\\
		&= 0
	\end{align*}
	since $\delta_{ip}\delta_{jp} = 0$. Next, we have
	\begin{align*}
		x_i^a\cdot (x_j^b\cdot x_*^c) - x_j^b\cdot (x_i^a\cdot x_*^c) &= 0 \mbox{ and}\\
		x_i^a\cdot (x_j^b\cdot y_*^c) - x_j^b\cdot (x_i^a\cdot y_*^c) &= x_i^a\cdot \delta_{bc} t_{j*} - x_j^b\cdot \delta_{ac} t_{i*}\\
		&= \delta_{bc} \delta_{ij}[t_{j*}, x_*^a] - \delta_{ac} \delta_{ij}[t_{i*}, x_*^b] = 0.
	\end{align*}
	This shows $x_i^a\cdot (x_j^b\cdot z) - x_j^b\cdot (x_i^a\cdot z) = 0$ for any $z\in \mathfrak{u}^f_{g,12\cdots n}$.
	\item \underline{$[x_i^a, x_0^b] = 0$}: We have
	\begin{align*}
		x_i^a\cdot (x_0^b\cdot t_{p*}) - x_0^b\cdot (x_i^a\cdot t_{p*}) &= x_i^a\cdot [x_*^b, t_{p*}] - x_0^b\cdot \delta_{ip}[t_{p*}, x_*^a]\\
		&= [x_i^a\cdot x_*^b, t_{p*}] + [x_*^b, x_i^a\cdot t_{p*}] - \delta_{ip}[x_0^b\cdot t_{p*}, x_*^a] - \delta_{ip}[t_{p*}, x_0^b\cdot x_*^a]\\
		&= 0 + [x_*^b, \delta_{ip}[t_{p*}, x_*^a]] - \delta_{ip}[[x_*^b, t_{p*}], x_*^a] - \delta_{ip}[t_{p*}, [x_*^b, x_*^a]]\\
		&= 0,\\
		x_i^a\cdot (x_0^b\cdot x_*^c) - x_0^b\cdot (x_i^a\cdot x_*^c) &= x_i^a\cdot [x_*^b, x_*^c] - 0 = 0,\mbox{ and}\\
		x_i^a\cdot (x_0^b\cdot y_*^c) - x_0^b\cdot (x_i^a\cdot y_*^c) &= x_i^a\cdot ([x_*^b, y_*^c] - \delta_{bc}\omega_*) - x_0^b\cdot \delta_{ac}t_{i*}\\
		&= [x_*^b, \delta_{ac}t_{i*}] - 0 - \delta_{ac}[x_*^b, t_{i*}] = 0.
	\end{align*}
	This shows $x_i^a\cdot (x_0^b\cdot z) - x_0^b\cdot (x_i^a\cdot z) = 0$ for any $z\in \mathfrak{u}^f_{g,12\cdots n}$.
	\item \underline{$[y_i^a, y_j^b] = 0$ ($i\neq j$)}: This is analogous to the above.
	\item \underline{$[y_i^a, y_0^b] = 0$}: This is alslo analogous to the above.
	\item \underline{$[x_k^a, t_{ij}] = 0$ ($\{i,j\}\cap\{k\}=\varnothing$)}: Since $\delta_{ip}\delta_{kp} = \delta_{jp}\delta_{kp} = \delta_{ik} = \delta_{jk} = 0$, we have
	\begin{align*}
		x_k^a \cdot (t_{ij} \cdot t_{p*}) &= x_k^a \cdot [t_{p*}, \delta_{ip}t_{j*} + \delta_{jp}t_{i*}]\\
		&= [\delta_{kp}[t_{p*}, x_*^a], \delta_{ip}t_{j*} + \delta_{jp}t_{i*}] + [t_{p*}, \delta_{ip}\delta_{jk}[t_{j*}, x_*^a]+ \delta_{jp}\delta_{ik}[t_{i*}, x_*^a]]\\
		&= 0,\mbox{ and}\\
		t_{ij} \cdot (x_k^a \cdot t_{p*}) &= t_{ij} \cdot \delta_{kp}[t_{p*}, x_*^a]\\
		&=  \delta_{kp}[[t_{p*}, \delta_{ip}t_{j*} + \delta_{jp}t_{i*}], x_*^a] + 0 = 0.
	\end{align*}
	Next, we have
	\begin{align*}
		x_k^a \cdot (t_{ij} \cdot x_*^c) - t_{ij}\cdot (x_k^a \cdot x_*^c) &= 0 \mbox{ and}\\
		x_k^a \cdot (t_{ij} \cdot y_*^c) - t_{ij}\cdot (x_k^a \cdot y_*^c) &= 0 - t_{ij}\cdot \delta_{ac} t_{k*}\\
		&= [t_{k*}, \delta_{ik}t_{j*} + \delta_{jk}t_{i*}] = 0.
	\end{align*}
	This shows $x_k^a \cdot (t_{ij} \cdot z) - t_{ij}\cdot (x_k^a \cdot z) = 0$ for any $z\in \mathfrak{u}^f_{g,12\cdots n}$.
	\item \underline{$[x_0^a, t_{ij}] = 0$}: We have
	\begin{align*}
		&x_0^a \cdot (t_{ij} \cdot t_{p*}) - t_{ij} \cdot (x_0^a\cdot t_{p*})\\
		&= x_0^a \cdot [t_{p*}, \delta_{ip}t_{j*} + \delta_{jp}t_{i*}] - t_{ij} \cdot [x_*^a, t_{p*}]\\
		&= [[x_*^a, t_{p*}], \delta_{ip}t_{j*} + \delta_{jp}t_{i*}] + [t_{p*}, \delta_{ip}[x_*^a, t_{j*}] + \delta_{jp}[x_*^a, t_{i*}]] - 0 - [x_*^a, [t_{p*}, \delta_{ip}t_{j*} + \delta_{jp}t_{i*}]]\\
		&= \delta_{ip} ([[x_*^a, t_{p*}], t_{j*}] + [t_{p*},[x_*^a, t_{j*}]] - [x_*^a, [t_{p*}, t_{j*}]] ) + \delta_{jp}([[x_*^a, t_{p*}],t_{i*}] + [t_{p*}, [x_*^a, t_{i*}]]  - [x_*^a, [t_{p*}, t_{i*}]])\\
		&= 0.
	\end{align*}
	Next, we have
	\begin{align*}
		x_0^a \cdot (t_{ij} \cdot x_*^c) - t_{ij} \cdot (x_0^a\cdot x_*^c) &= 0 - t_{ij}\cdot [x_*^a, x_*^c] = 0\mbox{ and}\\
		x_0^a \cdot (t_{ij} \cdot y_*^c) - t_{ij} \cdot (x_0^a\cdot y_*^c) &= 0 - t_{ij}\cdot ([x_*^a, y_*^c] - \delta_{ac}\omega_*) = 0.
	\end{align*}
	This shows $x_0^a \cdot (t_{ij} \cdot z) - t_{ij} \cdot (x_0^a\cdot z) = 0$ for any $z\in \mathfrak{u}^f_{g,12\cdots n}$.
	\item \underline{$[x_k^a, t_{i0}] = 0$ ($i\neq k$)}: We have
	\begin{align*}
		&x_k^a \cdot (t_{i0} \cdot t_{p*}) - t_{i0} \cdot (x_k^a \cdot t_{p*})\\
		&= x_k^a \cdot [t_{i*} + \delta_{ip}\omega_*, t_{p*}] - t_{i0} \cdot \delta_{kp}[t_{p*}, x_*^a]\\
		&= [\delta_{ik}[t_{i*}, x_*^a], t_{p*}] + [t_{i*} + \delta_{ip}\omega_*, \delta_{kp}[t_{p*}, x_*^a]] - \delta_{kp}[[t_{i*} + \delta_{ip}\omega_*, t_{p*}], x_*^a] - \delta_{kp}[t_{p*}, [t_{i*}, x^a_*]]\\
		&= \delta_{kp} [t_{i*}, [t_{p*}, x_*^a]] - \delta_{kp}[[t_{i*}, t_{p*}], x_*^a] - \delta_{kp}[t_{p*}, [t_{i*}, x^a_*]]\\
		&= 0,\\
		&x_k^a \cdot (t_{i0} \cdot x_*^c) - t_{i0} \cdot (x_k^a \cdot x_*^c)\\
		&= x_k^a \cdot [t_{i*}, x^c_*] - 0\\
		&= [\delta_{ik}[t_{i*}, x_*^a], x^c_*] + 0 = 0,\mbox{ and}\\
		&x_k^a \cdot [t_{i*}, y^c_*] - t_{i0} \cdot \delta_{ac}t_{k*}\\
		&= [\delta_{ik}[t_{i*}, x_*^a], y^c_*] + [t_{i*}, \delta_{ac} t_{k*}] - \delta_{ac}[t_{i*} + \delta_{ik}\omega_*, t_{k*}]\\
		&= [t_{i*}, \delta_{ac} t_{k*}] - \delta_{ac}[t_{i*}, t_{k*}] = 0.
	\end{align*}
	This shows $x_k^a \cdot (t_{i0} \cdot z) - t_{i0} \cdot (x_k^a \cdot z) = 0$ for any $z\in \mathfrak{u}^f_{g,12\cdots n}$.
	\item \underline{$[x_k^a, t_{00}] = 0$}: This is analogous to above since we have $x_k^a\cdot \omega_* = 0$.
	\item \underline{$[y_k^a, t_{ij}] = 0$ ($\{i,j\}\cap\{k\}=\varnothing$)}: This is analogous to the above.
	\item \underline{$[y_0^a, t_{ij}] = 0$}: This is analogous to the above.
	\item \underline{$[y_k^a, t_{i0}] = 0$ ($i\neq k$)}: This is analogous to the above.
	\item \underline{$[y_k^a, t_{00}] = 0$}: This is analogous to the above.
	\item \underline{$[x_i^a + x_j^a, t_{ij}] = 0$}: We have
	\begin{align*}
		&(x_i^a + x_j^a) \cdot (t_{ij} \cdot t_{p*}) - t_{ij} \cdot ((x_i^a + x_j^a) \cdot t_{p*})\\
		&= (x_i^a + x_j^a) \cdot [t_{p*}, \delta_{ip}t_{j*} + \delta_{jp}t_{i*}] - t_{ij} \cdot (\delta_{ip} + \delta_{jp}) [t_{p*}, x_*^a]\\
		&= [(x_i^a + x_j^a) \cdot t_{p*}, \delta_{ip}t_{j*} + \delta_{jp}t_{i*}] + [t_{p*}, \delta_{ip} (x_i^a + x_j^a) \cdot t_{j*} + \delta_{jp} (x_i^a + x_j^a) \cdot t_{i*}]\\
		&\qquad - (\delta_{ip} + \delta_{jp}) [t_{ij} \cdot t_{p*}, x_*^a] - (\delta_{ip} + \delta_{jp}) [t_{p*}, t_{ij} \cdot x_*^a]\\
		&= [(\delta_{ip} + \delta_{jp}) [t_{p*}, x_*^a], \delta_{ip}t_{j*} + \delta_{jp}t_{i*}] + [t_{p*}, \delta_{ip} (\delta_{ij} + \delta_{jj}) [t_{j*}, x_*^a] + \delta_{jp} (\delta_{ii} + \delta_{ij}) [t_{i*}, x_*^a]]\\
		&\qquad - (\delta_{ip} + \delta_{jp}) [[t_{p*}, \delta_{ip}t_{j*} + \delta_{jp}t_{i*}], x_*^a] - 0\\
		&= (\delta_{ip} + \delta_{jp}) [[t_{p*}, x_*^a], \delta_{ip}t_{j*} + \delta_{jp}t_{i*}] + (\delta_{ip} + \delta_{jp}) [t_{p*},  [\delta_{ip} t_{j*} + \delta_{jp} t_{i*}, x_*^a]]\\
		&\qquad - (\delta_{ip} + \delta_{jp}) [[t_{p*}, \delta_{ip}t_{j*} + \delta_{jp}t_{i*}], x_*^a]\\
		&= 0,\\
		&(x_i^a + x_j^a) \cdot (t_{ij} \cdot x_*^c) - t_{ij} \cdot ((x_i^a + x_j^a) \cdot x_*^c) = 0,\mbox{ and}\\
		&(x_i^a + x_j^a) \cdot (t_{ij} \cdot y_*^c) - t_{ij} \cdot ((x_i^a + x_j^a) \cdot y_*^c)\\
		&= 0 - t_{ij}\cdot \delta_{ac} (t_{i*} + t_{j*})\\
		&= - \delta_{ac} ([t_{i*}, \delta_{ii}t_{j*} + \delta_{ji}t_{i*}] + [t_{j*}, \delta_{ij}t_{j*} + \delta_{jj}t_{i*}])\\
		&= 0.
	\end{align*}
	This shows $(x_i^a + x_j^a) \cdot (t_{ij} \cdot z) - t_{ij} \cdot ((x_i^a + x_j^a) \cdot z) = 0$ for any $z\in \mathfrak{u}^f_{g,12\cdots n}$.
	\item \underline{$[x_i^a + x_0^a, t_{i0}] = 0$}: We have
	\begin{align*}
		&(x_i^a + x_0^a) \cdot (t_{i0} \cdot t_{p*}) - t_{i0} \cdot ((x_i^a + x_0^a) \cdot t_{p*})\\
		&= (x_i^a + x_0^a) \cdot [t_{i*} + \delta_{ip}\omega_*, t_{p*}] - t_{i0} \cdot (\delta_{ip} - 1)[t_{p*}, x_*^a]\\
		&= [x_i^a \cdot t_{i*} + x_0^a\cdot t_{i*}, t_{p*}] + [t_{i*} + \delta_{ip}\omega_*, (x_i^a + x_0^a) \cdot t_{p*}]\\
		&\qquad - (\delta_{ip} - 1)[t_{i0} \cdot t_{p*}, x_*^a] - (\delta_{ip} - 1)[t_{p*}, t_{i0} \cdot x_*^a] \\	
		&= [\delta_{ii}[t_{i*}, x_*^a] + [x_*^a, t_{i*}], t_{p*}] + [t_{i*} + \delta_{ip}\omega_*, \delta_{ip}[t_{p*}, x_*^a] + [x_*^a, t_{p*}] ]\\
		&\qquad - (\delta_{ip} - 1)[[t_{i*} + \delta_{ip}\omega_*, t_{p*}], x_*^a] - (\delta_{ip} - 1)[t_{p*}, [t_{i*}, x^a_*]] \\	
		&= (\delta_{ip} - 1) [t_{i*} + \delta_{ip}\omega_*, [t_{p*}, x_*^a]] - (\delta_{ip} - 1)[[t_{i*} + \delta_{ip}\omega_*, t_{p*}], x_*^a] - (\delta_{ip} - 1)[t_{p*}, [t_{i*}, x^a_*]] \\
		&=  (\delta_{ip} - 1) ([t_{i*}, [t_{p*}, x_*^a]] + [x_*^a, [t_{i*}, t_{p*}]] + [t_{p*}, [x^a_*, t_{i*}]]) + (\delta_{ip} - 1)\delta_{ip}([\omega_*, [t_{p*}, x_*^a]] - [[\omega_*, t_{p*}], x_*^a])\\
		&= 0
	\end{align*}
	since $ (\delta_{ip} - 1)\delta_{ip} = 0$. Next, we have
	\begin{align*}
		&(x_i^a + x_0^a) \cdot (t_{i0} \cdot x_*^c) - t_{i0} \cdot ((x_i^a + x_0^a) \cdot x_*^c)\\
		&= (x_i^a + x_0^a) \cdot [t_{i*}, x^c_*] - t_{i0} \cdot [x_*^a, x_*^c]\\
		&= [(x_i^a + x_0^a) \cdot t_{i*}, x^c_*] + [t_{i*},  (x_i^a + x_0^a) \cdot x^c_*] - [t_{i0} \cdot x_*^a, x_*^c] - [x_*^a,  t_{i0} \cdot x_*^c]\\
		&= [\delta_{ii}[t_{i*}, x_*^a] + [x_*^a, t_{i*}], x^c_*] + [t_{i*},  [x_*^a,x^c_*]] - [[t_{i*}, x^a_*], x_*^c] - [x_*^a,  [t_{i*}, x^c_*]]\\
		&= 0 \mbox{ and}\\
		&(x_i^a + x_0^a) \cdot (t_{i0} \cdot y_*^c) - t_{i0} \cdot ((x_i^a + x_0^a) \cdot y_*^c)\\
		&= (x_i^a + x_0^a) \cdot [t_{i*}, y^c_*] - t_{i0} \cdot (\delta_{ac} t_{i*} + [x_*^a, y_*^c] - \delta_{ac}\omega_*)\\
		&= [(x_i^a + x_0^a) \cdot t_{i*}, y^c_*] + [t_{i*}, (x_i^a + x_0^a) \cdot y^c_*] - (\delta_{ac}t_{i0} \cdot t_{i*} + [t_{i0} \cdot x_*^a, y_*^c] + [x_*^a, t_{i0} \cdot y_*^c] - \delta_{ac}t_{i0} \cdot \omega_*)\\
		&= [\delta_{ii}[t_{i*}, x_*^a] + [x_*^a, t_{i*}], y^c_*] + [t_{i*}, \delta_{ac} t_{i*} + [x_*^a, y_*^c] - \delta_{ac}\omega_*]\\
		&\qquad - (\delta_{ac}[t_{i*} + \delta_{ii}\omega_*, t_{i*}] + [[t_{i*}, x^a_*], y_*^c] + [x_*^a, [t_{i*}, y^c_*]])\\
		&= [t_{i*}, [x_*^a, y_*^c]] - ( [[t_{i*}, x^a_*], y_*^c] + [x_*^a, [t_{i*}, y^c_*]])\\
		&= 0.
	\end{align*}
	This shows $(x_i^a + x_0^a) \cdot (t_{i0} \cdot z) - t_{i0} \cdot ((x_i^a + x_0^a) \cdot z) = 0$ for any $z\in \mathfrak{u}^f_{g,12\cdots n}$.
	\item \underline{$[2x_0^a, t_{00}] = 0$}: This is analogous to above since we have $x_0^a\cdot \omega_* = 0$.
	\item \underline{$[y_i^a + y_j^a, t_{ij}] = 0$}: This is analogous to the above.
	\item \underline{$[y_i^a + y_0^a, t_{i0}] = 0$}: This is analogous to the above.
	\item \underline{$[2y_0^a, t_{00}] = 0$}: This is analogous to the above.
	\item \underline{$\sum_{1\leq a\leq g} [x_i^a, y_i^a] + \sum_{1 \leq j\leq n, j\neq i} t_{ij} + t_{i0} = (g-1) t_{ii}$}: We have
	\begin{align*}
		&\Big(\sum_{1\leq a\leq g} [x_i^a, y_i^a] + \sum_{1 \leq j\leq n, j\neq i} t_{ij} + t_{i0} - (g-1) t_{ii}\Big)\cdot t_{p*}\\
		&= \sum_{1\leq a\leq g} x_i^a\cdot (y_i^a\cdot t_{p*}) - y_i^a\cdot (x_i^a\cdot t_{p*}) + \sum_{1 \leq j\leq n, j\neq i} t_{ij}\cdot t_{p*} + t_{i0} \cdot t_{p*} - (g-1) t_{ii}\cdot t_{p*}\\
		&= \sum_{1\leq a\leq g} x_i^a\cdot \delta_{ip}[t_{p*}, y_*^a] - y_i^a\cdot \delta_{ip}[t_{p*}, x_*^a] + \sum_{1 \leq j\leq n, j\neq i} [t_{p*}, \delta_{ip}t_{j*} + \delta_{jp}t_{i*}]\\
		&\qquad  + [t_{i*} + \delta_{ip}\omega_*, t_{p*}] - (g-1)[t_{p*}, 2\delta_{ip}t_{i*}] \\
		&= \sum_{1\leq a\leq g} \delta_{ip}([x_i^a\cdot t_{p*}, y_*^a] + [t_{p*}, x_i^a\cdot y_*^a]) - \delta_{ip}([ y_i^a\cdot t_{p*}, x_*^a] + [t_{p*},  y_i^a\cdot x_*^a])\\
		&\qquad + \Big(\sum_{1 \leq j\leq n, j\neq i} \delta_{ip}[t_{p*}, t_{j*}]\Big) + (1 - \delta_{ip})[t_{p*}, t_{i*}] + [t_{i*} + \delta_{ip}\omega_*, t_{p*}] - 0\\
		&= \sum_{1\leq a\leq g} \delta_{ip}([\delta_{ip}[t_{p*}, x_*^a], y_*^a] + [t_{p*}, t_{i*}]) - \delta_{ip}([\delta_{ip}[t_{p*}, y_*^a], x_*^a] + [t_{p*}, -t_{i*}])\\
		&\qquad + \Big(\sum_{1 \leq j\leq n, j\neq i} \delta_{ip}[t_{p*}, t_{j*}]\Big) + (1 - \delta_{ip})[t_{p*}, t_{i*}] + [t_{i*} + \delta_{ip}\omega_*, t_{p*}]\\
		&= \delta_{ip}[t_{p*}, \omega_*] + [\delta_{ip}\omega_*, t_{p*}]\\
		&= 0,\\
		&\Big(\sum_{1\leq a\leq g} [x_i^a, y_i^a] + \sum_{1 \leq j\leq n, j\neq i} t_{ij} + t_{i0} - (g-1) t_{ii}\Big)\cdot x_*^c\\
		&= \sum_{1\leq a\leq g} x_i^a\cdot (y_i^a\cdot x_*^c) - y_i^a\cdot (x_i^a\cdot x_*^c) + \sum_{1 \leq j\leq n, j\neq i} t_{ij}\cdot x_*^c + t_{i0} \cdot x_*^c - (g-1) t_{ii}\cdot x_*^c\\
		&= \sum_{1\leq a\leq g} x_i^a\cdot (-\delta_{ac} t_{i*}) - 0 + 0 + [t_{i*}, x^c_*] - 0\\
		&= -x_i^c\cdot  t_{i*}  + [t_{i*}, x^c_*] \\
		&= 0,\mbox{ and}\\
		&\Big(\sum_{1\leq a\leq g} [x_i^a, y_i^a] + \sum_{1 \leq j\leq n, j\neq i} t_{ij} + t_{i0} - (g-1) t_{ii}\Big)\cdot y_*^c\\
		&= \sum_{1\leq a\leq g} x_i^a\cdot (y_i^a\cdot y_*^c) - y_i^a\cdot (x_i^a\cdot y_*^c) + \sum_{1 \leq j\leq n, j\neq i} t_{ij}\cdot y_*^c + t_{i0} \cdot y_*^c - (g-1) t_{ii}\cdot y_*^c\\
		&= \sum_{1\leq a\leq g} (0 - y_i^a\cdot \delta_{ac}t_{i*}) + 0 + [t_{i*}, y^c_*] - 0\\
		&= - y_i^c\cdot t_{i*}  + [t_{i*}, y^c_*] \\
		&= 0.
	\end{align*}
	This shows $\Big(\sum_{1\leq a\leq g} [x_i^a, y_i^a] + \sum_{1 \leq j\leq n, j\neq i} t_{ij} + t_{i0} - (g-1) t_{ii}\Big)\cdot z = 0$ for any $z\in \mathfrak{u}^f_{g,12\cdots n}$.
	\item \underline{$\sum_{1\leq a\leq g} [x_0^a, y_0^a] + \sum_{1 \leq j\leq n} t_{j0} = (g-1) t_{00}$}: We have
	\begin{align*}
		&\Big(\sum_{1\leq a\leq g} [x_0^a, y_0^a] + \sum_{1 \leq j\leq n} t_{j0} - (g-1) t_{00}\Big)\cdot t_{p*}\\
		&= \sum_{1\leq a\leq g} x_0^a\cdot (y_0^a\cdot t_{p*}) - y_0^a\cdot (x_0^a\cdot t_{p*}) + \sum_{1 \leq j\leq n} t_{j0}\cdot t_{p*} - (g-1) t_{00}\cdot t_{p*}\\
		&= \sum_{1\leq a\leq g} x_0^a\cdot [y_*^a, t_{p*}] - y_0^a\cdot [x_*^a, t_{p*}] + \sum_{1 \leq j\leq n} [t_{i*} + \delta_{jp}\omega_*, t_{p*}] - 2(g-1) [t_{p*},\omega_*]\\
		&= \sum_{1\leq a\leq g} [x_0^a\cdot y_*^a, t_{p*}] + [y_*^a, x_0^a\cdot t_{p*}] - [y_0^a\cdot x_*^a, t_{p*}] - [x_*^a, y_0^a\cdot t_{p*}] \\
		&\qquad + \Big(\sum_{1 \leq j\leq n} [t_{i*}, t_{p*}]\Big) + [\omega_*, t_{p*}] - 2(g-1) [t_{p*},\omega_*]\\
		&= \sum_{1\leq a\leq g} [[x_*^a, y_*^a] - \delta_{aa}\omega_*, t_{p*}] + [y_*^a, [x_*^a, t_{p*}]] - [[y_*^a, x_*^a] + \delta_{aa}\omega_*, t_{p*}] - [x_*^a, [y_*^a, t_{p*}]] \\
		&\qquad + \Big(\sum_{1 \leq j\leq n} [t_{i*}, t_{p*}]\Big) + [\omega_*, t_{p*}] - 2(g-1) [t_{p*},\omega_*]\\
		&= \sum_{1\leq a\leq g} [[x_*^a, y_*^a], t_{p*}] - [\omega_*, t_{p*}]+  [y_*^a, [x_*^a, t_{p*}]] - [[y_*^a, x_*^a], t_{p*}] - [\omega_*, t_{p*}] - [x_*^a, [y_*^a, t_{p*}]] \\
		&\qquad + \Big(\sum_{1 \leq j\leq n} [t_{i*}, t_{p*}]\Big) + [\omega_*, t_{p*}] - 2(g-1) [t_{p*},\omega_*]\\
		&= \Big(\sum_{1\leq a\leq g} [[x_*^a, y_*^a], t_{p*}]\Big) - 2g[\omega_*, t_{p*}] + \Big(\sum_{1 \leq j\leq n} [t_{i*}, t_{p*}]\Big) + [\omega_*, t_{p*}] - 2(g-1) [t_{p*},\omega_*]\\
		&= [\omega_*, t_{p*}] - 2g[\omega_*, t_{p*}] + [\omega_*, t_{p*}] + 2(g-1) [\omega_*, t_{p*}]\\
		&= 0 \mbox{ and}\\
		&\Big(\sum_{1\leq a\leq g} [x_0^a, y_0^a] + \sum_{1 \leq j\leq n} t_{j0} - (g-1) t_{00}\Big)\cdot x_*^c\\
		&= \sum_{1\leq a\leq g} x_0^a\cdot (y_0^a\cdot x_*^c) - y_0^a\cdot (x_0^a\cdot x_*^c) + \sum_{1 \leq j\leq n} t_{j0}\cdot x_*^c - (g-1) t_{00}\cdot x_*^c\\
		&= \sum_{1\leq a\leq g} x_0^a\cdot ([y_*^a, x_*^c] + \delta_{ac}\omega_*) - y_0^a\cdot [x_*^a, x_*^c] + \sum_{1 \leq j\leq n} [t_{j*}, x^c_*] - 2(g-1) [x_*^c,\omega_*]\\
		&= \sum_{1\leq a\leq g} ([x_0^a\cdot y_*^a, x_*^c] + [y_*^a, x_0^a\cdot x_*^c]) - ([y_0^a\cdot x_*^a, x_*^c] + [x_*^a, y_0^a\cdot x_*^c]) + \sum_{1 \leq j\leq n} [t_{j*}, x^c_*] - 2(g-1) [x_*^c,\omega_*]\\
		&= \sum_{1\leq a\leq g} ([[x_*^a, y_*^a] - \delta_{aa}\omega_*, x_*^c] + [y_*^a, [x_*^a, x_*^c]]) - ([[y_*^a, x_*^a] + \delta_{aa}\omega_*, x_*^c] + [x_*^a, [y_*^a, x_*^c] + \delta_{ac}\omega_*])\\
		&\qquad + \sum_{1 \leq j\leq n} [t_{j*}, x^c_*] - 2(g-1) [x_*^c,\omega_*]\\
		&= \Big(\sum_{1\leq a\leq g} ([[x_*^a, y_*^a], x_*^c - [\omega_*, x_*^c] + [y_*^a, [x_*^a, x_*^c]]) - ([[y_*^a, x_*^a], x_*^c] + [\omega_*, x_*^c] + [x_*^a, [y_*^a, x_*^c]])\Big) - [x_*^c,\omega_*]\\
		&\qquad + \sum_{1 \leq j\leq n} [t_{j*}, x^c_*] - 2(g-1) [x_*^c,\omega_*]\\
		&= \Big(\sum_{1\leq a\leq g} [[x_*^a, y_*^a], x_*^c  \Big) - 2g[\omega_*, x_*^c] - [x_*^c,\omega_*] + \sum_{1 \leq j\leq n} [t_{j*}, x^c_*] - 2(g-1) [x_*^c,\omega_*]\\
		&=  [\omega_*, x_*^c] - 2g[\omega_*, x_*^c] - [x_*^c,\omega_*] + 2(1-g) [x_*^c,\omega_*]\\
		&= 0.
	\end{align*}
	 It is similar for $y_*^c$. This shows $\Big(\sum_{1\leq a\leq g} [x_0^a, y_0^a] + \sum_{1 \leq j\leq n} t_{j0} - (g-1) t_{00}\Big)\cdot z = 0$ for any $z\in \mathfrak{u}^f_{g,12\cdots n}$.
\end{itemize}
Thus, we have exhausted all the relations, and this concludes the proof.\qed\\

\begin{lemma}
The linear map
\[
	F := \iota + \circ_{0}\id_{*0}\colon \mathfrak{u}^f_{g,12\cdots n} \rtimes \mathfrak{t}^f_{g,12\cdots n0} \to \mathfrak{t}^f_{g,12\cdots n*0}
\]
is a Lie algebra homomorphism.
\end{lemma}
\noindent Proof. First of all, $\iota$ is a Lie algebra homomorphism since $t_{**}\in \mathfrak{t}^f_{g,12\cdots n*0}$ is central. In addition, $\circ_{0}\id_{*0}$ is a Lie algebra homomorphism (which is already included in the fact that $\{\mathfrak{t}_{g,n}^f\}_{n\geq 1}$ is an operad module over the category of Lie algebras). We are done if the map $F$ preserves the $\mathfrak{t}^f_{g,12\cdots n0}$-action over $\circ_{0}\id_{*0}$, but the action is defined using the relations in $\mathfrak{t}^f_{g,12\cdots n*0}$, so this completes the proof.\qed\\

\begin{deflem}
We define the Lie algebra homomorphism $G\colon \mathfrak{t}^f_{g,12\cdots n*0} \to \mathfrak{u}^f_{g,12\cdots n} \rtimes \mathfrak{t}^f_{g,12\cdots n0}$ by
\begin{align*}
	G(t_{ij}) &= (0, t_{ij}),\quad G(t_{i*}) = (t_{i*}, 0),\quad G(t_{i0}) = (-t_{i*}, t_{i0}),\\
	G(t_{**}) &= (t_{**}, 0),\quad G(t_{*0}) = (-\omega_*, 0),\quad G(t_{00}) = (-t_{**} + 2\omega_*, t_{00}),\\
	G(x_i^a) &= (0, x_i^a),\quad G(x_*^a) = (x_*^a, 0),\quad G(x_0^a) = (-x_*^a, x_0^a),\\
	G(y_i^a) &= (0, y_i^a),\quad G(y_*^a) = (y_*^a, 0),\mbox{ and}\quad G(y_0^a) = (-y_*^a, y_0^a).
\end{align*}
for $1\leq i,j\leq n$ and $1\leq a \leq g$.
\end{deflem}
\noindent Proof. We only check the last relation in the definition of $\mathfrak{t}^f_{g,12\cdots n*0}$. Putting $I = \{1,2,\dotsc, n, * ,0\}$, we have
\begin{align*}
	&\sum_{1\leq a\leq g} [G(x_i^a), G(y_i^a)] + \sum_{p\in I\setminus\{i\}} G(t_{ip})\\
	&= \sum_{1\leq a\leq g} [(0, x_i^a), (0, y_i^a)] + \sum_{1\leq p\leq n, p\neq i} (0, t_{ip}) + (t_{i*}-t_{i*}, t_{i0}) = (0, (g-1)t_{ii}) = G((g-1) t_{ii}),\\
	&\sum_{1\leq a\leq g} [G(x_*^a), G(y_*^a)] + \sum_{p\in I\setminus\{*\}} G(t_{p*})\\
	&= \sum_{1\leq a\leq g} [(x_*^a,0), (y_*^a,0)] + \sum_{1\leq p\leq n} (t_{p*},0) + (-\omega_*, 0) =  ((g-1) t_{**}, 0) = G((g-1) t_{**})
\end{align*}
for $1\leq i\leq n$. For $i=0$, we have
\begin{align*}
	&\sum_{1\leq a\leq g} [G(x_0^a), G(y_0^a)] + \sum_{p\in I\setminus\{0\}} G(t_{p0})\\
	&= \sum_{1\leq a\leq g} [(-x_*^a, x_0^a), (-y_*^a, y_0^a)] + \sum_{1\leq p\leq n} (-t_{p*}, t_{p0}) + (-\omega_*, 0)\\
	&= \sum_{1\leq a\leq g} ([x_*^a, y_*^a] - x_0^a\cdot y_*^a + y_0^a\cdot x_*^a, [x_0^a, y_0^a]) + \sum_{1\leq p\leq n} (-t_{p*}, t_{p0}) + (-\omega_*, 0) \\
	&= \sum_{1\leq a\leq g} ([x_*^a, y_*^a] - ([x_*^a, y_*^a] - \omega_*) + ([y_*^a, x_*^a] + \omega_*), 0) + \sum_{1\leq p\leq n} (-t_{p*}, 0) + (0, (g-1)t_{00})+ (-\omega_*, 0)\\
	&= \sum_{1\leq a\leq g} ( 2\omega_* + [y_*^a, x_*^a], 0) + \sum_{1\leq p\leq n} (-t_{p*}, 0) + (0, (g-1)t_{00})+ (-\omega_*, 0)\\
	&= 2g\omega_* - (\omega_* + (g-1)t_{**}) + (0, (g-1)t_{00})+ (-\omega_*, 0)\\
	&=(2(g-1)\omega_* - (g-1) t_{**},  (g-1)t_{00})\\
	&= G((g-1)t_{00}),
\end{align*}
so the relation is respected by $G$. The rest is straightforward.\qed

\begin{lemma}
The Lie algebra homomorphisms $F$ and $G$ are isomorphisms.
\end{lemma}
\noindent Proof. This is straightforward.\qed

\begin{theorem}\label{thm:split}
The sequence
\[\begin{tikzcd}
	0 \arrow[r]& \mathfrak{u}^f_{g,12\cdots n} \arrow[r, "\iota"]& \mathfrak{t}^f_{g,12\cdots n*0} \arrow[r, "\varepsilon^*", yshift = 2pt]& \mathfrak{t}^f_{g,12\cdots n0} \arrow[r]\arrow[l, yshift = -2pt, "\circ_{0}\id_{*0}"]& 0
\end{tikzcd}\]
is split.
\end{theorem}
\noindent Proof. By the lemma above, $\iota(\mathfrak{u}^f_{g,12\cdots n})$ is an ideal of $ \mathfrak{t}^f_{g,12\cdots n*0}$. In addition, the map
\[
	 \mathfrak{t}^f_{g,12\cdots n*0}/\iota(\mathfrak{u}^f_{g,12\cdots n}) \to \mathfrak{t}^f_{g,12\cdots n0}
\]
induced by $G$ coincides with $\varepsilon^*$. Therefore, the pair $(F,G)$ corresponds exactly to the split diagram in the claim.\qed\\

Since the graded Lie algebras $\mathfrak{u}^f_{g,12\cdots n}$, $\mathfrak{t}^f_{g,12\cdots n*0}$ and $\mathfrak{t}^f_{g,12\cdots n0}$ are pro-nilpotent, taking the exponential yields the exact sequence of groups in the second row of the diagram in Lemma \ref{lem:exactness}. This concludes the proof of Lemma \ref{lem:exactness}.\\

\small
\bibliographystyle{alphaurl}
\bibliography{hgkva.bib}

\end{document}